\documentclass{gtart_a}
\pdfoutput=1

\title[SL singularities and the Lagrangian prescribed boundary problem]{Obstructions to special Lagrangian desingularizations\\ and
the Lagrangian prescribed boundary problem}

\author{Mark Haskins} 
\givenname{Mark}
\surname{Haskins}
\address{Department of Mathematics\\
Imperial College London\\\newline
South Kensington\\
London SW7 2AZ\\England}
\email{m.haskins@imperial.ac.uk}
\urladdr{}

\author{Tommaso Pacini}
\givenname{Tommaso}
\surname{Pacini}
\address{Georgia Institute of Technology\\
Department of Mathematics\\\newline
686 Cherry Street\\
Atlanta, GA 30332\\USA}
\email{t.pacini@imperial.ac.uk}
\urladdr{}

\volumenumber{10}
\issuenumber{}
\publicationyear{2006}
\papernumber{37}
\lognumber{0768}
\startpage{1453}
\endpage{1521}

\doi{}
\MR{}
\Zbl{}

\arxivreference{math.DG/0609352}


\received{2 July 2006}
\revised{}
\accepted{6 September 2006}
\published{21 October 2006}
\publishedonline{21 October 2006}
\proposed{Frances Kirwan}
\seconded{Jim Bryan, Simon Donaldson}
\corresponding{}
\editor{Yasha Eliashberg, Frances Kirwan, Gang Tian}
\version{}

\AtBeginDocument{\let\bar\wbar\let\tilde\wtilde\let\hat\what}
\def\S{Section }

\theoremstyle{definition}
\newtheorem{assump}{Assumptions about $\Sigma$}

\makeautorefname{assump}{Assumptions A--C}

\def\hatphi{\mskip2mu\hat{\mskip-2mu\vphantom{t}\smash{\phi}\mskip-.5mu}\mskip.5mu}

\makeatletter
\def\cnewtheorem#1[#2]#3{\newtheorem{#1}{#3}[section]
\expandafter\let\csname c@#1\endcsname\c@theorem}

\theoremstyle{plain}
\newtheorem{theorem}{Theorem}[section]
\cnewtheorem{lemma}[theorem]{Lemma}
\cnewtheorem{prop}[theorem]{Proposition}
\cnewtheorem{corollary}[theorem]{Corollary}
\cnewtheorem{conjecture}[theorem]{Conjecture}
\theoremstyle{definition}
\cnewtheorem{definition}[theorem]{Definition}
\cnewtheorem{example}[theorem]{Example}
\cnewtheorem{xca}[theorem]{Exercise}
\theoremstyle{remark}
\cnewtheorem{remark}[theorem]{Remark}
\numberwithin{equation}{section}

\makeatother  

\newcommand{\acslg}{ACSL\ }

\newcommand{\ka}{K\"ahler\ }

\newcommand{\gl}[1]{\textrm{GL({#1},\,$\C$)}}
\newcommand{\glr}[1]{\textrm{GL({#1},\,$\R$)}}
\newcommand{\unitary}[1]{\textrm{U({#1})}}
\newcommand{\sunitary}[1]{\textrm{SU({#1})}}
\newcommand{\sorth}[1]{\textrm{SO({#1})}}
\newcommand{\orth}[1]{\textrm{O({#1})}}
\newcommand{\slg}{special Lagrangian\ }

\newcommand{\vol}{\operatorname{vol}}

\newcommand{\RP}{\mathbb{RP}}
\newcommand{\CP}{\mathbb{CP}}
\newcommand{\Sph}{\mathbb{S}}

\newcommand{\ra}{\rightarrow}
\newcommand{\Hom}{\operatorname{Hom}}
\newcommand{\Ext}{\operatorname{Ext}}
\newcommand{\imag}{\operatorname{Im}}
\newcommand{\real}{\operatorname{Re}}
\newcommand{\diag}{\operatorname{diag}}
\newcommand{\Imag}{\operatorname{Im}}
\newcommand{\Real}{\operatorname{Re}}
\newcommand{\trace}{\operatorname{Tr}}
\newcommand{\dist}{\operatorname{dist}}

\begin{document}

\begin{asciiabstract}
We exhibit infinitely many, explicit special Lagrangian isolated
singularities that admit no asymptotically conical special Lagrangian
smoothings. The existence/ nonexistence of such smoothings is an
important component of the current efforts to understand which
singular special Lagrangians arise as limits of smooth special
Lagrangians.

We also use soft methods from symplectic geometry (the relative
version of the h--principle) and tools from algebraic topology to prove
(both positive and negative) results about Lagrangian
desingularizations of Lagrangian submanifolds with isolated
singularities; we view the Lagrangian desingularization problem as the
natural soft analogue of the special Lagrangian smoothing problem.
\end{asciiabstract}

\begin{abstract}
We exhibit infinitely many, explicit special Lagrangian isolated
singularities that admit no asymptotically conical special Lagrangian smoothings.
The existence/ nonexistence of such smoothings of special Lagrangian cones
is an important component of the current efforts to understand which singular special Lagrangians
arise as limits of smooth special Lagrangians.

We also use soft methods from symplectic geometry (the relative version
of the $h$--principle for Lagrangian immersions) and tools from
algebraic topology to prove (both positive and negative) results about
Lagrangian desingularizations of Lagran\-gian submanifolds with isolated singularities;
we view the (Maslov-zero) Lagrangian desingularization problem as the natural soft analogue of the special Lagrangian smoothing problem.
\end{abstract}

\maketitle

\section{Introduction}\label{S:intro}
Let $M$ be a Calabi--Yau manifold of complex dimension $n$ with
K\"ahler form $\omega$ and nonzero parallel holomorphic $n$--form
$\Omega$. Suitably normalized, $\real{\Omega}$ is a calibrated
form whose calibrated submanifolds are called special Lagrangian
(SL) submanifolds; see Harvey and Lawson \cite{harvey:lawson}. SL submanifolds are thus a
very natural class of volume-minimizing submanifolds in Calabi--Yau
manifolds. SL submanifolds also appear in string theory as ``supersymmetric
cycles'' and play a fundamental role in the Strominger--Yau--Zaslow
approach to mirror symmetry \cite{syz}. In SYZ-related work and
also in other problems in special Lagrangian geometry, singular SL
objects play a fundamental role. This has motivated a considerable
amount of recent work devoted to singular SL submanifolds
\cite{mcintosh:carberry,gross:fibrations1,gross:fibrations2,gross:examples,haskins:invent,haskins:ajm,haskins:kapouleas,joyce:syz,joyce:conifolds5,joyce:conifolds:summary,joyce:conifolds1,joyce:conifolds2,joyce:conifolds3,joyce:conifolds4,mcintosh}.

So far, the most fruitful approach to understanding singular SL objects has been
to focus on special types of singularities.
For example, one might study ruled SL
$3$--folds, $\textrm{U}(1)$--invariant SL $3$--folds, or SL $3$--folds with isolated
conical singularities.
In particular, a \textit{special Lagrangian $n$--fold with isolated
conical singularities} is a singular SL $n$--fold with an isolated
set of singular points each of which is modelled on a SL cone in
$\C^n$. We call such a singular SL $n$--fold a \textit{special
Lagrangian conifold}. In a series of five papers, Joyce developed
the basic foundations of a theory of compact SL conifolds in (almost) Calabi--Yau manifolds
\cite{joyce:conifolds5}.

In particular, Joyce studied the desingularization theory of SL conifolds. The desingularization theory
has both local and global aspects; since the global 
aspects are described in detail in \cite{joyce:conifolds5}, it remains to understand the local smoothing question. More concretely,
given a SL cone $C$ in $\C^n$ we would like to understand all asymptotically conical SL (ACSL) $n$--folds in
$\C^n$ which at infinity approach the given cone $C$. We call this the \textit{ACSL smoothing problem\/}.

The main aim of this paper is to study the ACSL smoothing problem and various
closely related questions. To begin with we prove a number of nonsmoothing results.
For example, we exhibit an explicit SL cone $C$ in $\C^6$
whose link is a $5$--manifold which is not nullcobordant.
Hence $C$ is an isolated SL singularity which for topological reasons cannot admit any smoothing whatsoever.

We then show that there are \textit{infinitely many\/} explicit
SL cones which admit no ACSL smoothings (with rate $\lambda\le 0$ -- see equation \eqref{E:acsl} for the precise definition). On the other hand, infinitely many of these SL cones have links which are nullcobordant. Therefore there is no differential-topological
obstruction to smoothing these singularities; in fact, our obstruction to ACSL smoothings arises from analytic constraints implicit in the SL condition.

Given such an SL cone $C$, it is natural to ask whether $C$ admits Lagrangian
rather than special Lagrangian smoothings.
In fact, since any SL submanifold has zero Maslov class, the natural
weakening of the SL smoothing question to the Lagrangian category is to ask whether $C$ admits Maslov-zero
Lagrangian smoothings. Therefore the main focus of the paper
is on the Lagrangian-topological aspects of these smoothing questions, viewed
as a preliminary step toward a better understanding of the existence/nonexistence
of ACSL smoothings and SL desingularizations of SL conifolds.

More specifically, in this paper we address questions
such as the following:
\begin{enumerate}
\item Let $C$ be a SL cone in $\C^n$. Does there exist a smooth complete oriented immersed Lagrangian submanifold
$Y$ which coincides with $C$ outside the sphere $\Sph^{2n-1}$?
What can we say about the topology of $Y$? Can we find such a $Y$ with zero Maslov class?
\item Let $X$ be a SL conifold in an (almost) Calabi--Yau manifold. Does there exist a Maslov-zero Lagrangian submanifold $Y$ which coincides with $X$ away from a neighbourhood of its singularities?
\end{enumerate}
We think of these two questions as soft analogues of the SL
desingularization problem. Notice, in particular, that in questions (1) and (2) the smoothings coincide with the original singular object
outside a compact set. Real analyticity would prevent this from
happening in the SL case: for example, a complete ACSL smoothing of a SL $n$--cone $C$ will never coincide with $C$ in an open set
(or even intersect it along a $(n-1)$--dimensional submanifold). We can take this behaviour into account by prescribing other types of asymptotics on the ``ends'' of $Y$. Hence one could also ask the following as a soft analogue of the existence question for ACSL smoothings.
\begin{itemize}
\item[(3)] Let $C$ be a SL cone in $\C^n$ with link $\Sigma$. Suppose we
are given a Maslov-zero Lagrangian conical end, ie a Maslov-zero
Lagrangian immersion $f$ of $\Sigma\times [1,\infty)$ which is asymptotic to $C$ at infinity. What Maslov-zero Lagrangian ``fillings'' of the conical
end $f$ exist?
\end{itemize}

All these Lagrangian smoothing problems can be viewed as particular
cases of a certain boundary value problem for Lagrangian submanifolds
 which we call the \textit{Prescribed Boundary Problem\/}
and study in \fullref{s:bdy}. The Prescribed Boundary Problem
is related to but different from the notion of Lagrangian cobordism
groups introduced by Arnold \cite{arnold} and studied in greater depth by
Audin \cite{audin:cobord,audin:book}.

At least in the category of \textit{immersed\/}
Lagrangian submanifolds it turns out that the Prescribed Boundary Problem
is a soft problem, in the sense that it obeys an $h$--principle;
using the relative version of the Gromov--Lees $h$--principle for Lagrangian
immersions we reduce the solvability of the Prescribed Boundary Problem to the
problem of extending a certain map from the boundary of a manifold with boundary
to its interior. This converts the solvability question into a problem in
homotopy theory, which we use obstruction theory to study.

We show that in general
there are obstructions to solving the Prescribed
Boundary Problem and therefore obstructions to smoothing our
Lagrangian singularities. To illustrate this we give explicit
examples of initial data for which the Prescribed Boundary Problem
has no solutions.
Moreover, in low dimensions we can identify the obstructions very
concretely and in particular identify situations in which all the obstructions vanish. An
easy corollary of this is that any SL conifold in dimensions $2$ or $3$ admits infinitely many topologically distinct oriented Maslov-zero Lagrangian
desingularizations (\fullref{t:lag:desing} and \fullref{c:SL:desing:CY}).
In particular, any SL cone in $\C^2$ and $\C^3$ admits infinitely many topologically
distinct oriented Maslov-zero smoothings (\fullref{c:SL:cones}).

The rest of this paper is organized as follows:

In \fullref{S:non:smooth} we prove various nonsmoothability results for
special Lagrangian singularities.
\fullref{SS:smooth:slg} serves as an introduction to
smoothing problems for isolated conical singularities of
minimal varieties, recalling known smoothing results
for volume-minimizing hypersurfaces.
\fullref{SS:slg}, \fullref{SS:slg:cones} and \fullref{SS:acsl} recall
basic definitions and facts about Lagrangian/special Lagrangian geometry,
special Lagrangian cones and asymptotically conical special Lagrangian $n$--folds respectively.
In \fullref{SS:non:smooth} we prove a first nonsmoothability result, \fullref{T:not:bounding}.
After recalling some basic results about $G$--invariant SL $n$--folds in \fullref{SS:slg:symmetry},
in \fullref{SS:non:smooth:2} we prove a further nonsmoothability result, \fullref{T:non:smooth}.

In \fullref{ss:hprinc:imm}--\fullref{ss:lagn:exist} we recall the Gromov--Lees
$h$--principle for Lagrangian immersions of closed manifolds
\cite{gromov:pdr,lees}, focusing on how it reduces an \textit{a
priori} geometric question to a topological one.  \fullref{s:maslov:class} discusses the Maslov class of a Lagrangian
immersion in $\C^n$ and uses the Lagrangian $h$--principle to prove
a Maslov-class realizability result, \fullref{p:anymaslov:ok}. \fullref{s:closed:examples} gives
examples to demonstrate how the $h$--principle provides an
effective tool to prove both existence and nonexistence results
for Lagrangian immersions. The material of this section is
well-known but is included both as an introduction to the
$h$--principle for $h$--principle novices and to establish notation,
terminology and some results that are needed in \fullref{s:bdy}. Those already comfortable with the Lagrangian
$h$--principle are advised to skim this section and proceed to
\fullref{s:bdy}.

In \fullref{ss:lagr:cobordism}--\fullref{ss:pbp} we describe the
\textit{Lagrangian Cobordism Problem\/}, some obstructions to solving it and
a refinement of it which we call the \textit{Prescribed Boundary
Problem}. Both problems are boundary value problems for Lagrangian
immersions in $\C^n$ of compact manifolds with boundary. In
\fullref{ss:pbp:rel:hprinc} we show how to use a relative
version of the $h$--principle for Lagrangian immersions (\fullref{T:hprinc:rel}) to analyze the solvability of the Prescribed
Boundary Problem. As a result we reduce the solvability of the
Prescribed Boundary Problem to a problem in algebraic topology
(see \fullref{L:extend:GL}); namely, does a certain map with
values in $\gl{n}$ extend continuously from the boundary to the
interior? In \fullref{S:obstruction} we recall a standard
obstruction theory framework for analyzing extension problems and
apply it to the case-at-hand. We pay particular attention to cases
in which either all obstructions vanish or where the obstructions
can be described particularly concretely. For example, in less
than $6$ complex dimensions there are exactly two obstructions to
extending the map in question; in \fullref{T:extend} we write
down these obstructions explicitly. One obstruction is closely
related to the Maslov class of the initial data; moreover in
dimensions $2$ and $3$ the other obstruction vanishes. This allows
us to prove particularly nice results for Maslov-zero initial data
in these dimensions.

In \fullref{ss:pbp:n=2}--\fullref{ss:pbp:n=4or5} we apply the
results of \fullref{s:bdy} to study the solvability of the
Prescribed Boundary Problem in complex dimensions $5$ and less,
taking care to show how our
\fullref{assumptions} can be satisfied. The main results are
\fullref{t:filling:n=2}, \fullref{t:filling:n=3} and
\fullref{t:filling:n=4,5}.

In \fullref{s:desing} we return to the soft Lagrangian desingularization problems
and address Questions (1)--(3) of \fullref{S:intro} using the results of \fullref{s:pbp:lowdim}.
We introduce a class of singular
oriented Lagrangian submanifolds called \textit{oriented
Lagrangian submanifolds with exact isolated singularities}. This
class includes the SL conifolds studied by Joyce in the series of
papers summarized in \cite{joyce:conifolds5}. In \fullref{t:lag:desing} we add the Maslov-zero assumption and apply our results
from \fullref{s:pbp:lowdim} to
prove desingularization and smoothing results for Lagrangian
submanifolds of this type in $\C^2$ and
$\C^3$. \fullref{s:general:ambient} extends these
results to almost Calabi--Yau manifolds, while \fullref{s:other:ends} addresses Question (3) above.

In \fullref{appendix:compare} we compare the work of this paper to existing work
in the Lagrangian and special Lagrangian literature.

The authors would like to thank IHES, Georgia Tech and Imperial College for their hospitality.
TP would like to thank NSF and MH would like to thank EPSRC for supporting their research.
The authors would also like to thank Rick Schoen for suggesting the problem, Yasha Eliashberg
for some $h$--principle assistance and Dominic Joyce for suggesting improvements to
an earlier version of the paper.

\section{Nonsmoothable special Lagrangian singularities}
\label{S:non:smooth}
The singularities of area-minimizing surfaces and codimension 1
volume-minimizing objects are by now fairly well-understood.
In higher dimension/codimension, examples show that
a far wider range of singular behaviour is possible and as a consequence
many codimension one results fail in higher codimension.
However, for special types of volume-minimizing objects, particularly
\textit{calibrated submanifolds\/} (or currents), one might still hope that some codimension one features survive.
In the case of smoothability of special Lagrangian singularities, however, we will see shortly that
such optimism would be misplaced.

\subsection{Smoothability of isolated conical singularities of minimal varieties}
\label{SS:smooth:slg}
Let $C$ be a minimal hypercone in $\R^n$ with $\textrm{Sing}(C)=(0)$.
Then the link of the cone, $\Sigma= C \cap \Sph^{n-1}$, is an embedded smooth compact $(n-2)$--dimensional minimal
manifold of $\Sph^{n-1}$.
It follows from the Maximum Principle that $\Sigma$ must be connected.
Since $\Sigma$ is an embedded hypersurface of $\Sph^{n-1}$ and the tangent bundle of $\Sph^{n-1}$
is stably trivial,
it follows that the tangent bundle of $\Sigma$ is also stably trivial.
This implies that the link of any minimal hypercone is nullcobordant
(see \fullref{SS:non:smooth})
and therefore admits some (topological) smoothing.

In fact, Hardt and Simon proved the following smoothing result
for volume-minimizing hypercones:
\begin{theorem}{\rm{\cite[Theorem 2.1]{hardt:simon}}}\qua
Let $C$ be a minimizing hypercone $C$ in $\R^n$ with $\textrm{Sing}(C)=(0)$ and
as above let $\Sigma$ denote the link of the cone.
Let $E$ denote one of the two components of $\R^n \setminus \overline{C}$.
Then $E$ contains a unique oriented embedded real analytic nonsingular volume-minimizing hypersurface $S$ with empty boundary
such that $\dist(S,0)=1$.
The hypersurface $S$ has the following additional properties:
\begin{itemize}
\item[(i)] For any $v\in E$, the ray $\{\lambda v: \lambda>0\}$  intersects
$S$ in a unique point, and the intersection is transverse.
\item[(ii)] (Up to orientation) any oriented minimizing hypersurface contained in $E$ coincides with some dilation of $S$.
\item[(iii)] $S$ is asymptotically conical in the sense
that there exists some $T>0$ and some compact subset $K \subset S$ so that
$S\setminus K$ is diffeomorphic to $(T,\infty) \times \Sigma$
and so that outside $K$, $S$ is graphical over the cone $C$ for some
$C^2$ graphing function $v$
(which satisfies some decay conditions at infinity -- see Hardt and Simon \cite[page 106 equation 1.9]{hardt:simon} for details).
\end{itemize}
\end{theorem}

Moving to higher codimension, in the spirit of the previous paragraphs
it is natural to ask whether calibrated cones have canonical desingularizations analogous to those of minimizing hypercones.
In particular, for special Lagrangian cones Rick Schoen originally posed
this question to the authors at the IPAM Workshop on Lagrangian submanifolds;
the authors would therefore like to thank him for suggesting this problem.
Unfortunately, it is not true that SL cones admit canonical SL smoothings;
we give examples below to show that there are SL
cones which admit no SL smoothings whatsoever.

To proceed further first we recall some basic definitions from Lagrangian and special Lagrangian geometry.

\subsection[Calibrations and special Lagrangian geometry in Cn]{Calibrations and special Lagrangian geometry in $\C^n$}
\label{SS:slg}
Let $(M,g)$ be a Riemannian manifold. Let $V$ be an oriented tangent $p$--plane on $M$, ie
a $p$--dimensional oriented vector subspace of some tangent plane $T_xM$ to $M$. The restriction
of the Riemannian metric to $V$, $g|_V$, is a Euclidean metric on $V$ which together with the
orientation on $V$ determines a natural $p$--form on $V$, the volume form $\vol_{V}$.
A closed $p$--form $\phi$ on $M$ is a \textit{calibration\/} \cite{harvey:lawson} on $M$ if for every oriented tangent $p$--plane $V$ on $M$
we have $\phi|_V \le \vol_V$. Let $L$ be an oriented submanifold of $M$ with dimension $p$. $L$ is a
\textit{$\phi$--calibrated submanifold\/} if $\phi|_{T_xL} = \vol_{T_xL}$ for all $x\in L$.

There is a natural extension of this definition to singular calibrated submanifolds using the
language of Geometric Measure Theory and rectifiable currents \cite[\S II.1]{harvey:lawson}.
The key property of calibrated submanifolds (even singular ones) is that they are \textit{homologically volume minimizing\/} \cite[Theorem II.4.2]{harvey:lawson}.
In particular, any calibrated submanifold is automatically \textit{minimal\/}, ie has vanishing mean curvature.

Let $z_1 = x_1 + i y_1, \ldots  ,z_n = x_n + i y_n$ be standard complex coordinates on $\C^n$ equipped with the Euclidean metric.
Let
\begin{equation}
\label{E:symplectic}
\omega = \frac{i}{2} \sum_{j=1}^n\!{dz_j {\wedge} d\bar{z}_j} = \sum_{j=1}^n\!{dx_j {\wedge} dy_j}  = d\lambda,
 \text{ where\ } \lambda = \frac{1}{2} \Big(\sum_{i=1}^{n}\! x_idy_i {-} y_idx_i\Big),
\end{equation}
be the standard symplectic $2$--form on $\C^n$.
Recall that an immersion  $i\co L^n \ra \R^{2n}$ of an $n$--manifold $L$ is \textit{Lagrangian\/} if $i^*\omega =0$ or, equivalently, if the
natural complex structure $J$ on $\R^{2n}=\C^n$ induces an
isomorphism between $TL$ and the normal bundle $NL$.

Define a complex $n$--form $\Omega$ on $\C^n$ by
\begin{equation}
\label{E:slg:form}
\Omega = dz_1 \wedge \ldots \wedge dz_n.
\end{equation}
Then $\Real{\Omega}$ and $\Imag{\Omega}$ are real $n$--forms on $\C^n$.
$\Real{\Omega}$ is a calibration on $\C^n$ whose calibrated submanifolds we call
\textit{\slg submanifolds\/} of $\C^n$, or SL $n$--folds for short.
There is a natural extension of \slg geometry to any Calabi--Yau manifold $M$ by replacing
$\Omega$ with the natural parallel holomorphic $(n,0)$--form on $M$.

More generally, let $f\co L \ra \C^n$ be a Lagrangian immersion of the oriented $n$--manifold $L$, and $\Omega$ be
the holomorphic $(n,0)$--form defined in \eqref{E:slg:form}.
Then $f^*\Omega$ is a complex $n$--form on $L$ satisfying $|f^*\Omega| = 1$ \cite[page 89]{harvey:lawson}. Hence we may write
\begin{equation}
\label{E:lagn:phase}
f^*\Omega = e^{i\theta}\vol_L\quad \text{on \ } L,
\end{equation}
for some \textit{phase function\/} $e^{i\theta}\co L \ra \Sph^1$.
We call $e^{i\theta}$ the \textit{phase of the oriented Lagrangian immersion $f$\/}.
$L$ is a SL $n$--fold in $\C^n$ if and only if the phase function $e^{i\theta}\equiv 1$.
Reversing the orientation of $L$ changes the sign of the phase function $e^{i\theta}$.
The differential $d\theta$ is a well-defined closed $1$--form on $L$ satisfying \cite[page 96]{harvey:lawson}
\begin{equation}
\label{mc:lagn:angle}
d\theta = \iota_H \omega,
\end{equation}
where $H$ is the mean curvature vector of $L$.
The cohomology class $[d\theta] \in H^1(L,\R)$
is closely related to the \textit{Maslov class\/} of the Lagrangian immersion defined in \fullref{s:maslov:class}.

In particular, \eqref{mc:lagn:angle} implies
that a connected component of $L$ is minimal if and only if the phase function $e^{i\theta}$
is constant. $e^{i\theta} \equiv 1$ is equivalent to $L$ being special Lagrangian;
if $e^{i\theta} \equiv e^{i\theta_0}$ for some constant $\theta_0 \in [0,2\pi)$
we will call $L$, \textit{$\theta_0$--special Lagrangian\/} or $\theta_0$-SL for short.
$\theta_0$-SL $n$--folds are calibrated submanifolds with respect to the real part of the rotated
holomorphic $(n,0)$--form, $\Omega_{\theta_0}= e^{-i\theta_0} \Omega$.

\subsection{Special Legendrian submanifolds and special Lagrangian cones}
\label{SS:slg:cones}
For any compact oriented embedded (but not necessarily connected) submanifold $\Sigma
\subset \Sph^{2n-1}(1)\subset \C^n$ define the \textit{cone on
$\Sigma$},
$$ C(\Sigma) = \{ tx: t\in \R^{\ge 0}, x \in \Sigma \}.$$
A cone $C$ in $\C^n$  (that is a subset invariant under dilations)
is $\textit{regular\/}$ if there exists
$\Sigma$ as above so that $C=C(\Sigma)$, in which case we call
$\Sigma$ the \textit{link\/} of the cone $C$. $C'(\Sigma):=C(\Sigma) - \{0\}$ is an
embedded smooth submanifold, but $C(\Sigma)$ has an isolated
singularity at $0$ unless $\Sigma$ is a totally geodesic sphere.

Let $r$ denote the radial coordinate on $\C^n$ and let
$X$ be the Liouville vector field
$$X= \frac{1}{2}r \frac{\partial}{\partial r} = \frac{1}{2}\sum_{j=1}^{n}{x_j \frac{\partial}{\partial x_j} + y_j \frac{\partial}{\partial y_j}}.$$
The unit sphere $\Sph^{2n-1}$ inherits a natural contact form
$$\gamma = \iota_X\omega | _{\Sph^{2n-1}} = {\sum_{j=1}^n{ x_j dy_j - y_j dx_j}}\Bigr|_{\Sph^{2n-1}}$$ from
its embedding in $\C^n$.
A regular cone $C$ in $\C^n$ is Lagrangian if and only if its link $\Sigma$ is Legendrian.
We call a submanifold $\Sigma$ of $\Sph^{2n-1}$ \textit{special Legendrian\/}
if the cone over $\Sigma$, $C'(\Sigma)$
is special Lagrangian in $\C^n$.

\subsection[Asymptotically conical SL n--folds and ACSL smoothings]{Asymptotically conical SL $n$--folds and ACSL smoothings}
\label{SS:acsl}
Roughly speaking, an \textit{asymptotically conical special Lagrangian $n$--fold\/} $L$ (ACSL) in $\C^n$
is a nonsingular SL $n$--fold which tends to a regular SL cone $C$ at infinity.
In particular, if $L$ is an ACSL submanifold then $\lim_{t\ra 0_+}tL=C$ for
some SL cone $C$. Hence ACSL submanifolds give rise to local models for how nonsingular
SL submanifolds can develop isolated singularities modelled on the cone $C$.
Conversely, as in \cite{joyce:conifolds5}, one can use ACSL submanifolds to desingularize
SL conifolds.

More precisely, we say that a closed, nonsingular SL submanifold $L$ of $\C^n$ is
\textit{asymptotically conical (AC) with decay rate $\lambda<2$ and cone $C$\/} if for some compact
subset $K \subset L$ and some $T>0$, there is a diffeomorphism
$\phi\co \Sigma \times (T,\infty) \ra L \setminus K$ such that
\begin{equation}
\label{E:acsl}
|\nabla^k(\phi-i)|=O(r^{\lambda-1-k}) \ \ \textrm{as\ } r\rightarrow\infty \textrm{\ for\ } k=0,1,
\end{equation}
where $i(r,\sigma) = r\sigma$ and $\nabla$ and $|\cdot|$ are defined using the cone metric
$g'$ on $C$.

Given a regular SL cone $C$ in $\C^n$ an important question is: does
$C$ admit an ACSL submanifold $L$ asymptotic to the given SL cone $C$?
We call this the \textit{asymptotically conical special Lagrangian smoothing problem for the SL cone $C$\/}
or for short the \textit{ACSL smoothing problem\/}.

From the point of view
of the global smoothing theory developed by Joyce in \cite{joyce:conifolds5}
it is natural to ask whether $C$ admits an ACSL smoothing with decay rate $\lambda \le 0$.
We call this the \textit{ACSL smoothing problem with decay\/}.

We will see shortly that neither version of the ACSL smoothing problem is always solvable.

\subsection{Nonsmoothable SL singularities I}
\label{SS:non:smooth}
Suppose that $C$ is a SL cone with link $\Sigma$ for which the ACSL smoothing problem is solvable.
Then clearly the smooth oriented manifold $\Sigma$ bounds a compact smooth oriented manifold $M$.
In other words, $\Sigma$ is \textit{oriented cobordant to the empty set\/}.
This is a special case of the oriented version of the Cobordism
Problem solved by Wall \cite{wall}, following the pioneering work
of Thom in the unoriented case \cite{thom}.

Wall determined the structure of the oriented cobordism ring
$\Omega_* = (\Omega_0, \Omega_1, \Omega_2, \ldots )$ \cite{wall}.
He showed that the oriented cobordism groups $\Omega_n$ are
trivial for $n=1, 2, 3, 6, 7$ and are nonzero in all other
dimensions (see \cite[page 203]{milnor:stasheff} for a table of
$\Omega_n$ up to $n=11$). Hence for $n\in \{1,2,3,6,7\}$ any
oriented $n$--manifold $\Sigma$ bounds an oriented $(n+1)$--manifold
$L$.

More generally, Wall \cite[page 306]{wall} proved that $[\Sigma]=0$
in the oriented cobordism ring if and only if all the Pontrjagin
and Stiefel--Whitney numbers of $\Sigma$ vanish (see Milnor and Stasheff \cite[Chapters
4 and 16]{milnor:stasheff} for a definition of these characteristic
numbers). For instance it follows from \eqref{e:pont:cpn} that
$\CP^2$ does not bound any oriented $5$--manifold. On the other
hand, if all the Pontrjagin and Stiefel--Whitney classes of
$T\Sigma$ are zero then of course so are the Pontrjagin and
Stiefel--Whitney numbers. Hence every stably parallelizable
manifold bounds.

Thus if we can find a SL cone $C$ whose link $\Sigma$ is not nullcobordant
then for purely topological reasons neither version of the ACSL smoothing
problem for $C$ is solvable. We will give such an example below and hence see that
there is no obvious SL analogue of the Hardt--Simon desingularization theorem for minimizing
hypercones.
\begin{theorem}
\label{T:not:bounding}
Let $A$ denote an $n \times n$ complex matrix. Define a map $\phi$ from $\sunitary{n}$ to $\text{Sym}(n,\C)$,
the symmetric $n \times n$ complex matrices, by
$$A \mapsto \frac{1}{\sqrt{n}} AA^T.$$
Then $\phi$ passes to a well-defined map $\hatphi$ on the quotient $\sunitary{n}/\sorth{n}$.
The induced map $\hatphi\co  \sunitary{n}/\sorth{n} \ra \text{Sym}(n,\C)$ gives a special Legendrian
embedding of the quotient $\Sigma=\sunitary{n}/\sorth{n}$ into the unit sphere of $\text{Sym}(n,\C)$
(and hence the cone over $\Sigma$ is special Lagrangian).

For $n=3$, $\Sigma = \sunitary{3}/\sorth{3}$ does not bound any compact
$6$--manifold
and in particular the ACSL smoothing problem with cone $C=C(\Sigma)$
has no solutions.
\end{theorem}
\begin{proof}
The fact that the above embedding is special Lagrangian
was proven by Cheng  \cite[Theorem 1.2]{cheng}.
Audin showed that the Stiefel--Whitney number $w_2
w_3(T\Sigma)=1\neq 0$ for $n=3$ \cite[page 191]{audin}. Since $\Sigma$ has a nonzero
Stiefel--Whitney number it does not bound any smooth $6$--manifold.
\end{proof}

\begin{remark}
In fact, Ohnita \cite[Theorem 2.2]{ohnita:slg} has recently studied the \textit{stability index\/}
of these (very symmetric) SL cones. He showed that the SL cone over $\sunitary{3}/\sorth{3}$  is \textit{strictly
stable} in the sense of \cite{haskins:invent,joyce:conifolds5}.
Hence this SL cone would be a good candidate for a ``common singularity type''.
On the other hand what we have just said shows that there can be no SL
smoothings of $X$.
\end{remark}

Next we exhibit infinitely many SL cones for which there is no solution
to the ACSL smoothing problem with decay. These nonsmoothing results will also show that
there are obstructions to the solvability of the ACSL smoothing problem
beyond the cobordism class of the link of the cone.

First we need some auxiliary results on SL $n$--folds with symmetry and moment maps.

\subsection[G--invariant SL submanifolds]{$G$--invariant SL submanifolds}
\label{SS:slg:symmetry}

In symplectic geometry continuous
group actions which preserve $\omega$ are often induced by a collection of functions,
the \textit{moment map\/} of the action. We describe the
situation in $\C^n$, where the situation is particularly simple.

Let $(M,J,g,\omega)$ be a \ka manifold, and let $\text{G}$ be a Lie group
acting smoothly on $M$ preserving both $J$ and $g$ (and hence $\omega$).
This induces a linear map $\phi$ from the Lie algebra $\mathfrak{g}$ of $\text{G}$ to
the vector fields on $M$. Given $x\in \mathfrak{g}$, let $v=\phi(x)$ denote the
corresponding vector field on $M$. Since $\mathcal{L}_v\omega = 0$, it follows
from Cartan's formula that $\iota_v\omega$ is a closed $1$--form on $M$. If
$H^1(M,\R)=0$ then there exists a smooth function $\mu^x$ on $M$, unique up to a constant, such that
$d\mu^x = \iota_v\omega$. $\mu^x$ is called a moment map for $x$, or a Hamiltonian function
for the Hamiltonian vector field $v$.
We can attempt to collect all these functions $\mu^x$, $x\in \mathfrak{g}$, together to make a moment map for the whole action.
A smooth map $\mu\co M \ra \mathfrak{g}^{*}$ is called a \textit{moment map\/} for the action of $\text{G}$ on $M$ if
\begin{enumerate}
\item[(i)]
$\iota_{\phi(x)}\omega = \langle x,d\mu \rangle$ for all $x\in \mathfrak{g}^{*}$ where $\langle,\rangle$ is the natural pairing of $\mathfrak{g}$ and $\mathfrak{g}^*$;
\item[(ii)]
$\mu\co M \ra \mathfrak{g}^*$ is equivariant with respect to the $\text{G}$--action on $M$ and the coadjoint $\text{G}$--action on $\mathfrak{g}^*$.
\end{enumerate}
In general there are obstructions
to a symplectic $\text{G}$--action admitting a moment map.
The subsets $\mu^{-1}(c)$ are the \textit{level sets\/} of the moment map. The \textit{centre\/} $Z(\mathfrak{g}^*)$ is the subspace
of $\mathfrak{g}^*$ fixed by the coadjoint action of $\text{G}$. Property (ii) of $\mu$ implies
that a level set $\mu^{-1}(c)$ is $\text{G}$--invariant if and only if $c\in Z(\mathfrak{g}^*)$.

When a moment map for a symplectic group action does exist, it is
a very useful tool for studying $G$--invariant isotropic submanifolds of $(M,\omega)$.
Using property (i) above it is easy to see that any $G$--invariant isotropic submanifold
of $(M,\omega)$ must be contained in some level set $\mu^{-1}(c)$ of the moment map.
Using property (ii) one can check that $c\in Z(\mathfrak{g}^*)$. Hence we have
the useful fact that: any $G$--invariant isotropic submanifold of $(M,\omega)$ is contained in
the level set $\mu^{-1}(c)$ for some $c\in Z(\mathfrak{g}^*)$.

The group of automorphisms of $\C^n$ preserving $g$, $\omega$ and $\Omega$
is $\text{SU(n)} \ltimes \C^n$, where $\C^n$ acts by translations.
Since $\C^n$ is simply connected, each element $x$ in the Lie algebra $\mathfrak{su}(n) \ltimes \C^n$ has a moment map $\mu^x\co \C^n\ra \R$. More concretely, we can describe these functions as follows:
a function is said to be a \textit{harmonic Hermitian quadratic\/} if it is of the form
\begin{equation}
\label{hermitian:harmonic}
f = c + \sum_{i=1}^n{(b_i z_i + \bar{b}_i\bar{z}_i)} + \sum_{i,j=1}^n{a_{ij}z_i\bar{z}_j}
\end{equation}
for some $c\in \R, \ b_i, \ a_{ij} \in \C \textrm{\ with \ } a_{ij}=\bar{a}_{ji} \textrm{\ and\ } \sum_{i=1}^n{a_{ii}}=0.$
A harmonic Hermitian quadratic with $c=0$, $a_{ij}=0$ for all $i$ and $j$ corresponds to the moment
map of a translation.
A harmonic Hermitian quadratic with $c=0$, $b_i=0$ for all $i$ corresponds to the moment map
of the element $A=\sqrt{-1}\,(a_{ij})\in \mathfrak{su}(n)$.
Moreover, these functions satisfy conditions (i) and (ii) above, and
thus give rise to a moment map for the whole group action.

Since this group action also preserves the condition to be special Lagrangian,
the restriction to \slg submanifolds of the moment map of any $\mathfrak{su(n)} \ltimes \C^n$ vector field
enjoys some special properties.
\begin{lemma}{\rm{(Joyce \cite[Lemma 3.4]{joyce:conifolds2} and Fu \cite[Theorem 3.2]{fu})}}\qua
\label{harmonic:momentmap}
Let $\mu\co  \C^n \ra \R$ be a moment map for a vector field $x$ in $\mathfrak{su}(n) \ltimes \C^n$.
Then the restriction of $\mu$ to any \slg $n$--fold in $\C^n$ is a harmonic function
on $L$ with respect to the metric induced on $L$ by $\C^n$.

Conversely, a function $f$ on $\C^n$ is harmonic on every SL $n$--fold
in $\C^n$ if and only if
\begin{equation}
\label{fu:eqn}
d( \iota_{X_f}\Imag(\Omega)) = 0
\end{equation}
where $X_f = -J\nabla f$ is the Hamiltonian vector field associated with $f$.
For $n\ge3$, $f$ satisfies \eqref{fu:eqn} if and only if $f$ is a harmonic Hermitian quadratic.
\end{lemma}
\noindent
In other words, any \slg $n$--fold $L$ in $\C^n$ automatically has certain distinguished harmonic functions.
We can use \fullref{harmonic:momentmap} together with the Maximum
Principle to conclude that a compact \slg submanifold with boundary
inherits all the symmetries of the boundary.
\begin{prop}{\rm{(Haskins \cite[Proposition 3.6]{haskins:invent})}}\qua
\label{slg:bdy:symmetry}
Let $L^n$ be a compact connected \slg submanifold of $\C^n$ with boundary $\Sigma$.
Let $\text{G}$ be the identity component of the subgroup of $\text{SU(n)} \ltimes \C^n$ which
preserves $\Sigma$. Then $\text{G}$ admits a moment map $\mu\co \C^n\ra \mathfrak{g}^*$,
both $\Sigma$ and $L$ are contained in $\mu^{-1}(c)$ for some $c\in Z(\mathfrak{g}^*)$
and $\text{G}$ also preserves $L$.
\end{prop}
Here is an asymptotically conical analogue of this result:
\begin{prop}{\rm{(Haskins \cite[Propostion B.1]{haskins:invent})}}\qua
\label{acslg:symmetry}
Let $L$ be an \acslg submanifold $L$ with rate $\lambda$ and regular \slg cone $C$.
Let $\text{G}$ be the identity component of the subgroup of $\text{SU(n)}$ preserving $C$. Then
\begin{enumerate}
\item[(i)] $\text{G}$ admits a moment map $\mu$ and $C\subset \mu^{-1}(0)$;

\item[(ii)] if $\lambda<0$, then $L\subset \mu^{-1}(0)$ and $\text{G}$ also preserves $L$;

\item[(iii)] if $\lambda =0$ and $L$ has one end, then $L\subset \mu^{-1}(c)$ for some $c\in Z(\mathfrak{g}^*)$ and $\text{G}$ preserves $L$.
\end{enumerate}
\end{prop}

The following lemma shows that associated to any SL cone $C$ there is a natural $1$--parameter family of ACSL submanifolds $L_t$,
asymptotic with rate $\lambda=2-n$ to the union of the two SL cones $C \cup e^{i\unfrac{\pi}{n}}C$.
Hence there always exist natural ``two-ended'' desingularizations of any SL cone $C$. It may be useful to emphasize that this construction does not provide a solution to the ACSL smoothing problem for $C$, unless $C$ is $(e^{i\unfrac{\pi}{n}})$--invariant.
Notice that this lemma
is compatible with the fact that the disjoint union of any oriented manifold $\Sigma$ together with $-\Sigma$
is always nullcobordant.
\begin{lemma}{\rm{\cite[Lemma B.3]{haskins:invent}}}\qua
\label{L:acsl}
Let $C$ be any regular SL cone in $\C^n$ with link $\Sigma$. For any $t \neq 0 \in \R$ define $L_t$ as
$$L_t= \left\{ \sigma\in \C^n: \sigma \in \Sigma, z\in \C, \textrm{\ with\ } \imag{z^n}=t, \ \arg{z}\in (0,\tfrac{\pi}{n}) \right\}.$$
Then $L_t$ is an ACSL $n$--fold with rate $\lambda=2-n$, and cone $C  \cup e^{i\unfrac{\pi}{n}}C$.
Moreover, $L_t$ has the same symmetry group as $C$.
\end{lemma}

\subsection{Nonsmoothable SL singularities II}
\label{SS:non:smooth:2}
\fullref{acslg:symmetry} and \fullref{L:acsl} show that any ACSL smoothing (with appropriate decay)
of a $G$--invariant SL cone $C$ is also $G$--invariant, and that $C$ has a canonical $1$--parameter family of
$G$--invariant two-ended ACSL smoothings.
For certain group actions $G$ we will see that
the only $G$--invariant SLG $n$--folds are the SL cones $C$ themselves and their canonical
two-ended ACSL smoothings.
Hence we will be able to prove that certain $G$--invariant SL cones $C$
admit no (one-ended) ACSL smoothings (with decay).

The following three infinite families of very symmetric SL cones have been studied previously by
a number of authors \cite{castro:urbano:sym,cheng,naitoh,ohnita:slg}.
\begin{example}
\label{e:slg:su} For $n\ge 3$,
let $\mathfrak{gl}(n,\C)$ denote the space of all $n \times n$ complex matrices equipped
with the hermitian inner product $\langle\langle A,B\rangle\rangle:=\trace{AB^*}$,
and let $\sunitary{n}$ act on $\mathfrak{gl}(n,\C)$ by $(A,B) \mapsto AB$.
Then the map $\phi\co  \sunitary{n} \ra \mathfrak{gl}(n,\C)$ given by
$$\phi(A):= \frac{1}{\sqrt{n}}A$$
gives an $\sunitary{n}$--invariant special Legendrian submanifold of the unit
sphere in $\mathfrak{gl}(n,\C)$ diffeomorphic to $\sunitary{n}$.
\end{example}
\begin{example}
\label{e:slg:su:so}
For $n\ge3$, let $\textrm{Sym}(n,\C)$ denote the space of all symmetric $n \times n$ complex matrices equipped
with the hermitian inner product induced from $\mathfrak{gl}(n,\C)$.
Let $\sunitary{n}$ act on $\textrm{Sym}(n,\C)$ by $(A,B) \mapsto ABA^T$.
Then the map $\phi\co  \sunitary{n} \ra \textrm{Sym}(n,\C)$ given by
$$\phi(A):= \frac{1}{\sqrt{n}}AA^T$$
induces an $\sunitary{n}$--invariant special Legendrian submanifold of the unit
sphere in $\textrm{Sym}(n,\C)$ diffeomorphic to $\sunitary{n}/\sorth{n}$.
Note, these are the same examples that already appeared in \fullref{T:not:bounding}.
\end{example}
\begin{example}
\label{e:slg:su:sp}
For $n\ge 2$, let $\mathfrak{so}(2n,\C)$ denote the space of all skew-symmetric $2n \times 2n$ complex matrices equipped
with the hermitian inner product induced from $\mathfrak{gl}(2n,\C)$.
Let $\sunitary{2n}$ act on $\mathfrak{so}(2n,\C)$ by $(A,B) \mapsto ABA^T$.
Then the map $\phi$ from $\sunitary{2n}$ to $\mathfrak{so}(2n,\C)$ given by
$$\phi(A):= \frac{1}{\sqrt{n}}A J_n A^T,
\quad \text{where\ }
J_n=\left(\begin{array}{c|c}
0 & -I_n \\
\hline
I_n & 0
\end{array}\right),
$$
induces an $\sunitary{2n}$--invariant special Legendrian submanifold of the unit
sphere in $\mathfrak{so}(2n,\C)$ diffeomorphic to $\sunitary{2n}/\textrm{Sp}(n)$.
\end{example}

\begin{theorem}
\label{T:non:smooth}
Let $C$ be any of the SL cones described in
\fullref{e:slg:su}, \fullref{e:slg:su:so} and \fullref{e:slg:su:sp}. Then
the ACSL smoothing problem with decay is not solvable for the cone $C$.
\end{theorem}
\begin{proof}
For \fullref{e:slg:su}--\fullref{e:slg:su:sp}, Castro--Urbano \cite[\S 2]{castro:urbano:sym} proved that the only $G$--invariant SL submanifolds
(with the group $G$ and its action as described in the examples) are the cones themselves and their canonical two-ended ACSL smoothings.
However, by \fullref{acslg:symmetry}(iii) any ACSL smoothing of $C$ with decay rate $\lambda \le 0$ is
also $G$--invariant. Hence no such ACSL smoothing of $C$ with decay rate $\lambda \le 0$ exists.
\end{proof}
Here is a way to think about why we should expect a result like \fullref{T:non:smooth} to hold. In the examples above, the common feature is that we have a low
cohomogeneity action of the nonabelian
group $G$ for which $Z(\mathfrak{g}^*)=(0)$. Let $\mu$ denote the associated moment map.
Then we know that all $G$--invariant Lagrangian submanifolds are contained in $\mu^{-1}(0)$.
Thus there are rather few $G$--invariant Lagrangian submanifolds for these groups $G$.

If $C$ is a $G$--invariant SL cone and $G$ has cohomogeneity one (as do the above examples)
then any ACSL smoothing with decay will also be $G$--invariant with cohomogeneity one.
Thus to obtain (one-ended) ACSL smoothings of $C$ it is necessary that there be $\omega$--isotropic orbits
of $G$ whose dimension is neither maximal (corresponding to the link of the cone) nor minimal (when the orbit is just the origin of $\C^n$).
If isotropic orbits cannot collapse in this intermediate way then there is no way to obtain $G$--invariant
one-ended ACSL $n$--folds. For instance, in \fullref{e:slg:su} it is easy to see that
$\sunitary{n}$ acts freely on $\mu^{-1}(0)$, except at the origin. Hence all the $G$--orbits in $\mu^{-1}(0)$
except the origin are diffeomorphic to $\sunitary{n}$, and thus there are no ``intermediate'' isotropic orbits.

However, if instead $G$ is abelian, then any level set of the corresponding moment map $\mu$
is $G$--invariant, and not just the zero level set $\mu^{-1}(0)$. Thus in the abelian case there are many more $G$--invariant Lagrangian
submanifolds and many more isotropic $G$--orbits.
In particular, if $G$ is the maximal abelian subgroup $\mathbb{T}^{n-1} \subset \sunitary{n}$ acting in the obvious way on $\C^n$,
then there are ``intermediate''
isotropic $G$--orbits. As a result the $\mathbb{T}^{n-1}$--invariant SL cone $C$ over the Legendrian Clifford torus
$T^{n-1}$ does admit one-ended ACSL smoothings.
In fact, one can classify all the one-ended ACSL smoothings with decay of the
Legendrian Clifford torus $\mathbb{T}^{n-1} \subset \Sph^{2n-1}$ this way;
see \cite[\S 10]{joyce:conifolds5} for a detailed description
of the $n=3$ case.

One can use the same idea for other group actions $G$. Provided $G$ is ``large enough'' then we can expect
to classify \textit{all\/} $G$--invariant SL $n$--folds and hence to classify all ACSL smoothings with decay of
any such $G$--invariant SL cone.

In \fullref{e:slg:su} the link of the SL cone is diffeomorphic to $\sunitary{n}$.
Since any Lie group is parallelizable, all its Pontrjagin and Stiefel--Whitney classes vanish
and hence it is nullcobordant. Hence we have infinitely many SL cones where there is no
differential-topological obstruction to smoothing the singularity, but which nevertheless
admit no ACSL smoothing (with decay).

In this case a natural weaker question to ask is: do the nullcobordant
special Legendrian links $\Sigma$ above admit Lagrangian or
Maslov-zero Lagrangian smoothings? We will see that we can answer the immersed
version of this question using a relative version of the Gromov--Lees $h$--principle
for Lagrangian immersions. We will see that in general there are obstructions
to Lagrangian smoothing that are of an algebro-topological nature which go
beyond the cobordism type of the link.

In order to describe these obstructions
we need to review the $h$--principle for Lagrangian immersions, beginning with the
version for closed manifolds. Although much of this material is well-known to symplectic
topologists, we expect that many readers whose background is in minimal submanifolds or special Lagrangian
geometry will be unfamiliar with it.

\section[The h--principle for Lagrangian immersions of closed manifold]{The $h$--principle for Lagrangian immersions \\of closed manifolds}
\label{S:h:princ:closed}
Let $L^n$ be a compact connected (not necessarily orientable) manifold without boundary.
A natural question in symplectic geometry is:

When does $L$ admit a Lagrangian immersion into $\R^{2n}$?

Gromov and Lees showed that Lagrangian immersions of closed
manifolds satisfy the so-called $h$--principle, where h stands for
homotopy. The idea is that for a Lagrangian immersion to exist
certain conditions of an algebro-topological nature must be
satisfied, and saying that the $h$--principle holds for Lagrangian
immersions means that these necessary topological conditions are
in fact sufficient.

Hence the previous, apparently geometric question about the
existence of Lagrangian immersions reduces to a topological
problem which in many cases can be solved.

\subsection[The h--principle for immersions of closed manifolds]{The $h$--principle for immersions of closed manifolds}
\label{ss:hprinc:imm} We begin with a description of the
$h$--principle for immersions of a closed manifold, ie
for the moment we drop the Lagrangian condition. Let $V$ and $W$
be smooth manifolds of dimension $n$ and $q$ respectively. Suppose
that $n \le q$ and that there exists an immersion $f\co  V \ra W$.
Then the differential $df\co TV \ra TW$ is a smooth bundle map
(ie a smooth map which fibrewise is linear) which is
injective on each fibre. This motivates the following definition:
\begin{definition}
A \textit{monomorphism\/} $F\co TV \rightarrow TW$ is a smooth
bundle morphism which is fibrewise injective.  We will usually
denote by $f=\textrm{bs}(F)$ the underlying ``base map" $V\rightarrow W$.
Notice that $f$ is a smooth map, but it need not be an immersion.
\end{definition}
Eliashberg and Mishachev \cite{eliashberg} call such a
monomorphism a \textit{formal solution\/} of the immersion problem.
Clearly, the existence of a monomorphism $F\co TV \ra TW$ (a formal
solution) is a necessary condition for the existence of an
immersion $f\co V \ra W$ (a genuine solution). A monomorphism of the
form $F=df$ is also sometimes called \textit{holonomic\/}. One says
that a differential relation $\mathcal{R}$ (eg the
immersion relation) satisfies the $h$--principle if every formal
solution of $\mathcal{R}$ (a monomorphism $TV\ra TW$ in our case)
is homotopic in the space of formal solutions to a genuine
solution of $\mathcal{R}$ (an immersion $V\ra W$). Hirsch
\cite{hirsch} extending work of Smale \cite{smale} proved that the
$h$--principle holds for the differential relation associated with
immersions of a closed $n$--manifold $V$ into a $q$--manifold $W$
provided $q>n$. In the case where $V$ is an $n$--manifold and
$W=\R^{q}$ we know from the results of Whitney \cite{whitney:hard}
that any smooth $n$--manifold ($n>1$) can be immersed in $\R^{q}$ for any
$q\ge 2n-1$. In particular, this gives a very indirect argument
that we must always be able to construct a monomorphism from $TV$
to $T\R^{q}$, for $q\ge 2n-1$ whatever the $n$--manifold $V$.
Perhaps more importantly it shows that the Smale--Hirsch
$h$--principle is (in the case of immersions into Euclidean space)
of greatest interest in the case of relatively low codimension.

\subsection[The h--principle for Lagrangian immersions]{The $h$--principle for Lagrangian immersions}
\label{ss:hprinc:lagn} To understand the statement of the
$h$--principle for Lagrangian immersions first we have to
understand what the appropriate notion of a formal solution to the
Lagrangian immersion problem is.

Suppose as in the previous section that $V$ and $W$ are smooth
manifolds of dimension $n$ and $q$ respectively and that $q>n$.
Any monomorphism $F\co TV \ra TW$ induces a natural map
$GF\co  V \ra\textrm{Gr}_n(W)$, where $\textrm{Gr}_n(W)$ denotes the
Grassmannian bundle of $n$--planes in $W$. In particular, for any
immersion $f\co V \ra W$ we have an associated (tangential) Gauss map
written $Gdf$. If $(W,\omega)$ is a symplectic manifold then one
can define a number of natural subbundles of the Grassmannian
$n$--plane bundles. In particular, if $q=2n$ we can define the
subset $\textrm{Gr}_{\textrm{lag}}(W)$ of all Lagrangian
$n$--planes in $(W,\omega)$. Clearly, if $f\co V^n \ra (W^{2n},\omega)$
is a Lagrangian immersion then the image of the
associated (tangential) Gauss map $Gdf$ is contained in the
subbundle $\textrm{Gr}_{\textrm{lag}}(W)$ of $\textrm{Gr}_n(W)$.

Similarly, one can also define a subbundle
$\textrm{Gr}^+_{\textrm{lag}}(W) \subset \textrm{Gr}^+_n(W)$
corresponding to the set of all oriented Lagrangian planes inside
the set of all oriented $n$--planes. If $f\co V^n \ra (W^{2n},\omega)$
is a Lagrangian immersion of an oriented manifold $V$ then the
associated Gauss map $Gdf$ is contained in
$\textrm{Gr}^+_{\textrm{lag}}(W)$.
Similarly, if $(W,J)$ is an almost-complex manifold then there are other natural subbundles
of the Grassmannian $n$--plane bundles. For instance, a linear subspace $S$ of a complex vector space $(W,J)$
is \textit{totally real\/} if $JS \cap S =(0)$. Hence, in an almost complex manifold $(W^{2n},J)$ there is
a subbundle $\textrm{Gr}_{\textrm{real}}(W) \subset \textrm{Gr}_n(W)$ such that any totally real immersion
$f\co V^n \ra (W^{2n},J)$ has tangential (Gauss) map $Gdf$ contained in $\textrm{Gr}_{\textrm{real}}(W)$.
If $J$ and $\omega$ are compatible then $\textrm{Gr}_{\textrm{lag}}(W) \subset \textrm{Gr}_{\textrm{real}}(W)$,
so that any Lagrangian immersion is totally real. The converse is false.
There is also an oriented version of this totally real Grassmannian bundle.
\begin{definition}
Let $V$ be a smooth (not necessarily orientable) $n$--manifold and let $(W,\omega)$ be a symplectic $2n$--manifold.
A \textit{Lagrangian monomorphism\/} $F$ is a monomorphism $F\co TV \ra TW$
such that $GF(V) \subset \textrm{Gr}_{\textrm{lag}}(W)$.
Similarly, one can define a \textit{totally real monomorphism\/} in an almost-complex manifold $(W,J)$
by replacing $\textrm{Gr}_{\textrm{lag}}(W)$ with $\textrm{Gr}_{\textrm{real}}(W)$.
\end{definition}
A Lagrangian monomorphism is the correct notion of a formal
solution of the Lagrangian immersion problem. For the case of
Lagrangian immersions of closed manifolds into $(\R^{2n},\omega)$
with its standard symplectic structure the following version of
the $h$--principle, due to Gromov \cite{gromov:pdr} and Lees
\cite{lees}, holds.
\begin{theorem}[Lagrangian $h$--principle for closed manifolds]
\label{T:h:princ:closed}
Let $L$ be a smooth closed $n$--manifold. Suppose there is a Lagrangian monomorphism $F\co TL\rightarrow T\R^{2n}$. Then there exists a family of Lagrangian monomorphisms
$F_t\co [0,1]\times TL\rightarrow T\R^{2n}$ such that $F_0=F$ and
$F_1$ is holonomic; ie $F_1=d\tilde{f}$.
In particular, the base map $\textrm{bs}(F_1)=\tilde{f}\co L\rightarrow\C^n$ is a Lagrangian
immersion. Furthermore, $\tilde{f}$ is exact.
\end{theorem}
\begin{remark}[Exactness of the resulting Lagrangian immersion]
Let $f\co \Sigma\rightarrow \R^{2n}$ be an immersion of a manifold $\Sigma^k$,
of dimension $k\leq n$. Recall that $f$ is said to be an \textit{exact\/} immersion if the $1$--form
$f^*\lambda$ is exact, ie $f^*\lambda=dz$, for some $z\in C^\infty(\Sigma,\R)$.
Equivalently $f$ is exact if $\int_\gamma f^*\lambda=0$ for all closed curves
$\gamma$ in $\Sigma$.
Since $\omega=d\lambda$, any exact immersion is automatically $\omega$--isotropic.

If $f$ is $\omega$--isotropic, $f^*\lambda$ is a closed $1$--form on
$\Sigma$ and the class
$[f^*\lambda]\in H^1(L;\R)$ is called the \textit{symplectic area class\/} of the
isotropic immersion $f$.
The isotropic immersion $f$ is exact if and only if the
symplectic area class of $f$ vanishes.
Hence any Lagrangian immersion of a simply
connected manifold is exact. Also, any Lagrangian submanifold is
``locally exact" near any smooth point, because any closed $1$--form
is locally exact. In \fullref{s:desing} we
will see that this is also true near certain types of
isolated singularities.

There is an important reformulation of the condition that an
immersion $f\co  \Sigma \ra \R^{2n}$ be exact; $f$ is exact if and
only if it can be lifted to an immersion in $\R^{2n+1}$ which is
isotropic with respect to its standard contact structure $\alpha =
dz - \lambda$. In particular, (up to constants) there is a
one-to-one correspondence between exact Lagrangian immersions in
$(\R^{2n},\omega)$ and Legendrian immersions in
$(\R^{2n+1},\alpha)$. The proof of \fullref{T:h:princ:closed}
rests on the fact that any Lagrangian monomorphism into $\R^{2n}$
lifts to a Legendrian monomorphism into $\R^{2n+1}$; one can then
apply the corresponding result for Legendrian monomorphisms; see
\cite[\S 16.1.3]{eliashberg}. Hence the previous observation
explains why the Lagrangian immersion that we obtain from the
$h$--principle in \fullref{T:h:princ:closed} must be exact.
\end{remark}
Recall that two Lagrangian immersions are said to be
\textit{regularly homotopic\/} if there exists a homotopy through
Lagrangian immersions. We will also need the following parametric
version of the Lagrangian $h$--principle (also due to Gromov and
Lees).
\begin{theorem}[Parametric Lagrangian $h$--principle]
\label{T:h:princ:parametric}
Let $f_0,f_1\co L\ra\R^{2n}$ be Lagrangian immersions. If the monomorphisms $F_0:=df_0$ and $F_1:=df_1$
are homotopic through Lagrangian monomorphisms, then they are also homotopic through holonomic Lagrangian monomorphisms;
ie\ $f_0$ and $f_1$ are regularly homotopic.
\end{theorem}
\begin{remark} Notice that if $f_0,f_1$ are exact, a proof of \fullref{T:h:princ:parametric}
can be obtained along the same lines as before: we can lift the
curve of Lagrangian monomorphisms $F_s$ into $\R^{2n+1}$ and get a
curve of Legendrian monomorphisms which are holonomic for $s=0$
and $s=1$. The parametric Legendrian $h$--principle now yields a
curve of holonomic Legendrian monomorphisms which coincide with
the given ones for $s=0$ and $s=1$; this curve projects back down
to $\R^{2n}$ yielding the desired curve of holonomic Lagrangian
monomorphisms.

If $f_0$ and/or $f_1$ are not exact, the proof is more
complicated. One first needs to prove the equivalent of \fullref{T:h:princ:closed} for general exact ambient spaces
$(N,\omega=d\alpha)$. Then one proves an analogue of \fullref{T:h:princ:parametric} for exact $f_0,f_1\co L\ra N$ and
$s$--dependent $1$--forms $\alpha_s$ (this time we lift the given data
into $N\times\R$ endowed with the appropriate $s$--dependent
contact structure). It turns out, using the parametric
$h$--principle for immersions, that the problem of \fullref{T:h:princ:parametric} is equivalent to the following.

Let $N^{2n}$ be a rank $n$ vector bundle over $L$, endowed with a
curve of exact symplectic forms $\omega_s=d\lambda_s$ for
$s\in[0,1]$. Let $i\co L\ra N$ be the zero section and $\pi\co N\ra L$
be the projection. Assume we are given a curve of monomorphisms
$F_s\co TL\ra TN$ such that $F_s$ is Lagrangian with respect to
$\omega_s$ and $F_0=F_1=di$. Then $F_s$ can be homotoped, rel
$s=0$ and $s=1$, to a curve of holonomic Lagrangian monomorphisms.

The proof of this last claim is rather simple. Let us choose a
curve of closed $1$--forms $\tilde{\lambda}_s$ which coincide with
$\pi^*(\lambda_s|_{i(L)})$ when $s=0$ and $s=1$. Let
$\alpha_s:=\lambda_s-\tilde{\lambda}_s$. Notice that
$d\alpha_s=\omega_s$ so the symplectic structure has not changed
and the $F_s$ are still Lagrangian. However $F_0=di$ and $F_1=di$
are now exact with respect to the $1$--forms $\alpha_s$, so we can
apply the simpler parametric $h$--principle seen above.
\end{remark}
\begin{remark}[Regular homotopies vs Hamiltonian isotopies]\label{r:reg:ham:isotopies}
Suppose $f\co L^n \ra \R^{2n}$ is a nonexact Lagrangian
immersion. Since $df$ is a Lagrangian monomorphism, we can first
apply \fullref{T:h:princ:closed} to find an exact immersion
$\tilde{f}$; we can then apply \fullref{T:h:princ:parametric}
to prove that $f$ is regularly homotopic to $\tilde{f}$. There is
however a finer classification of Lagrangian immersions than
that given by the regular homotopy classes, namely the
classification up to Hamiltonian isotopy. The condition that a
Lagrangian immersion be exact or nonexact is invariant under
Hamiltonian isotopy. This shows that using the $h$--principle to
perturb a given Lagrangian submanifold will typically change the
Hamiltonian isotopy class of that submanifold.
\end{remark}
\begin{remark}[Lagrangian and totally real embeddings]
\label{r:embeddings} It is natural to ask whether the
$h$--principle also holds for Lagrangian embeddings and not just
immersions. By Whitney's hard embedding theorem any smooth
$n$--manifold $L$ admits an embedding into $\R^{2n}$. There is,
however, a simple topological obstruction to a closed $n$--manifold
admitting a Lagrangian embedding in $\R^{2n}$ which we now
describe.

Suppose that $L^n$ admits a Lagrangian embedding into $\R^{2n}$.
The self-intersection number of any compact submanifold $L^n$ of
$\R^{2n}$ is zero since we may use translations to disjoin $L$
from itself. Since $L$ is embedded we can identify a neighbourhood
of $L\subset\R^{2n}$ with a neighbourhood of $L$ inside its normal
bundle $NL$; thus the Euler number of the normal bundle of $L$ is
also zero. Since the Lagrangian condition allows us to identify
$NL$ with $TL$, this proves that the Euler number of the tangent
bundle of $L$, ie the Euler characteristic of $L$, must
also be zero. This shows, for example, that the only orientable
surface which could admit a Lagrangian embedding into $\R^4$ is
the torus $T^2$, and the product $\Sph^1 \times \Sph^1 \subset \C
\times \C$ gives such a Lagrangian embedding. Notice that we
cannot conclude anything from this Euler characteristic
information about Lagrangian embeddings of closed orientable
$3$--manifolds in $\R^6$, and in particular of $S^3$.

In fact, Gromov proved that Lagrangian embeddings do not satisfy
the $h$--principle, as an application of his theory of
$J$--holomorphic curves in symplectic manifolds
\cite{gromov:jholo}. In particular, he proved that there are no
exact Lagrangian embeddings of closed $n$--manifolds in
$(\R^{2n},\omega)$. For example, there is no Lagrangian embedding
of $S^3$ in $\R^6$.

This also shows that if we apply \fullref{T:h:princ:closed} to
perturb a (nonexact) embedded Lagrangian submanifold, the
resulting exact Lagrangian submanifold will typically only be
immersed. On the other hand, there are totally real embeddings of
$S^3$ in $\R^6$. In fact both totally real immersions and totally
real embeddings satisfy versions of the $h$--principle (see
Gromov \cite{gromov:pdr} or Eliashberg and Mishachev \cite[\S 19.3.1, \S 19.4.5]{eliashberg}).
\end{remark}
\subsection{Existence and classification of Lagrangian immersions}
\label{ss:lagn:exist} As we have just seen, the $h$--principle
reduces questions about the existence of Lagrangian immersions to
questions about the existence of Lagrangian monomorphisms. In this
section we show that it suffices to prove the existence of totally
real monomorphisms; this last problem has a straightforward
solution which leads to a classification of Lagrangian immersions
in terms of topological data.

Since $T\R^{2n}=\R^{2n}\times\R^{2n}$, we can think
of any monomorphism $F\co TL \ra T\R^{2n}$ as the data $(f,\phi)$,
where $f\co L\rightarrow \R^{2n}$ is the underlying base map and
$\phi\co TL\rightarrow \R^{2n}$ is a fibrewise injective
homomorphism. In what follows we will often use $\phi$ to pull
back the standard metric $g$ on $\R^{2n}$, obtaining a metric
$\phi^*g$ on $TL$.

Let us start by showing that totally real monomorphisms are
closely related to the \textit{complexified tangent bundle\/}
$TL^\C$; as usual, we will say that $TL^\C$ is trivial if there
exists a complex vector bundle isomorphism (a \textit{complex
trivialization}) $TL^\C \simeq L\times\C^n$. Any map $\phi$ as
above has a natural \textit{complexification\/}
$\phi^\C\co TL^\C\rightarrow \R^{2n}$, defined by
$$\phi^\C(v+iw):=\phi(v)+J\phi(w)$$
where $J$ is the standard complex structure on $\R^{2n}$.
\begin{lemma}
\label{l:real}
There is a one-to-one correspondence between
totally real injections $TL\rightarrow \R^{2n}$ and
trivializations of $TL^\C$.
\end{lemma}
\begin{proof}
Let $\phi\co TL\rightarrow\R^{2n}$ be a fibrewise linear map. We will
first show that $\phi$ is totally real if and only if
$\phi^\C\co TL^\C\rightarrow \R^{2n}$ is a complex isomorphism (with
respect to $J$). Since $\phi^{\C}$ is a complex linear map between
complex spaces of the same dimension, it suffices to show
surjectivity of $\phi^{\C}$. But $\phi^{\C}$ is surjective at
$p\in L$ if and only if the two real $n$--planes $\phi(T_pL)$ and
$J\phi(T_pL)$ in $\R^{2n}$ intersect transversally, ie
if and only if $\phi$ is a totally real monomorphism. Hence any
totally real injection naturally induces a  complex trivialization
of $TL^{\C}$.

Conversely, any complex trivialization $\phi^{\C}\co  TL^{\C} \ra L
\times \C^n$ determines a unique totally real injection $\phi: TL
\ra \C^n$ defined by $\phi:= \phi^{\C} \circ i$ where $i\co TL \ra
TL^{\C}$ denotes the natural (totally real) inclusion of $TL$ in
$TL^{\C}$.
\end{proof}
\fullref{l:real} indicates the existence of a distinguished subset of trivializations of $TL^\C$: the ones obtained by complexifying Lagrangian injections $TL\ra \C^n$. In \fullref{l:trivs:maps} we will refer to these as the \textit{Lagrangian trivializations\/} of $TL^\C$.

Even though Lagrangian monomorphisms form a proper subset of the
set of all totally real monomorphisms (analogously, the Lagrangian
trivializations are a proper subset of the set of all
trivializations of $TL^\C$), \fullref{c:monos} below will
show that the difference between these two categories disappears
when one works in terms of homotopy classes. This is a direct
consequence of the following result.

\begin{lemma}[Polar decomposition]
\label{l:polar} There exists a strong deformation retraction
$\rho\co \gl{n} \times I \ra \gl{n}$ onto $\unitary{n}$ which is
equivariant with respect to $\unitary{n}$--multiplication on the
left or the right; ie $$\rho_t(AM)=A\rho_t(M), \
\rho_t(MA)=\rho_t(M)A,  \  \text{for all} \, A\in \textrm{U}(n), \ M\in
\gl{n}, \ t\in [0,1].$$
\end{lemma}
\begin{proof} The standard polar decomposition theorem proves the
existence of a unique decomposition $M=PU$ with $U\in
\textrm{U}(n)$ and $P = \sqrt{M M^*} \in \mathcal{P}^+$, where $\mathcal{P}^+$
denotes the set of positive self-adjoint matrices. Note if
$P\in \mathcal{P}^+$ then the matrix $P_t:=t\, \text{Id}+(1-t)P\in\mathcal{P}^+$.
Let us define $\rho_t(M):=P_tU$. Clearly $\rho_0(M)=M$ and
$\rho_1(M)=U$, so $\rho_t$ is a strong retraction as claimed. Now
suppose that $A$ is in $\textrm{U}(n)$. The equivariance of $\rho_t$
with respect to the right action of $\unitary{n}$ is easy to see.
To show the left-equivariance, notice that
$AM=APU=APA^{-1}AU=\tilde{P}AU$, where $\tilde{P}:=APA^{-1}$ is
positive self-adjoint and $AU\in \textrm{U}(n)$. Uniqueness of the
decomposition implies that the polar decomposition of $AM$ is
$\tilde{P}(AU)$. Notice also that $\tilde{P}_t=AP_tA^{-1}$. This
shows that $\rho_t(AM)=\tilde{P}_tAU=AP_tU=A\rho_t(M)$.
\end{proof}
\begin{remark} \label{r:real:polar}
If we restrict $\rho$ to the set of real matrices $\glr{n} \subset
\gl{n}$, we obtain corresponding statements for the \textit{real\/}
polar decomposition theorem: for any $M\in \glr{n}$, $M=PO$ where
$P$ is positive symmetric and $O\in \orth{n}$; furthermore,
$\glr{n}$ retracts $\orth{n}$--equivariantly onto $\orth{n}$.
\end{remark}
\begin{corollary}
\label{c:monos} A totally real monomorphism $F\co TL \ra T\R^{2n}$
is homotopic through a path of totally real monomorphisms to a
Lagrangian monomorphism $\wtilde{F}\co  TL \ra T\R^{2n}$.
Moreover, two totally real monomorphisms $F_1, F_2\co  TL \ra
T\R^{2n}$ are homotopic via totally real monomorphisms if and only
if the corresponding Lagrangian monomorphisms $\wtilde{F}_1, \wtilde{F}_2$ are
homotopic through Lagrangian monomorphisms.
\end{corollary}
\begin{proof}
The proof is basically a direct consequence of the \fullref{l:polar}.
For completeness, we provide here the details.

Let $F=(f,\phi)$. For any $p\in L$, let $\{v_i\}_{i=1}^n$ be an
orthonormal basis of $T_pL$ (with respect to the pullback metric
on $L$ induced by $F$). Let $M=M(p)$ denote the matrix expressing
the homomorphism $\phi^\C\co T_pL^\C\ra \C^n$ in terms of the basis
$\{v_i\}$ of $T_pL^\C$ and the standard basis of $\C^n$.
Since $F$ is totally real, $M\in \gl{n}$. Define $M_t:= \rho_t(M)$
where $\rho$ is the retraction defined in the
\fullref{l:polar}. $M_t$ gives a homotopy through $\gl{n}$
matrixes between $M$ and the matrix $L:=\rho_1(M) \in
\unitary{n}$. We now define $\phi_t\co T_pL\ra\C^n$ to be the
homomorphism corresponding to the matrix $M_t$, and hence obtain a
one-parameter family of monomorphisms $F_t:=(f,{\phi_t})$. Since
$M_t \in \gl{n}$ and $L\in \unitary{n}$, $F_t$ (for $t\in [0,1$))
are totally real monomorphisms and $\tilde{F}:=F_1$  is a
Lagrangian monomorphism (furthermore, the pullback metrics induced
by $F,F_1$ coincide). Notice that if we begin with a different
orthonormal basis $\{w_1,\dots, w_n\}$ then $w_i=v_ja_{ji}$ for some
$A=(a_{ij})\in \textrm{O}(n)$; now $\phi$ gives rise to the matrix $MA$ and hence,
by the right-equivariance of $\rho$, $\phi_1$ corresponds to $LA$.
Hence the construction is independent of the choice of basis of
$T_pL$, so that $\wtilde{F}$ is well-defined. Since all this
also works parametrically, it also proves the second statement.
\end{proof}
\begin{remark}
\label{r:retract} The same methods apply to show that polar
decomposition induces a natural retraction from
$\textrm{Gr}_{\textrm{real}}(\R^{2n})$ to
$\textrm{Gr}_{\textrm{lag}}(\R^{2n})$. Using the
$\unitary{n}$--invariance of polar decomposition mentioned in \fullref{l:polar} these results can be extended to totally real
monomorphisms from $TL$ into any symplectic ambient manifold
$(W,\omega)$ endowed with a compatible almost-complex structure
$J$ and metric $g$. The proof is as above: we choose an
orthonormal basis $v_i$ for $T_pL$ and a unitary basis $e_i$ for
$T_{f(p)}W$ and reduce the question to matrices, then prove that
the result is independent of the choice of both sets of bases
using the left and right equivariance of polar decomposition.
\end{remark}
Hence we have the following result:
\begin{corollary}
\label{c:hprinc:reform}
Let $L$ be a closed smooth (not necessarily orientable) $n$--manifold. Then the following conditions are equivalent.
\begin{enumerate}
\item $L$ admits a Lagrangian immersion $f\co L\rightarrow \R^{2n}$.
\item $TL^\C$ is trivial.
\item $L$ admits a totally real monomorphism $F\co TL\rightarrow T\R^{2n}$.
\item $L$ admits a Lagrangian monomorphism $F\co TL\rightarrow T\R^{2n}$.
\end{enumerate}
Furthermore, there is a one-to-one correspondence between the homotopy classes of the above objects.
\end{corollary}
\begin{proof}
Notice that there is an obvious map from the set $\mathcal{I}$ of
regular homotopy classes of Lagrangian immersions to the set
$\mathcal{M}$ of homotopy classes of Lagrangian monomorphisms,
namely,
$$\mathcal{I}\ra\mathcal{M},\quad [f]\mapsto[df].$$
\fullref{T:h:princ:closed} shows that this map is surjective;
this also proves the equivalence between statements (1) and (4).
\fullref{T:h:princ:parametric} shows that this map is
injective. \fullref{c:monos} proves the equivalence of (3) and
(4) and of their homotopy classes. To prove the equivalence between
(2) and (3), let $(f,\phi)$ be a totally real monomorphism. \fullref{l:real} then proves that $\phi$ defines a complex
trivialization of $TL^\C$. Conversely,  given $\phi$ we can use
any smooth map $f\co L\rightarrow \R^{2n}$ (eg a constant
map) to obtain a totally real monomorphism $(f,\phi)$. Since
$\R^{2n}$ is topologically trivial all such $f$ are homotopic, so
the correspondence between homotopy classes is simple.
\end{proof}
\fullref{c:hprinc:reform} part (2) thus completely solves
the ``Lagrangian immersion in $\R^{2n}\,$ problem" for closed
manifolds in terms of a topological condition on $L$. There is one
last thing we can do, which is to classify these homotopy classes.
This is based on the following simple fact.
\begin{lemma}
\label{l:trivs:maps} Any given trivialization $\phi_0$ of $TL^{\C}$ induces
a bijective correspondence between the set of all trivializations
of $TL^{\C}$ and the set of all maps $L\rightarrow \gl{n}$. Two trivializations are homotopic if and only if the
corresponding maps are homotopic.

Furthermore we have the following:
\begin{enumerate}
\item
Let $h$ and $g$ denote the standard Hermitian and Riemannian metrics on $\C^n$. Then the above correspondence, restricted to the set of maps $L\ra U(n)$, yields the set of all other trivializations whose pullback Hermitian metric $\phi^*h$ on $TL^\C$ coincides with $\phi_0^*h$.
\item Now assume $\phi_0$ is a Lagrangian trivialization of $TL^\C$ (as defined above). Then the above correspondence, restricted to the set of maps $L\ra U(n)$, yields the set of all other Lagrangian trivializations whose pullback Riemannian metric $\phi^*g$ on $TL$ coincides with $\phi_0^*g$.
\end{enumerate}
\end{lemma}
\begin{proof}
For $i=0, \, 1$ let $\phi_i\co  TL^{\C} \ra L \times \C^n$ be complex
trivializations of $TL^{\C}$. Then $\phi_1 \circ \phi_0^{-1}:
L\times \C^n \ra L \times \C^n$ is a fibrewise complex
isomorphism, or equivalently, a map from $L$ into
$\textrm{Aut}(\C^n)$. We identify $\textrm{Aut}(\C^n)$ with
$\gl{n}$ by choosing the standard basis $e_1, \ldots ,e_n$ of
$\C^n$. We write the map $L\ra \gl{n}$ determined by $\phi_0$ and
$\phi_1$ (and the choice of our standard basis for $\C^n$) as
$M_{\phi_1,\phi_0}$. Hence if we choose a reference trivialization
$\phi_0$ of $TL^{\C}$, then any other trivialization $\phi$ of
$TL^{\C}$ determines a map $L \ra \gl{n}$, namely
$M_{\phi,\phi_0}$. Conversely, a map $M\co  L \ra \gl{n}$ determines
a fibrewise complex-linear automorphism $A(M)$ of $L \times \C^n$,
and hence another complex trivialization of $TL^{\C}$, denoted
$\phi_0^M$ and given by $\phi_0^M:= A(M) \circ \phi_0 $.

Now let $v_i:=\phi_0^{-1}e_i$. Notice that $\phi^*h=\phi_0^*h$ if
and only if $\phi^*h(v_i,v_j)=\phi_0^*h(v_i,v_j)$. However, the
left hand side of this last expression is
$\smash{(\phi\circ\phi_0^{-1})^*h(e_i,e_j)}$ while the right hand side is
$h(e_i,e_j)$, so we conclude that $\phi^*h=\phi_0^*h$ if and only
if $\phi\circ\phi_0^{-1}$ is unitary, ie
$M_{\phi,\phi_0}\in \unitary{n}$. This proves (1). Now recall that $h=g+i\omega$, so if $\phi_0$ is Lagrangian then
$\phi_0^*h=\phi_0^*g$. In this case $\phi^*h=\phi_0^*h$ if and
only if $\phi^*h=\phi_0^*g$; in particular this implies that
$\phi^*\omega=0$, ie that $\phi(TL)$ is Lagrangian in
$\C^n$, and that $\phi^*g=\phi_0^*g$. This proves (2).
\end{proof}
We can now classify the homotopy classes of Lagrangian immersions as follows (cf \cite[page 274]{audin}).
\begin{prop}
\label{p:class:lagr:imm} Given a Lagrangian immersion
$L\rightarrow \C^n$, the set of regular homotopy classes of all
Lagrangian immersions $L \ra \C^n$ (or of homotopy classes of Lagrangian monomorphisms) is in bijective correspondence
with the set $[L;\unitary{n}]$ of homotopy classes of maps from $L$
into $\unitary{n}$ (or equivalently into $\gl{n}$).
\end{prop}
\begin{proof}
\fullref{c:hprinc:reform} allows us to rephrase the
statement in terms of Lagrangian monomorphisms; the proof of
\fullref{c:hprinc:reform} shows that all base maps are
homotopic, allowing us to further reduce to Lagrangian
trivializations of $TL^\C$. We can now prove the statement as
follows. The given Lagrangian immersion defines a reference
Lagrangian trivialization $\phi_0$ and a pullback metric
$\phi_0^*g$ on $TL$. Choose another Lagrangian trivialization
$\phi$; for any point $p\in L$, let $\pi_p$ denote the image
Lagrangian plane in $\C^n$, with the metric induced from $g$. Then
$\phi(p)\co T_pL\ra \pi_p$ is an isomorphism, and by the real polar
decomposition lemma (see \fullref{r:real:polar}) it can be
homotoped to an isometry $\tilde{\vphantom{i}\smash{\phi}}$ (as in \fullref{c:monos}) without changing the image plane $\pi_p$. Since
$\tilde{\vphantom{i}\smash{\phi}}$ has the same image $\pi_p$, $\tilde{\vphantom{i}\smash{\phi}}$  still
defines a Lagrangian trivialization of $TL^\C$; thus part (2) of
\fullref{l:trivs:maps} shows that $\tilde{\vphantom{i}\smash{\phi}}$ corresponds to
a map $L\ra \unitary{n}$.
\end{proof}
\subsection{Gauss maps and the Maslov class}\label{s:maslov:class}
To simplify the discussion, in this section we will assume
that $L$ is oriented; at the end of the section we will add a few comments on the nonorientable case.

Since $\R^{2n}$ is contractible, $\smash{\textrm{Gr}^+_{\textrm{lag}}}(\R^{2n})$
is a trivial bundle over $\R^{2n}$ which has fibre
$\textrm{U}(n)/\textrm{SO}(n)$. Hence given any Lagrangian
monomorphism $F\co TL \ra T\R^{2n}$ the tangential Gauss map $GF$
gives us a map
$$GF\co L \ra \textrm{U}(n)/\textrm{SO}(n).$$ More explicitly, fix $x\in L$ and consider the
homomorphism $F_x\co T_xL\ra \C^n$. We can use any positive
orthonormal basis $\{v_i\}_{i=1}^n$ (with respect to the pullback
metric) of $T_xL$ and the standard basis of $\C^n$ to represent
$F_x$ via a matrix in $\unitary{n}$. This matrix depends on the
choice of $\{v_i\}$ but its equivalence class in
$\unitary{n}/\sorth{n}$ does not; this is the Gauss map $GF(x)$.
Composition with the complex determinant yields an associated
\textit{determinant map\/}
$$\textrm{det}_{\C} \, \circ \,GF\co  L \ra\Sph^1.$$
This construction induces a well-defined map $\mu$ on the set of
homotopy classes $\mathcal{M}$ of Lagrangian monomorphisms, which
we call the \textit{Maslov map\/} $\mu\co \mathcal{M}\ra [L;\Sph^1]$
and is given by
$$\mu_{F}:=[\textrm{det}_{\C} \circ GF].$$
When the Lagrangian monomorphism $F$ arises from an oriented Lagrangian immersion
then the determinant map defined above coincides with the Lagrangian phase
function $e^{i\theta}$ defined in \fullref{SS:slg}.

We now recall some well-known facts about homotopy classes of maps
to $\Sph^1$. Since $\Sph^1$ is an abelian Lie group the homotopy
classes of maps into $\Sph^1$ inherit the structure of an abelian
group. We will let $d\theta$ denote the standard harmonic $1$--form
on $\Sph^1$, normalized so that $[d\theta]$ is the generator of
$H^1(\Sph^1,\Z)\simeq\Z$.
\begin{lemma}
\label{l:homotopy:s1}
 For any CW complex $K$ the map 
$[K;\Sph^1]\ra H^1(K,\Z)$ induced by $$f\mapsto f^*[d\theta]$$
is an isomorphism of abelian groups. Also, there is an isomorphism
between $H^1(K,\!\Z)$ and $\Hom(\pi_1(K),\Z)$. In particular, a map
$f\co K\ra\Sph^1$ is homotopic to a constant if and only if
$f|_\gamma$ is homotopic to a constant, for all closed loops
$\gamma\subset K$.
\end{lemma}
\begin{proof}
If $K$ is any CW complex and $Y$ is an Eilenberg--Mac Lane space of
type $K(\pi,n)$, there is a natural isomorphism between the
homotopy classes of maps $[K;Y]$ and $H^n(K;\pi)$ \cite[Theorem
VII.12.1]{bredon}. The first part now follows since $\Sph^1$ is an
Eilenberg--Mac Lane space $K(\Z,1)$. From the Universal Coefficient
Theorem we have $H^1(K,\Z) \simeq \Hom(H_1(K),\Z)
+\Ext(H_0(K),\Z)$ \cite[page 282]{bredon}. Since $H_0(K)$ is always
free, we have $H^1(K,\Z) \simeq \Hom(H_1(K),\Z)
\simeq\Hom(\pi_1(K),\Z)$ since $\Z$ is abelian.
\end{proof}
Using \fullref{l:homotopy:s1} we can identify $\mu_F$ with the
element $(\textrm{det}_{\C} \circ GF)^*[d\theta]\in H^1(L,\Z)$;
this is the \textit{Maslov class\/} of $F$. Furthermore, $\mu_F$ is
equivalent to a homomorphism $m_F\co  \pi_1(L)\ra\Z$. Given a loop
$\gamma\co \Sph^1\ra L$, $m_F([\gamma])\in \Z$ simply calculates the
degree of the map $\textrm{det}_\C\circ GF\circ\gamma\co \Sph^1\ra
\Sph^1$; we call this integer the \textit{Maslov index\/} of the
curve $\gamma$. Notice that according to \fullref{c:hprinc:reform} we can always represent $\mu_{F}$ by a
Lagrangian immersion $f$ in the same homotopy class; we will then
use the notation $\mu_f$.
We now want to describe a second way to compute $\mu_F$, along the
lines of \fullref{p:class:lagr:imm}. Let us fix a
Lagrangian monomorphism $F_0=(f_0,\phi_0)$. As explained in \fullref{p:class:lagr:imm}, any other Lagrangian
monomorphism $F=(f,\phi)$ can be homotoped so that $f=f_0$ (see the proof of \fullref{c:hprinc:reform}) and the induced pullback metrics $\phi^*g$, $\phi_0^*g$ on $TL$ coincide, so $F$ will correspond to a map $M_{F,F_0}\co L\ra \unitary{n}$.
\begin{lemma} \label{l:matrix:maslov}
In the situation just described, $\mu_F=[\textrm{det}_\C
M_{F,F_0}]\cdot \mu_{F_0}$, where $\cdot$ denotes group
multiplication in the abelian group $[L;\Sph^1]$. In particular,
if $F_0$ has zero Maslov class then $\mu_F=[\textrm{det}_\C
M_{F,F_0}]$.
\end{lemma}
\begin{proof}
Let $F_0=(f,\phi_0)$ and $F=(f,\phi)$. Fix $p\in L$; let
$\{v_i\}_{i=1}^n$ be an orthonormal basis of $T_pL$ with respect
to $\phi_0^*g$. Let $M(\phi_0), M(\phi)$ be the matrices
expressing $\smash{\phi^\C_0}$ and $\smash{\phi^\C}$ with respect to the basis
$\{v_i\}$ of $T_pL^\C$ and the standard basis of $\C^n$. Notice
that $\phi^\C=\phi^\C\circ(\phi_0^\C)^{-1}\circ\phi_0^\C$.
 Thus $\textrm{det}_\C\circ GF=\textrm{det}_\C\,M(\phi)=\textrm{det}_\C\,M_{F,F_0}\cdot (\textrm{det}_\C\circ GF_0)$ so
 $\mu_F=[\textrm{det}_\C\,M_{F,F_0}]\cdot \mu_{F_0}$.
\end{proof}
\begin{prop}\label{p:anymaslov:ok}
Let $L$ be a closed orientable $n$--manifold $L$ with $TL^\C$ trivial, and let $\mu$
be an arbitrary element of $H^1(L,\Z)$.
\begin{enumerate}
\item $L$ admits a Lagrangian monomorphism $F\co TL\ra T\C^n$ with $\mu_F=\mu$.
\item $L$ admits an exact Lagrangian immersion $f\co L \ra \R^{2n}$ with $\mu_f=\mu$.
\end{enumerate}
\end{prop}
\begin{proof}
Since $TL^\C$ is trivial there exists at least one Lagrangian
monomorphism $F_0\co TL\rightarrow T\R^{2n}$ which we can fix as our
reference monomorphism. As in \fullref{l:trivs:maps}, all other
such monomorphisms are determined (up to homotopy) by maps from $L$ to 
$\unitary{n}$ and the corresponding Maslov classes can be
calculated as in \fullref{l:matrix:maslov}. The inclusion
$\Sph^1 \hookrightarrow \unitary{n}$, given by
$$e^{i\theta} \mapsto \diag{(e^{i\theta},1,\ldots,1)}$$ shows that any map $L\rightarrow \Sph^1$ is
the determinant of some map $L\rightarrow \unitary{n}$; in other
words, the map
$$[L,\unitary{n}]\rightarrow[L,\Sph^1]\simeq H^1(L,\Z)$$
induced by $\textrm{det}_\C$ is surjective. Thus for any $\mu\in
H^1(L,\Z)$ we can find $M\co L\rightarrow \unitary{n}$ such that
$\mu=[\textrm{det}_\C M]\cdot\mu_{F_0}$. Let $F$ be the Lagrangian
monomorphism determined by $F_0$ and $M$. It follows from  \fullref{l:matrix:maslov} that $\mu_F=\mu$ as required.

The second statement is a direct consequence of the Lagrangian
$h$--principle.
\end{proof}

It follows from \fullref{p:anymaslov:ok} that in the orientable case there are no
extra obstructions to finding Maslov-zero monomorphisms (or immersions).

It is an interesting issue how one can use information on the Gauss map to study a given Lagrangian monomorphism. To investigate this, let us introduce a new way of thinking about Lagrangian trivializations. Let $\phi\co TL^\C\ra L \times \C^n$ be a Lagrangian trivialization. Let $SO(L)$ denote the bundle of positive orthonormal frames on $L$ induced by the pullback metric. Recall that $SO(L)$ is a $\sorth{n}$--principal fibre bundle over $L$ with respect to the right action of $\sorth{n}$ defined as follows: if $(v_1,\dots,v_n)\in SO(L)$ and $A\in \sorth{n}$, then
$(v_1,\dots,v_n)\cdot A:=(w_1,\dots,w_n)$ where $w_i:=v_j a_{ji}$. We will denote by $p\co SO(L)\ra L$ the obvious projection.

There is a natural map
$$\Phi\co SO(L)\ra \unitary{n},\quad (v_1,\dots,v_n)\in p^{-1}(x)\mapsto \Phi(v_1,\dots,v_n)$$
where $\Phi(v_1,\dots,v_n)$ is the matrix which represents $\phi_x$ with respect to the basis $(v_1,\dots,v_n)$ and the standard basis of $\C^n$. This map is $\sorth{n}$--equivariant with respect to the action of $\sorth{n}$ on $\unitary{n}$ determined by right multiplication. Denote the natural projection by $q\co \unitary{n}\ra \textrm{U}(n)/\textrm{SO}(n)$. Notice that the Gauss map satisfies the relation $p\circ G\phi=\Phi\circ q$. In other words, $\Phi$ is a map between the principal fibre bundles $SO(L)$ and $\unitary{n}$ which covers the corresponding Gauss map on the base spaces $L$ and $\unitary{n}/\sorth{n}$.

More generally, any Lagrangian trivialization inducing the same metric defines an equivariant map $SO(L)\ra \unitary{n}$. Furthermore, it is clear that this procedure defines a one-to-one correspondence between all such trivializations and all such maps. This correspondence respects homotopy classes in the following sense: two such trivializations are homotopic if and only if the corresponding maps are homotopic.

Notice that a given map $L\ra \textrm{U}(n)/\textrm{SO}(n)$ is a Gauss map if and only if it admits a lift to an equivariant map $SO(L)\ra \unitary{n}$. It follows from the general theory of principal fibre bundles that the existence of such lifts is a homotopically invariant property: if $g\co L\ra \textrm{U}(n)/\textrm{SO}(n)$ is a Gauss map, ie $g=G\phi$, and $g$ is homotopic to $g'$ then $g'$ is the Gauss map of some $\phi'$ and the corresponding $\Phi,\Phi'$ are themselves homotopic through equivariant maps $SO(L)\ra \unitary{n}$. We can use this as follows.

\begin{definition}
Let $\textrm{Gr}_{\textrm{SL}}(\C^n)$ denote the Grassmannian of special Lagrangian planes in $\C^n$; it is a trivial subbundle of $\smash{\textrm{Gr}^+_{\textrm{lag}}(\R^{2n})}$, with fibre $\sunitary{n}/\sorth{n}$.
A Lagrangian monomorphism $F$ is \textit{special Lagrangian\/} if its Gauss map takes values in $\textrm{Gr}_{\textrm{SL}}(\C^n)$.
\end{definition}
SL monomorphisms are the obvious formal analogue of SL immersions.
\begin{prop}\label{p:MZisSL}
Let $L$ be an oriented manifold. Then any Maslov-zero monomorphism $F$ is homotopic through Lagrangian monomorphisms to an SL monomorphism $F'$.
\end{prop}
\begin{proof} Write $F=(f,\phi)$ and let $\Phi$ denote the corresponding equivariant map from $SO(L)$ to $\unitary{n}$ (with respect to the induced metric). Recall that $\unitary{n}\simeq \Sph^1\times \sunitary{n}$ and that $\unitary{n}/\sorth{n}\simeq \Sph^1\times \sunitary{n}/\sorth{n}$. Thinking of $G\phi$ as a map from $L$ to $\Sph^1\times \sunitary{n}/\sorth{n}$, we see that $F$ is Maslov-zero if and only if $G\phi$ is homotopic to a map $g'\co L\ra \sunitary{n}/\sorth{n}$. By the homotopic invariance of the lifting property, $g'$ is the Gauss map of some $\Phi'$ homotopic to $\Phi$. By construction, $\Phi'$ is a map $SO(L)\ra \sunitary{n}$; the corresponding monomorphism $F'=(f,\phi')$ is a SL monomorphism homotopic to $F$.
\end{proof}
From the viewpoint of SL geometry, this proposition justifies our interest in Maslov-zero monomorphisms and immersions. It also shows that there are no extra obstructions to finding SL monomorphisms. However there do exist strong obstructions to finding SL immersions. For example, since $\C^n$ contains no closed minimal submanifolds, under our assumptions ($L$ compact) there cannot
exist SL immersions $f\co L\rightarrow \C^n$. In particular this
shows that SL monomorphisms do not satisfy an $h$--principle.

We conclude this section with a few comments on the nonorientable case. The nonoriented Lagrangian Grassmannian has fibre $\unitary{n}/\orth{n}$, so to get a well-defined determinant map it is necessary to replace $\textrm{det}_\C$ with $\textrm{det}_\C^2$. Unoriented analogues of the Maslov map, the Maslov class and the Maslov index of a curve can then be defined as before.
The following lemma shows that if the unoriented Maslov data of a Lagrangian
 monomorphism $F\co  TL^{\C} \ra \C^n$ satisfies certain conditions then in fact $L$ must be orientable.
\begin{lemma}
Suppose that $L$ (not assumed to be orientable) admits a Lagrangian
mono\-morphism $F\co TL^\C\ra \C^n$ such that the (unoriented) Maslov index
$$m_F(\gamma):= \textrm{deg} (\textrm{det}_\C^2 \circ GF \circ \gamma)$$ of every
loop $\gamma$ in $L$ is even. Then $L$ is orientable.
\end{lemma}
\begin{proof}
Consider the two-to-one covering map $p\co \unitary{n}/\sorth{n} \ra \unitary{n}/\orth{n}$ obtained by forgetting the orientation of an oriented Lagrangian $n$--plane
in $\R^{2n}$.
One can show that $p_*(\pi_1(\unitary{n}/\sorth{n})$ is an index
two subgroup of $\pi_1(\unitary{n}/\orth{n})$ isomorphic to $2\Z \subset
\Z$.
By standard covering space theory the map $GF\co L \ra \unitary{n}/\orth{n}$
lifts to the two-fold cover $\unitary{n}/\sorth{n}$ if and only if
$GF_*\pi_1(L) \subseteq p_*\pi_1(\unitary{n}/\sorth{n})$. But this is equivalent
to the fact that the Maslov index of every loop in $L$ is even.
Hence $GF$ lifts to a map $GF^+\co L \ra \unitary{n}/\sorth{n}$, and this
gives an orientation of each tangent space of $L$ as needed.
\end{proof}
An immediate corollary of the previous Lemma is the following:
\begin{corollary} Suppose $L$ admits a Lagrangian monomorphism $F\co TL^\C\ra \C^n$ such that $\textrm{det}_\C^2 GF\co  L\ra \Sph^1$ is homotopic to a constant. Then $L$ is orientable, and the corresponding oriented Gauss map is Maslov-zero.
\end{corollary}
The above shows that \fullref{p:anymaslov:ok} cannot hold for closed nonorientable manifolds $L$; in particular, $0\in H^1(L,\Z)$ cannot be realized as a Maslov class of any nonorientable manifold. The analogous result, which can be proved via the same methods, is as follows.
\begin{prop} \label{p:somemaslov:ok}
Let $L$ be a closed, nonorientable manifold. Assume $TL^\C$ is trivial, so that there exists a Lagrangian monomorphism $F$ with (unoriented) Maslov class $\mu_F\in H^1(L,\Z)$. Then any other Lagrangian monomorphism from $L$
has Maslov class in the set $2H^1(L,\Z)+\mu_F\subset H^1(L,\Z)$, and any such cohomology class can be realized this way.
\end{prop}

\subsection{Examples}
\label{s:closed:examples} We conclude this section by using
\fullref{c:hprinc:reform} to give examples of manifolds
which do and do not admit Lagrangian immersions into $\R^{2n}$.
\subsubsection{Low-dimensional cases} \label{s:closed:low:dim}
For $n\le 3$, one can use
the previous results to give a very good description of which
closed $n$--manifolds admit Lagrangian immersions into $\C^n$ and
to describe the regular homotopy classes of Lagrangian immersions.

$n=1$\qua In this case $L=S^1$ and any immersion of $L$ into $\C$ is
Lagrangian. Let us fix one immersion $f_0\co S^1\rightarrow \C$. As
explained in \fullref{l:trivs:maps}, any other immersion $f$
now defines an element in $[S^1;\unitary{1}] = \pi_1(\unitary{1})
= \Z$, ie an integer. Notice that the Gauss maps take
values in the Grassmannian $\unitary{1}/\textrm{SO}(1) \simeq
\Sph^1$ and that $\textrm{det}_\C=\text{Id}$. Choosing $f_0$ to have
zero-Maslov class is equivalent to making $Gdf_0$ homotopically
trivial; for example, we could choose $f_0$ to be a ``figure eight
curve" inside $\C$. With such a choice of $f_0$ the above integer
is the turning number of $f$ in $\C$. Hence the classification of
regular homotopy classes of Lagrangian immersions of $S^1$ in $\C$
given by \fullref{p:class:lagr:imm} reduces to Whitney's
classification of immersions of $S^1$ in $\R^2$ according to their
turning number.

$n=2$\qua Let $L$ be a closed surface, not necessarily orientable. In
this case one can show that the complexified tangent bundle
$TL^{\C}$ is trivial if and only if the Euler characteristic
$\chi(L)$ is even \cite[page 274]{audin}. Hence by \fullref{c:hprinc:reform} every orientable surface admits Lagrangian
immersions in $\C^2$, while if $L$ is the connected sum of $k$
copies of $\RP^2$ then it admits Lagrangian immersions into $\C^2$
if and only if $k$ is even. For example, $\RP^2$ has no Lagrangian
immersions into $\C^2$, while the Klein bottle $K$ does.

In the oriented case, let us now fix one such immersion. Recall
from the proof of \fullref{p:anymaslov:ok} that the map
$[L;\unitary{2}] \ra [L;\Sph^1]$ induced by $\textrm{det}_\C$ is
surjective. Since the map $\det_{\C}\co  \unitary{n} \ra \Sph^1$
induces isomorphisms on $\pi_1$ and $\pi_2$, it follows using
techniques similar to those of \fullref{S:obstruction} that in
dimension 2 this map is also injective  \cite[Corollary
VII.11.13]{bredon}. Hence the map $M\mapsto
(\det_{\C}M)^*[d\theta]$ induces a bijection from
$[L;\unitary{2}]$ to $H^1(L,\Z)$. In other words, for any closed
oriented surface we can strengthen \fullref{p:anymaslov:ok} as follows: the regular homotopy class of the
immersion $f$ is completely determined by $\mu$.

$n=3$\qua We shall assume that $L$ is a closed orientable $3$--manifold.
In this case, it is a theorem originally stated by Stiefel that $L$
is parallelizable \cite[page 148]{milnor:stasheff}. Complexification
gives an obvious complex parallelization of $L$. Hence any closed
orientable $3$--manifold admits exact Lagrangian immersions into $\C^3$.
In the special case where $L=S^3$ then the regular homotopy classes of
Lagrangian immersions are in one-to-one correspondence with
$\pi_3(\unitary{3}) = \Z$. In particular, unlike the case of curves and surfaces,
 the Maslov class of a Lagrangian immersion $f$ of $L^3$ no longer determines its regular homotopy class.
However, one can prove that there is a bijection
between $[L;\unitary{3}]$ and $H^1(L,\Z) {\times} H^3(L,\Z)$
\cite[Proposition 1]{auckly}.

On the other hand, it is easy to find nonorientable closed
$3$--manifolds which do not admit Lagrangian immersions;
$\RP^2\times S^1$ is probably the simplest example.
\subsubsection{Manifolds with trivial or stably trivial tangent
bundles} \label{s:stably:trivial} In the previous example we saw
that if there exists a trivialization of $TL$ then complexification gives an obvious
trivialization of $TL^{\C}$. Hence any closed parallelizable
$n$--manifold admits Lagrangian immersions into $\C^n$. In
particular any compact Lie group admits Lagrangian immersions into
$\C^n$ and so do the spheres and real projective spaces in
dimensions $1$, $3$ and $7$ (which are well-known to be the only
parallelizable spheres or projective spaces).
More generally, one can show that if $TL$ is
only stably parallelizable (ie the direct sum of $TL$
with some trivial bundle is itself trivial)
then $TL^\C$ is still trivial \cite[page 273]{audin}.
As we have already mentioned, any embedded hypersurface in $\R^{n+1}$ has this property.
This gives an alternative proof that any compact orientable surface
admits Lagrangian immersions into $\R^4$.
It also proves that $S^n$ admits
Lagrangian immersions into $\R^{2n}$ for any $n$
(on the other hand we will
see that for most $n$, $\RP^n$ does not admit Lagrangian immersions).
Furthermore, the regular homotopy classes of Lagrangian immersions are in
one-to-one correspondence with $[S^n,\unitary{n}] =\pi_n(\unitary{n})$, which equals $0$ or $\Z$ depending on
whether $n$ is even or odd respectively.

There are many more stably parallelizable manifolds than
parallelizable manifolds. For example, the class of stably
parallelizable manifolds is closed under taking products and
connected sums \cite[page 187]{kosinski} and also contains
all manifolds which are homotopy spheres \cite[page 191]{kosinski}.
\subsubsection{Nonexistence of Lagrangian immersions via
characteristic classes} \label{s:nonexist} One can sometimes give
simple proofs of Lagrangian nonimmersion results by
characteristic class arguments. Below we give representative
examples of this type of argument using the Stiefel--Whitney and
Pontrjagin classes of $TL$.

Given any real vector bundle $E^r\rightarrow L$ over any closed
manifold $L^n$,  one can define Stiefel--Whitney classes $w_i(E)\in
H^i(L;\Z/2)$, with $w_0(E)=1$ \cite[\S 4]{milnor:stasheff}.
$w(E):=w_0(E)+\dots+w_r(E)$ is called the total Stiefel--Whitney
class of the bundle. Recall that $w(E)$ has the following
fundamental properties: $w(E\oplus F)=w(E)\cdot w(F)$ and, if $E$
is trivial, then $w(E)=1$. The total Stiefel--Whitney class of $T
\RP^n$ is well-known  to be
\begin{equation}
\label{e:sw:rpn}
w(\RP^n) = (1+a)^{n+1},
\end{equation}
where $a$ denotes the nonzero element of $H^1(\RP^n,\Z/2)$ which
under cup product generates the full $\Z/2$ cohomology ring of
$\RP^n$ \cite[Theorem 4.5]{milnor:stasheff}.
As unoriented real vector bundles $TL^{\C} \simeq TL \oplus TL$.
Hence if $TL^{\C}$ is trivial, then $w(TL)^2 = 1$. Using
\eqref{e:sw:rpn} it is straightforward to show that $w(\RP^n)^2
=1$ if and only if $n+1$ is a power of $2$. Hence there are no
Lagrangian immersions of $\RP^n$ in $\C^n$ unless $n=2^k -1$ for
some $k\in \Z$. For $k=1, 2, 3$ we already saw in \fullref{s:stably:trivial} that such Lagrangian immersions do exist.

The Pontrjagin classes of a real vector bundle $E^r \ra L$ may be
obtained from the Chern classes of its complexification
$E\otimes_{\R} \C$. The standard definition of the $i$--th
Pontrjagin class is $p_i(E): = (-1)^i c_{2i}(E\otimes \C) \in
H^{4i}(L,\Z)$ \cite[page 174]{milnor:stasheff}. The total Pontrjagin
class $p(E)$ is defined to be
\begin{equation}
\label{e:pont} p(E) = 1 + p_1(E) + \ldots + p_{[n/2]}(E) \in
H^*(L).
\end{equation}
The class $p(E)$ equals the total Chern class of $E\otimes \C$, ignoring the
odd Chern classes $c_{2i+1}(E\otimes \C)$ which vanish in de Rham
cohomology and are order $2$ in the integral cohomology
(Bott and Tu give a different definition of $p(E)$
which does not ignore these odd Chern classes \cite[page
289]{bott:tu}).
If the complex vector bundle $E\otimes \C$ is trivial, then from
the basic properties of Chern classes we have $p(E)=1$. In
particular, if $L$ admits a Lagrangian immersion in $\R^{2n}$ then
$TL^{\C}$ is trivial and hence $p(TL)=1$. For instance, the total
Pontrjagin class of $\CP^n$ is well-known \cite[page
177]{milnor:stasheff} to be
\begin{equation}
\label{e:pont:cpn} p(\CP^n) = (1+a^2)^{n+1} = 1 + (n+1)a^2 +
\ldots \end{equation} where $a\in H^2(\CP^n;\Z)$ generates
$H^*(\CP^n;\Z)$ under cup product. In particular, $p(\CP^n) \neq
1$ for $n>1$, and hence admits no Lagrangian immersions in
$\C^{2n}$. For an oriented $4$--manifold $M$ there is only one
nontrivial Pontrjagin class $p_1(M) \in H^4(M,\Z)$ and hence only
one Pontrjagin number $p_1[M]:= \langle p_1(TM),[M]\rangle \in \Z$.
If $\Sigma_d$ is a nonsingular algebraic hypersurface of degree
$d>0$ in $\CP^3$, then a standard characteristic class computation
(see Donaldson and Kronheimer \cite[\S 1.1.7]{don:kron}) shows that the Pontrjagin number
of $\Sigma_d$ is given by
\begin{equation}
\label{e:pont:hyp:d} p_1[\Sigma_d]= (4-d^2)d.
\end{equation}
In particular, $\Sigma_d$ does not admit Lagrangian immersions
into $\C^4$ unless $d=2$, in which case $\Sigma_2 \simeq S^2
\times S^2$ which is stably parallelizable and hence admits a
Lagrangian immersion.

\section{Prescribed Boundary Problem}\label{s:bdy}
We now turn to Lagrangian immersions of compact manifolds with boundary. Specifically, we will prescribe an immersion along the boundary $\Sigma^{n-1}$ and try to find a ``Lagrangian
filling" of this data. In \fullref{ss:pbp} we will describe this ``Prescribed Boundary Problem" more precisely. In \fullref{ss:lagr:cobordism} we begin with a preliminary, related question. Throughout this section, our focus will be on the orientable case.

\subsection{Preliminary considerations}\label{ss:lagr:cobordism}
Let $\Sigma^{n-1}$ be a compact oriented (not necessarily
connected) manifold without boundary, and $i\co \Sigma\rightarrow
\R^{2n}$ be an immersion. Consider the following question.

\medskip
\textbf{Lagrangian Cobordism Problem}\qua Does there exist
a compact oriented $n$--manifold $L$ bounding $\Sigma$ and
a Lagrangian immersion $f\co L\rightarrow \R^{2n}$ extending $i\co \Sigma \ra \R^{2n}$?
\medskip

Our main interest is actually in a stronger version of this question, which we call the Prescribed Boundary Problem. Our goal now is thus not to solve the Lagrangian
Cobordism Problem, but simply to highlight some of the restrictions that it imposes on $\Sigma$ and $i$.

\medskip
\textbf{Topological restrictions}\qua The first question is
whether $\Sigma$ bounds any compact oriented smooth manifold $L$.
As we have already discussed in \fullref{SS:non:smooth} this is equivalent to the vanishing of
all the Pontrjagin and Stiefel--Whitney numbers of $\Sigma$.
\medskip

The second topological issue is that $L$ must admit Lagrangian
immersions. As in \fullref{S:h:princ:closed}, it is easy to
prove that this implies that $TL^\C$ must be
trivial. The existence of such an $L$ is an additional constraint on $\Sigma$.

\medskip
\textbf{Geometric restrictions}\qua Suppose that we have
overcome the topological restrictions described above,
ie $\Sigma$ bounds some compact $n$--manifold $L$ with
$TL^{\C}$ trivial. If there exists a Lagrangian immersion $f: L^n
\ra \R^{2n}$ which extends $i\co \Sigma \ra \R^{2n}$ then $i$ must be
isotropic, ie $i^*\omega=0$. Thus, from now on we
assume that $i$ is an isotropic immersion.
\medskip

A second necessary condition on $i$ comes from cohomological
considerations. First consider the following:
\begin{example} \label{e:1d:exact}
For dimensional reasons, any immersion $i\co S^1\ra\R^4$ is
isotropic. Assume $i$ admits an oriented Lagrangian extension
$f\co L\ra\R^4$. Then by Stokes' Theorem,
$$\int_{S^1}i^*\lambda=\int_{S^1}f^*\lambda=\int_L f^*\omega=0,$$
so $i$ is exact. For instance, this shows the
isotropic immersion $i\co  \Sph^1\ra\C^2$ given by
$$i(e^{i\theta})=(e^{i\theta},0)$$
does not admit any oriented Lagrangian filling.
\end{example}

More generally, suppose $i$ admits an oriented Lagrangian extension
$f\co L \ra \R^{2n}$. Then, for any closed curve
$\gamma\subset\Sigma$ which is homologically trivial in $L$,
Stokes' Theorem shows that $\int_\gamma\lambda=0$. Another way of
saying this is as follows. If $\gamma\subset\Sigma$ is
homologically trivial in $L$, there exists $\alpha\in
H_2(L,\Sigma)$ such that  $\partial\alpha=\gamma$. On the other
hand, since $i$ is isotropic any extension $f$ defines an element
$[f^*\omega]$ in the relative de Rham cohomology $H^2(L,\Sigma)$.
If $f$ is Lagrangian then $[f^*\omega]=0$. By duality this occurs
if and only if $\langle f^*\omega,\alpha \rangle = 0$ for all
$\alpha \in H_2(L,\Sigma)$. Now, using the particular $\alpha$
defined above, we see that
$$\int_\gamma\lambda=\int_\alpha f^*\omega=\langle f^*\omega,\alpha \rangle = 0.$$
This condition is not automatic from the fact that $i$ is an
isotropic immersion.

However if we assume that the immersion $i$ is not merely
isotropic but \textit{exact\/}, then we have $\int_\gamma\lambda=0$
for any closed curve $\gamma$ in $\Sigma$. \fullref{e:1d:exact} shows that when $n=2$ (and $\Sigma$ is connected) exactness is actually a
necessary condition. Moreover, we will see that (as in the closed
case) Lagrangian submanifolds with boundary produced using the
relative version of the Lagrangian $h$--principle are automatically
exact and hence so are their boundaries. Thus, if we want to make
use of the relative $h$--principle then the exactness of $\Sigma$
is a necessary assumption and not a mere convenience. Finally,
when we come to apply our results on the Prescribed Boundary
Problem to Lagrangian desingularization problems we will find that
the local geometry near the singular points often forces $\Sigma$
to be exact and not just isotropic.

Hence from now on we make the following assumptions about $\Sigma$.

\begin{assump}\setobjecttype{Ass}~
\label{assumptions}

\begin{itemize}
\item[A.] $\Sigma$ bounds a compact oriented manifold $L$.
\item[B.] $TL^{\C}$ is trivial.
\item[C.] The immersion $i\co \Sigma \ra \R^{2n}$ is exact.
\end{itemize}
\end{assump}

\subsection{Lagrangian thickenings and the Prescribed Boundary Problem}\label{ss:pbp}
In \fullref{s:desing} we investigate the existence of Lagrangian
smoothings of a singular Lagrangian object with only isolated
singular points. We attempt to find our Lagrangian smoothings by
removing a small neighbourhood $U$ of a singular point and gluing
in a smooth Lagrangian immersion of some compact manifold $L$ with
boundary $\Sigma:=\partial U$. To ensure smoothness of the new
submanifold along the boundary $\partial U$ we need to modify the
Lagrangian Cobordism Problem by assigning as initial data, instead
of just $\Sigma$, a neighbourhood of $\Sigma$ to be ``filled" by
$L$.

Following \cite{eliashberg,gromov:pdr} let us introduce the
notation $\mathcal{O}\!p\,\Sigma$ to denote some open neighbourhood
of $\Sigma\subset L$, which can be varied as appropriate by
restriction to a smaller neighbourhood of $\Sigma$. In our
applications we will usually choose $\mathcal{O}\!p\,\Sigma$ to be
topologically $[0,\epsilon)\times\Sigma$.

\begin{definition}\label{d:initial:data}
Let $L^n$ be a compact connected oriented manifold with oriented
boundary $\Sigma$, such that $TL^\C$ is trivial. Let
$f\co \mathcal{O}\!p\,\Sigma\rightarrow \R^{2n}$ be an exact Lagrangian
immersion of an open neighbourhood of $\Sigma\subset L$ into
$\R^{2n}$. We say that the triple $(\Sigma, L, f)$ is
\textit{(exact) initial data\/} for the Prescribed Boundary Problem.
A \textit{solution to the Prescribed Boundary Problem with initial
data $(\Sigma, L, f)$} is any exact Lagrangian immersion $\hat{f}:
L\rightarrow \R^{2n}$ which agrees with $f$ on some (possibly
smaller) open neighbourhood of $\Sigma$.
\end{definition}
\begin{remark}\label{r:pbp:exact}
As in \fullref{r:reg:ham:isotopies}, it will follow from (the
parametric version of) \fullref{T:hprinc:rel} below that any
Lagrangian immersion $\hat{f}\co L \ra \C^n$ which agrees with the
initial data $(\Sigma,L,f)$ on some neighbourhood of $\Sigma$ is
regularly homotopic to an exact Lagrangian immersion with the same
initial data. For this reason there is no loss of generality in
incorporating exactness directly into the definition of
``solution".
\end{remark}

Any solution of the Prescribed Boundary Problem yields a
solution of the Lagrangian Cobordism Problem, defined by setting
$i=f|_\Sigma$. The prescription of a neighbourhood of $\Sigma$ adds however an important cohomological constraint to the problem, as follows.

Assume for simplicity that $\mathcal{O}\!p\,\Sigma\simeq [0,\epsilon)\times\Sigma$. Notice that the initial data (let us call it a ``Lagrangian thickening" of $i$) determines a Maslov class $\mu_f\in
[\Sigma\times[0,\epsilon),\Sph^1]\simeq [\Sigma,\Sph^1]\simeq
H^1(\Sigma,\Z)$ defined, as before, via the corresponding Gauss
map $Gdf$. Now suppose that $(L,\hat{f})$ is a solution to the induced
Lagrangian Cobordism Problem. If $(L,\hat{f})$ also solves the
Prescribed Boundary Problem, then in particular the Maslov classes $\mu_f$ and
$\mu_{\hat{f}}$ must satisfy
$$i^*\mu_{\hat{f}}=\mu_{f},$$ where $i^*\co H^1(L,\Z)\ra
H^1(\Sigma,\Z)$ is the map induced by the inclusion $\Sigma\times
[0,\epsilon)\subset L$. Hence, the choice of initial data $f$ with a different Maslov class $\mu_f$ will impose a
different condition on $\mu_{\hat{f}}$ but will not change the
induced Lagrangian Cobordism Problem.

In particular, we see that if a solution to the Prescribed
Boundary Problem does exist then $\mu_f$ must belong to the image of $i^*$.
Thus for example if $\mu_f\neq 0$ and $H^1(L,\Z)=0$, the
Prescribed Boundary Problem will not admit solutions even though
by hypothesis the Lagrangian Cobordism Problem does (eg
we could use $\Sigma = S^1$ and $L=D^2$).

The above argument supposes that the same $\Sigma$ can support
many different Lagrangian thickenings, with different Maslov
classes. One can think of several different constructions to show
that this is indeed true. Two of the simplest are outlined in the
examples below. Notice that since Lagrangian thickenings have the
same topology as $\Sigma$, they are exact if and only if $\Sigma$
is.

\begin{example}
\label{E:slg:extend} If the isotropic immersion
$i\co \Sigma^{n-1}\rightarrow \R^{2n}$ is not just smooth but real
analytic, we can produce Lagrangian
thickenings of $i$ as follows. Let $L:=(-\epsilon,
\epsilon)\times\Sigma$. Then, for $\epsilon$ sufficiently small,
and any $e^{i\theta}\in \Sph^1$, an application of Cartan--K\"ahler
theory proves that there exists a unique $\theta$--special
Lagrangian immersion of $L$ extending $i$ \cite[Theorem 
III.5.5]{harvey:lawson}.
\end{example}

\begin{example}\label{e:maslov:thickening}
Suppose $i$ is actually a Lagrangian immersion of $\Sigma$ into
$\C^{n-1}$ so that it defines a Maslov class
$\mu_i\in[\Sigma,\Sph^1]$. Let $f$ be a Lagrangian thickening of
$i$ in $\C^n$ and $\mu_f\in [\Sigma,\Sph^1]$ be its Maslov class. For
example, if $f\co  \Sigma\times [0,\epsilon)\ra\C^n$ is the
``cylindrical thickening" of $i$ given by
$$(x,t)\mapsto (i(x),t)$$
then $\mu_f=\mu_i$. On the other hand, if $i$ is real analytic
then \fullref{E:slg:extend} shows how to build
Maslov-zero thickenings of $i$. This shows that the same $i$ may
admit Lagrangian thickenings with different Maslov classes.
\end{example}

\subsection[The Prescribed Boundary Problem and the relative
h--principle]{The Prescribed Boundary Problem and the relative
$h$--principle} \label{ss:pbp:rel:hprinc}
Our main tool to solve the Prescribed Boundary Problem described
in \fullref{d:initial:data} will again be the Gromov--Lees
$h$--principle, but this time in its relative form.
\begin{theorem}[Lagrangian $h$--principle, relative version, Eliashberg and Mishachev \mbox{\cite[\S 6.2.C, \S 16.3.1]{eliashberg}}]
\label{T:hprinc:rel} Suppose $\Sigma=\partial L$ and that there
exists a Lagrangian monomorphism $F\co TL\rightarrow T\R^{2n}$ which
is holonomic over an open neighbourhood $\mathcal{O}\!p\,\Sigma$ of
$\Sigma\subset L$; ie $F=(f,\phi)$, where $\phi$
satisfies
$\phi|_{\mathcal{O}\!p\,\Sigma}=df|_{\mathcal{O}\!p\,\Sigma}$. Assume
furthermore that $f|_{\mathcal{O}\!p\,\Sigma}$ is exact. Then there
exists a family of Lagrangian monomorphisms $F_t:[0,1]\times
TL\rightarrow T\R^{2n}$ such that $F_0=F$,
$F_{t}|_{\mathcal{O}\!p\,\Sigma}=F_{0}|_{\mathcal{O}\!p\,\Sigma}$ and $F_1$
is holonomic; ie $F_1=d\tilde{f}$. In particular, the
base map $bs(F_1)=\tilde{f}\co L\rightarrow\R^{2n}$ is a Lagrangian
immersion with the prescribed boundary data
$f|_{\mathcal{O}\!p\,\Sigma}$. Furthermore, $\tilde{f}$ is exact.
\end{theorem}
\begin{remark}
\fullref{T:hprinc:rel} still holds without assuming
orientability of $\Sigma$ or $L$.
\end{remark}
\begin{remark}
As in the closed case, the proof of the relative $h$--principle for
Lagrangian immersions relies on lifting the Lagrangian
monomorphism to a Legendrian one in $\R^{2n+1}$ and applying the
relative $h$--principle for Legendrian immersions; again this
explains why the resulting Lagrangian immersion is exact. However,
in order to apply the Legendrian $h$--principle we need this lift
to be holonomic on $\mathcal{O}\!p\,\Sigma$\,: this explains the
extra assumption in the statement.
\end{remark}
\fullref{T:hprinc:rel} shows that to solve the Prescribed
Boundary Problem with exact initial data, the only obstruction one
needs to overcome is the construction of a Lagrangian monomorphism
with the prescribed initial data. The triviality condition on
$TL^\C$ is not sufficient to guarantee that such monomorphisms
exist: the boundary data imposes an additional constraint, which
we formalize as follows.
\begin{definition} \label{d:comp:triv}
Let $(\Sigma,L,f)$ be initial data for the Prescribed Boundary Problem.
A trivialization $TL^\C\rightarrow L \times \C^n$ is
\textit{compatible\/} with the initial data if it extends the
trivialization of $T(\mathcal{O}\!p\,\Sigma)^\C$ induced by
$df^{\C}$.
\end{definition}
The theory now proceeds exactly as in \fullref{S:h:princ:closed}. In particular, the following statement can be proved using those same methods.
\begin{corollary}\label{C:h:princ:rel}
Let $(\Sigma,L,f)$ be exact initial data for the Prescribed
Boundary Problem. Then the following conditions are equivalent.
\begin{enumerate}
\item The Prescribed Boundary Problem with initial data $(\Sigma,L,f)$ admits a solution.
\item $TL^\C$ admits a compatible trivialization.
\item $L$ admits a totally real monomorphism $F\co TL\rightarrow T\R^{2n}$ which extends $df$.
\item $L$ admits a Lagrangian monomorphism $F\co TL\rightarrow T\R^{2n}$ which extends $df$.
\end{enumerate}
Furthermore, there is a one-to-one correspondence between the
homotopy classes of the above objects. In particular, if a
solution to the Prescribed Boundary Problem exists, it induces a
bijective correspondence between the set of all compatible
trivializations of $TL^{\C}$ and the set of all maps $L\rightarrow
\gl{n}$ which map a neighbourhood of $\Sigma$ to the identity
matrix. Two trivializations are homotopic if and only if the
corresponding maps are homotopic.
\end{corollary}
The theory is not yet completely satisfactory: we would still like
a constructive procedure to verify the existence of compatible
trivializations. To this end, let $(\Sigma, L, f)$ be initial data
for the Prescribed Boundary Problem. Notice that since $\R^{2n}$
is contractible there is no obstruction to extending $f$ to a
smooth map $\smash{\hat{f}}\co  L\rightarrow \R^{2n}$; we can thus fix one
such extension $\smash{\hat{f}}$ and use it as the base map for all our
monomorphisms.

Let us now fix a reference trivialization $F_0=(\hat{f},\phi_0)$
of $TL^\C$ (not necessarily compatible with the initial data). As
seen in \fullref{l:trivs:maps} this choice determines a
one-to-one correspondence between all other trivializations of
$TL^\C$ and maps $L\rightarrow \gl{n}$. In particular, $df^\C$
determines a map $M_{df^\C,F_0}\co  \mathcal{O}\!p\,\Sigma\rightarrow
\gl{n}$. The proof of the following result is a direct consequence
of \fullref{d:comp:triv}.
\begin{lemma}[Existence of compatible trivializations]\label{L:extend:GL}
There exists a trivialization compatible with the initial data
$(\Sigma,L,f)$ if and only if the map $M:= M_{df^\C,F_0}$ from 
$\mathcal{O}\!p\,\Sigma$ to $\gl{n}$ determined by $df^{\C}$
and by the choice of reference trivialization $F_0$ can be
extended to a map $\widehat{M}\co  L \ra \gl{n}$.
\end{lemma}
\begin{remark}The existence of a compatible trivialization should be independent of any choices made in the construction of the matrix-valued map $M$; in particular if we change the reference trivialization we will obtain a different map $M'$, and $M$ should be extensible if and only if $M'$ is. Notice that the change of trivialization will be described by a matrix-valued map defined on the whole $L$; using this it is simple to see that the above indeed holds.
\end{remark}
In other words, the fact that $TL^\C$ is trivial allows us to
translate the Prescribed Boundary Problem into an extension
problem for maps into $\gl{n}$. In \fullref{S:obstruction} we will describe a standard
obstruction theory framework for dealing with such extension
problems using algebraic topology, following the treatments in
Bredon \cite{bredon} and Hatcher \cite{hatcher}. In the meantime we extend
the definitions and results of \fullref{s:maslov:class} to manifolds with boundary.

Given a manifold with boundary $L$ and a Lagrangian monomorphism
$F\co TL\ra T\C^n$, we define its Gauss map $GF\co L\ra
\unitary{n}/\sorth{n}$ and determinant map $\det_\C\circ GF\co L\ra
\Sph^1$ exactly as before. This allows us to continue to use the
same definition also for the Maslov map $\mu\co \mathcal{M}\ra
[L,\Sph^1]\simeq H^1(L,\Z)$, where again $\mathcal{M}$ denotes the
set of homotopy classes of Lagrangian monomorphisms.

\begin{lemma} \label{l:maslov:bdy}
Let $L$ be a compact oriented $n$--manifold with boundary. Suppose $TL^\C$
is trivial and fix a reference Lagrangian monomorphism $F_0$.
Then, for any other Lagrangian monomorphism $F$ (homotoped so that the corresponding matrix map $M_{F,F_0}$ takes values in $\unitary{n}$), $\mu_F=[\textrm{det}_\C
M_{F,F_0}]\cdot \mu_{F_0}$. In particular, if $F_0$ has zero
Maslov class then $\mu_F=[\textrm{det}_\C M_{F,F_0}]$.

Furthermore, let $\mu$
be an arbitrary element of $H^1(L,\Z)$. Then $L$ admits a Lagrangian monomorphism $F\co TL\ra T\C^n$ with $\mu_F=\mu$. Finally, any Maslov-zero Lagrangian monomorphism can be homotoped to a SL monomorphism.
\end{lemma}
\begin{proof}
The proofs are exactly the same as those of \fullref{l:matrix:maslov}, \fullref{p:anymaslov:ok} and \fullref{p:MZisSL}.
\end{proof}
Notice that the second part of \fullref{p:anymaslov:ok}
does not extend so easily because the $h$--principle for manifolds
with boundary requires additional assumptions about the
monomorphism near the boundary. Fortunately the above Lemma will
be sufficient for our purposes.

\begin{example}\label{e:maslov:mobius} Concerning the nonorientable case, suppose for example that $L$ is the M\"{o}bius strip. Recall that $L$ is a subset of the Klein bottle $K$, which from \fullref {s:closed:examples} has $TK^\C$ trivial; thus $TL^\C$ is trivial. The analogue of \fullref{p:somemaslov:ok} for nonorientable manifolds with boundary shows that the set of unoriented Maslov classes of $L$ must be either $2\Z$ or $2\Z+1$ inside $H^1(L,\Z)\simeq \Z$. However, since $L$ is nonorientable, $0$ cannot be a Maslov class so the subset realized by Maslov classes is $2\Z+1 \subset
\Z$.
\end{example}

\subsection{The Extension Problem and obstruction theory}\label{S:obstruction}
Let $(\Sigma,L,f)$ be initial data for the Prescribed Boundary
Problem.  By \fullref{L:extend:GL} the existence of a
trivialization compatible with the initial data (in the sense of
\fullref{d:comp:triv}) is equivalent to the extensibility
to $L$ of a map $M\co  \mathcal{O}\!p\, \Sigma \ra \gl{n}$. In this
section first we recall a standard approach to the abstract
extension problem for continuous maps of topological spaces. Then
we apply these results to prove various existence and
nonexistence results for compatible trivializations.

Let $A$ be a subspace of the topological space $X$ and let $M:A
\ra Y$ be a continuous map. The Extension Problem asks:

\medskip
\textbf{Extension Problem}\qua
When can $M$ be extended to a continuous map $\widehat{M}\co X \ra Y$?
\medskip

It is not always possible to give a continuous extension of $M$.
For example, take $(X,A) = (D^{n+1},S^{n})$ and $Y=S^n$. Then
$M$ extends to $D^{n+1}$
if and only if the class $[M] = 0 \in \pi_{n}(S^n) \equiv \Z$. In
particular, $M\co S^n \ra S^n$ extends to $D^{n+1}$ if and only if
$\deg(M)=0$. Thus there is a homotopy-theoretic obstruction to
extending $M$. If we keep $(X,A) = (D^n,S^{n-1})$ and use
$Y=\gl{n}$, then this is simplest case of interest in \fullref{L:extend:GL}. Once again $M\co S^{n-1} \ra \gl{n}$ extends to
$D^n$ if and only if
$[M]=0 \in \pi_{n-1}(\gl{n})$. In particular, if $n$ is odd then
this homotopy group vanishes (see below) and hence from \fullref{C:h:princ:rel} and  \fullref{L:extend:GL} we obtain:
\begin{corollary}
If $n$ is odd, then the Prescribed Boundary Problem is solvable
for any initial data of the form $(D^n,S^{n-1},f)$.
\end{corollary}
On the other hand, if $n$ is even, then $\pi_{n-1}(\gl{n}) = \Z$.
For the case $n=2$, we will show in the next section that the
Prescribed Boundary Problem is solvable if and only if the Maslov
class of the initial data $(D^2,S^1,f)$ is zero; see \fullref{t:filling:n=2} for this result in a slightly more general
context. This is again a homotopy-theoretic obstruction to solving
the Prescribed Boundary Problem. Furthermore, since $S^{n-1}$ is
simply connected for  $n\ge 3$, we see that these obstructions
will be sensitive to more than the Maslov class of the initial
data.

\medskip
\textbf{Note}\qua For the remainder of this section a map between two
topological spaces means a continuous map even if we do not say so
explicitly.
\medskip

Now we would like to understand the Extension Problem for more
general pairs $(X,\!A)$. In the simple extension problems considered
above the extensibility of $M$ depended only on the homotopy class
of $M$. It is natural to ask if this is always the case.

An equivalent reformulation of this question is: given two
homotopic maps $M_0$ and $M_1$ from $A$ to $Y$ and an extension
$\widehat{M}_0\co X \ra Y$ of $M_0$, is there always an extension
$\widehat{M}_1\co X \ra Y$ of $M_1$ which is homotopic to
$\widehat{M}_0$? If for fixed $(X,A)$ and $Y$ we can always find
such a homotopic extension $\widehat{M}_1$ then $(X,A)$ is said to
have the \textit{homotopy extension property\/} with respect to $Y$.

In other words, if $(X,A)$ has the homotopy extension property
with respect to $Y$ then the extensibility of maps $M\co A \ra Y$
depends only on the homotopy class of $M$; in particular in this
case, any map $M\co A \ra Y$ which is homotopic to a constant map has
an extension $\widehat{M}\co X \ra Y$ which is also homotopic to a
constant map ($M$ may also have other homotopically distinct
extensions).

Not all spaces satisfy the homotopy extension property \cite[page
430]{bredon}. However, under rather mild assumptions on $(X,A)$ we
get the much stronger conclusion that $(X,A)$ satisfies the
homotopy extension property with respect to any space $Y$
\cite[VII.1.1]{bredon}. For example, if $(X,A)$ is a CW-pair,
ie $A$ is a subcomplex of a CW-complex $X$, then the
homotopy extension property holds for any space $Y$ and hence the
Extension Problem is a homotopy-theoretic one \cite[Corollary 
VII.1.4]{bredon}.

If, in addition, the space $Y$ is path-connected and simple (that
is, $\pi_1(Y,y_0)$ acts trivially on $\pi_n(Y,y_0)$ for all $y_0
\in Y$ and $n\ge 1$) then one can use a Postnikov decomposition of
$Y$ \cite[page 501]{bredon} to solve the Extension Problem as follows:
\begin{theorem}[Extension/Obstruction Theorem \mbox{\cite[page 507]{bredon}}]
\label{T:extend}  Let $(X,A)$ be a CW-pair and
$Y$ be path-connected and simple and let $M\co A \ra Y$ be a
(continuous) map. Then there exists a sequence of obstructions
$$ c_M^{n+1} \in H^{n+1}(X,A;\pi_n(Y)), \quad n\ge 1,$$
(where $\smash{c_M^{n+1}}$ is defined only when all the previous
obstructions vanish and depends on the previous liftings made)
such that there is a solution of the Extension Problem for $M$ if
and only if there is a complete sequence of obstructions
$\smash{c_M^{n+1}}$ all of which are zero. Furthermore, if $\widehat{M}_1$
and $\widehat{M}_2$ are two different extensions of $M$ then there
exists a sequence of obstructions
$$d^n(\widehat{M}_1,\widehat{M}_2) \in H^n(X,A;\pi_n(Y))$$
to the existence of a homotopy rel $A$ between $\widehat{M}_1$ and
$\widehat{M}_2$.
\end{theorem}
We will discuss later how to construct these obstruction
cocycles $c_M^i$ in some cases.

In our applications $(X,A)= (L,\partial L)$, where $L$ is a smooth
compact oriented manifold with oriented boundary $\partial L$. It
is a classical result of J\,H\,C\,Whitehead \cite{whitehead:1940}
(see also Milnor and Stasheff \cite[page 240]{milnor:stasheff}) that a smooth compact
manifold with boundary can be triangulated as a simplicial complex
with the boundary as a subcomplex. In particular, $(L,\partial L)$
is a CW-pair. Furthermore, the space $Y$ will be a connected Lie
group,
in which case $Y$ is known to be simple \cite[pages 88--89]{steenrod},
and hence \fullref{T:extend} applies. For the application of
obstruction theory to \fullref{L:extend:GL} and the construction
of compatible trivializations we are interested in the case where
$Y=\gl{n}$; the retraction $\gl{n}\rightarrow \unitary{n}$ allows us, from now on, to concentrate on the case $Y=\unitary{n}$.

For our applications of \fullref{T:extend} we are only
interested in $\pi_i(\unitary{n})$ for $i\le n$. These fall within
the ``stable range'' where $[(i+2)/2] \le n$ and hence by Bott
Periodicity are equal to $0$ when $n$ is even and $\Z$ when $n$ is
odd \cite[pages 467--8]{bredon}. Hence we have the following
immediate corollary of \fullref{T:extend}:
\begin{corollary}\label{c:un:obstruc}
A map $M\co \Sigma^{n-1} \ra \unitary{n}$ extends to $L$ if and only
there is a sequence of $[n/2]$ obstructions $$c_M^{2i}
\in H^{2i}(L,\Sigma;\Z)\quad  \textrm{for\ }  i \in \{1, \ldots ,
[n/2] \}$$ which all vanish (where
$\smash{c_M^{2i}}$ is defined only when all the previous obstructions
vanish and depends on the previous extensions made)
. Furthermore if $\widehat{M}_1$ and
$\widehat{M}_2$ are two different extensions of $M$ there is
a sequence of obstructions
$$d^{2i+1}(\widehat{M}_1,\widehat{M}_2) \in H^{2i+1}(L,\Sigma;\Z), \quad \textrm{for\ }
i=0, \ldots ,[(n-1)/2],$$ to the existence of a homotopy rel $A$
between $\widehat{M}_1$ and $\widehat{M}_2$.
\end{corollary}
To analyze further the extensibility of maps into $\unitary{n}$ it
is convenient first to break the problem up into two separate
extension problems. There are two advantages of splitting the
problem into these separate extension problems. The first
advantage is that it shows that the higher obstructions cocycles
$\smash{c^{2i}_M}$ ($i>1$) for a map into $\unitary{n}$ can be studied
independently of the first obstruction cocycle $\smash{c_M^2}$. The second
advantage is that it will allow us in \fullref{c:extend:s1}
to identify the first obstruction cocycle $\smash{c_M^{2}}$ in a very
concrete manner. The splitting of the $\unitary{n}$ extension
problem is achieved as follows:

Given $U \in \unitary{n}$ define a map $\pi_{\Sph^1}:\unitary{n}
\ra \Sph^1$ by $\pi_{\Sph^1}(U)= \det_{\C}U$ and a map
$\pi_{\sunitary{n}}\co \unitary{n} \ra \sunitary{n}$ by
$\pi_{\sunitary{n}}(U) = \diag(\pi_{\Sph^1(U)},1, \ldots, 1)^{-1}U$.
Since any $U$ in $\unitary{n}$ may be written uniquely as a product
$U = \diag(\pi_{\Sph^1}(U),1,\ldots ,1) \cdot
\pi_{\sunitary{n}}(M)$, we see that $\unitary{n}$ is diffeomorphic
to $\Sph^1 \times \sunitary{n}$. In particular, any map $M$ from
any topological space $X$ into $\unitary{n}$ induces two composition maps
$M_{\Sph^1}:= \pi_{\Sph^1} \circ M$ and $M_{\sunitary{n}}:=
\pi_{\sunitary{n}} \circ M$ from $X$ into $\Sph^1$ and
$\sunitary{n}$ respectively. Conversely any pair of maps from $X$
into $\Sph^1$ and $\sunitary{n}$ determines a unique map from $X$
into $\unitary{n}$. In particular, a  map $M\co \Sigma \ra
\unitary{n}$ extends to $\widehat{M}\co L \ra \unitary{n}$ if and
only if the maps $M_{\Sph^1}\co \Sigma \ra \Sph^1 $ and
$M_{\sunitary{n}}\co \Sigma \ra \sunitary{n}$ both extend to $L$.

\begin{remark}\label{r:homotopy:U}
Using the observation above it follows that for any CW-complex $K$
the homotopy classes of maps $[K;\unitary{n}]$ can be identified
(as a set but not necessarily as a group) with $[K;\Sph^1] \times
[K;\sunitary{n}] \simeq H^1(K,\Z) \times [K;\sunitary{n}]$ (using
\fullref{l:homotopy:s1} for the final identification).
\end{remark}
We must now apply the \fullref{T:extend}
to maps with target $\sunitary{n}$ or $\Sph^1$ and compare the
results obtained to those of \fullref{c:un:obstruc}, where
the target is $\unitary{n}$. In particular, we want to see how the
$[n/2]$ obstructions which occur in the $\unitary{n}$ case split
into two parts according to the decomposition of the map into its
$\Sph^1$ and $\sunitary{n}$ parts.

 Since $\unitary{n} \simeq
\Sph^1 \times \sunitary{n}$ as manifolds (but not as Lie groups),
$\pi_i(\unitary{n}) = \pi_i(\sunitary{n})$ for
$i\ge 2$. We also have $\pi_1(\sunitary{n})=0$. Therefore, given a
    map $M\co \Sigma \ra \sunitary{n}$, by \fullref{T:extend} there is a 
    sequence of $[n/2]-1$ obstructions $\smash{c_M^{2i}} \in
H^{2i}(L,\Sigma;\Z)$ for $i \in \{2, \ldots , [n/2] \}$ to
extending $M$ to $\widehat{M}\co L \ra \sunitary{n}$. In particular,
the first obstruction to extending a map $M$ into $\unitary{n}$
disappears if instead $M$ maps into $\sunitary{n}$. In particular,
if $n=2$ or $n=3$ then there are no obstructions to extending $M$.
In other words, we have the following extension result for maps
into $\sunitary{2}$ or $\sunitary{3}$.
\begin{corollary}
\label{c:SU:nice} Let $L$ be a compact smooth $n$--manifold with
boundary $\Sigma$ and let $M\co \Sigma \ra \sunitary{n}$ be any map.
Then for $n=2,3$, $M$ extends to a map $\widehat{M}\co  L \ra
\sunitary{n}$. Furthermore, the homotopy classes of maps $L \ra
\sunitary{n}$ rel $\Sigma$ which equal $M$ on $\Sigma$ are in
one-to-one correspondence with $0$ and $\Z$ for $n=2$ and $n=3$
respectively.
\end{corollary}
\begin{remark} \fullref{c:SU:nice} is already sufficient for a proof of \fullref{c:SL:n=2} (see the second method of proof given there) and \fullref{c:SL:n=3}, and (via \fullref{r:PBP:MZcase}) for the applications of \fullref{s:desing}.
\end{remark}
The difference between the extensibility of maps $M\co  \Sigma \ra
\unitary{n}$ and $M_{SU}\co  \Sigma \ra \sunitary{n}$ will be
captured by the extensibility of the map $M_{\Sph^1}\co  \Sigma \ra
\Sph^1$. Since $\Sph^1$ has trivial homotopy groups for $i>1$,
given a map $M\co  \Sigma \ra \Sph^1$ then \fullref{T:extend}
gives us precisely one obstruction cocycle $c_M^{2} \in
H^2(L,\Sigma;\Z)$ to extending it to $L$.

Hence as expected in total we have the same number of obstruction
cocycles (namely $[n/2]$ of them) which measure the
nonextensibility for a map into $\unitary{n}$ or for a pair of
maps into $\sunitary{n}$ and $\Sph^1$; one of these cocycles
encodes the extensibility of the map into $\Sph^1$, while the
remaining $[n/2]-1$ encode the extensibility of the map into
$\sunitary{n}$.

To proceed further we need to understand better the obstruction
cocycles. For the construction of these cocycles $c_M^{i+1}$ in the
general case we refer the reader to Bredon \cite[VII.13]{bredon}.
However, in the case that the only possible nonzero obstruction is
the so-called primary one it is simple to describe this
obstruction. This will allow us to understand very concretely the
obstruction to extending a map into $\Sph^1$.
\begin{theorem}{\rm{\cite[Corollary VII.13.13]{bredon}}}\qua
\label{t:primary:obstruct} Let $(X,A)$ be a CW-pair and $Y$ be
simple and $(k-1)$--connected. Suppose that
$H^{i+1}(X,A;\pi_i(Y))=0$ for all $i>k$. Then a map $M\co A \ra Y$
can be extended to a map $\widehat{M}\co X \ra Y$ if and only if the
homomorphism
$$\delta^* M^*\co  H^k(Y;\pi) \ra H^{k+1}(X,A;\pi)$$ is trivial,
where $\pi = \pi_k(Y)$ and $\delta^*$ is the usual map $\delta^*:
H^k(A;\pi) \ra H^{k+1}(X,A;\pi)$ which appears in the exact
cohomology sequence of the pair $(X,A)$.
\end{theorem}
\begin{remark}
\label{r:obstruct} There is an easy way to see that the
nontriviality of the map $\delta^* M^*$ of \fullref{t:primary:obstruct}
is an obstruction to extending $M$. Suppose
$\widehat{M}\co X \ra Y$ is an extension of $M$. Then $\smash{\widehat{M}}
\circ i = M$ and hence $\delta^* M^* = \delta^* i^* \smash{\widehat{M}}^*
=0$ since $\delta^* i^* = 0$ by the exactness of the cohomology
sequence of any pair $(X,A)$.
\end{remark}
The following elementary lemma about $\delta^*$ is useful for
applications of \fullref{t:primary:obstruct}.

\begin{lemma}
\label{l:delta:star} Let $L$ be a compact connected $n$--manifold
$M$ with boundary $\partial L=\Sigma$. The map
$\delta^*\co H^{n-1}(\partial L;G) \ra H^n(L,\partial L;G)$ is surjective for
any coefficient group $G$. If additionally, $L$ is orientable and
$\partial L$ is connected, then $\delta^*$ is an isomorphism.
\end{lemma}
\begin{proof}
Any compact $n$--manifold with nonempty boundary has the homotopy
type of a CW complex of dimension at most $n-1$ \cite[page 170]
{kosinski}. In particular, $H^n(L)=0$ for any coefficient group
$G$. Surjectivity of $\delta^*\co H^{n-1}(\partial L) \ra
H^n(L,\partial L)$ now follows immediately from the exactness of
the cohomology sequence of the pair $(L,\partial L)$.
If $L$ is orientable, then we may use Poincar\'e--Lefschetz duality
\cite[page 357]{bredon} to identify $\ker{\delta^*}$ with
$\ker{i_*}$, where $i_*\co H_0(\partial L) \ra  H_0(L)$ is the map
induced by the inclusion $i\co \partial L \ra L$. If $\partial L$ is
connected then $i_*$ is an isomorphism and hence so is $\delta^*$.
\end{proof}
\begin{remark}
If $L$ is orientable but $\Sigma$ is disconnected with $k\ge 2$
components, then $\delta^*\co  H^{n-1}(\Sigma) \ra H^n(L,\Sigma)$
 fails to be an isomorphism (since $H^n(L,\Sigma;G)=G$
while $H^{n-1}(\Sigma;G) \simeq H_0(\Sigma;G) = G^k$).
\end{remark}
It remains to use \fullref{t:primary:obstruct} to study the
obstructions to extending maps into $\Sph^1$ and $\sunitary{n}$.
We begin with maps into $\Sph^1$.
\begin{corollary}
\label{c:extend:s1} Let $L$ be a compact connected smooth
$n$--manifold with boundary $\Sigma$ and $i\co \Sigma \ra L$ denote
the natural inclusion map. A map $M\co  \Sigma^{n-1} \ra \Sph^1$
extends to a map $\widehat{M}\co L^n \ra \Sph^1$ if and only if
$$ M^*[d\theta] \in \imag(i^*\co H^1(L,\Z)\ra H^1(\Sigma,\Z)).$$
Hence if $M$ is homotopic to a constant map then $M$ extends to
$L$. Moreover, in this case there is an extension $\widehat{M}$
which is homotopic to a constant map from $L$, and the possible
homotopy classes of different extensions of $M$ are parameterized
by $H^1(L,\Sigma;\Z)$.

Conversely, if either
\begin{itemize}
\item[(i)] $L$ is an orientable surface with connected boundary
$\Sigma$, or
\item[(ii)] $H^1(L;\Z)=0$,
\end{itemize}
then $M$ extends to $L$ if and only if $M$ is homotopic to a
constant.
\end{corollary}
\begin{proof}
Since $\Sph^1$ is $0$--connected and its higher homotopy groups
$\pi_i$ for $i\ge 2$ all vanish then $M$ extends if and only if the map
$\delta^* M^*\co H^1(\Sph^1;\Z) \ra H^2(L,\Sigma;\Z)$ is zero by \fullref{t:primary:obstruct}. Since
$H^1(\Sph^1;\Z) = \Z$ and is generated by $[d\theta]$, this is
true if and only if $\delta^* M^*[d\theta]=0$. But by the
exactness of the long exact sequence in cohomology of the pair
$(L,\Sigma)$, $M^*[d\theta] \in \ker{\delta^*}$ if and only if
$M^*[d\theta] \in \imag{i^*}$. Clearly, if $H^1(L,\Z)=(0)$ then
$M^*[d\theta] \in \imag{i^*}$ if and only if $M^*[d\theta]=0$,
which by \fullref{l:homotopy:s1} is equivalent to $M$ being
homotopic to a constant. In case (i), it follows from \fullref{l:delta:star} that $\delta^*\co H^1(\Sigma;\Z) \ra H^2(L,\Sigma;
\Z)$ is an isomorphism. Hence $\delta^* M^*=0$ if and only if
$M^*[d\theta]=0$, as in case (ii). If $M$ is homotopic to a
constant then since the pair $(L,\Sigma)$ satisfies the homotopy
extension property for any space $Y$ then there is an extension
$\widehat{M}\co L \ra \Sph^1$ which is homotopic to a constant map.
The statement about the different homotopy classes of extensions
follows from the second part of \fullref{T:extend} since the
only nonzero obstruction space $H^i(L,\Sigma,\pi_i(\Sph^1))$
occurs when $i=1$.
\end{proof}
\begin{remark}\label{r:maslov:mobius}
It is interesting to check what \fullref{c:extend:s1} says in the nonorientable case. For example, suppose $L$ is the M\"{o}bius strip with boundary $\Sigma=S^1$. Then $H^1(L,\Z)\simeq \Z \simeq H^1(\Sigma,\Z)$ and $i^*$ is the map $n\mapsto 2n$, so it has image $2\Z\subset H^1(\Sigma,\Z)$. Thus a map $M\co \Sigma\ra \Sph^1$ extends to $L$ if and only if $M^*[d\theta]\in 2\Z$. From \fullref{e:maslov:mobius}, the image under $i^*$ of the possible Maslov classes of $L\subset\C^2$ is the set $4\Z+2\subset H^1(\Sigma,\Z)$. On the other hand, the standard immersion of $\Sigma$ in $\C$ has Maslov class $2\in H^1(\Sigma,\Z)$ (when calculated via $\textrm{det}^2_\C$); this is also the Maslov class of any exact perturbation of this immersion, and of the corresponding ``cylindrical thickening" $f$ in $\C^2$ (see \fullref{e:maslov:thickening}). One could use these facts as a basis for proving that the Prescribed Boundary Problem determined by this $(\Sigma,L,f)$ is solvable. The solution is actually known explicitly \cite{audin:cobord}.
\end{remark}

Similarly, if we apply \fullref{t:primary:obstruct} with
$Y=\sunitary{n}$ then we obtain the following extension of
\fullref{c:SU:nice}.
\begin{corollary}
\label{c:SU:obstruct} Let $L$ be a compact connected smooth
$n$--manifold with boundary $\Sigma$ and $i\co \Sigma \ra L$ denote
the natural inclusion map. Let $M\co \Sigma \ra \sunitary{n}$ be a
map. Then for $2\le n \le 5$, $M$ extends to a map of $L$ if and
only if
$$M^* \Theta \in  \imag(i^*\co H^3(L,\Z)\ra H^3(\Sigma,\Z))$$
where $\Theta$ is the generator of $H^3(\sunitary{n},\Z) =\Z$.

Furthermore, if $n=4$ and $L$ is orientable with connected
boundary $\Sigma$ then $M$ extends if and only if $M$ is homotopic
to a constant map.
\end{corollary}
\begin{proof}
The proof of the first part is entirely analogous to the proof of
the previous Corollary, the only difference being that
$\sunitary{n}$ is $2$--connected. If $n=4$, $L$ is orientable and
$\Sigma$ is connected, then by \fullref{l:delta:star}
$\delta^*\co H^3(\Sigma;\Z) \ra H^4(L;\Sigma,\Z)$ is an isomorphism.
Hence $\delta^* M^*=0$ if and only if $M^*\Theta=0$. But since
once can prove that the map $M\mapsto M^*\Theta$ induces an
isomorphism between the groups $[\Sigma^3,\sunitary{4}]$ and
$H^3(\Sigma^3;\Z)$, $M^*\Theta=0$ if and only if $M$ is homotopic
to a constant map.
\end{proof}
If $n>5$ the situation is not as simple since there are
further obstructions to extending the map $M_{\sunitary{n}}$,
eg a second obstruction $c_M^6\in
H^6(L,\Sigma,\pi_5(SU(n)))=H^6(L,\Sigma,\Z)$ can now be nonzero.

We now summarize the main results of this section in terms of the
set-up introduced in the context of \fullref{L:extend:GL}.
\begin{theorem}\label{t:obstruction:summary}
Let $(\Sigma^{n-1},L^n,f)$ be initial data for the Prescribed
Boundary Problem and assume $2\le n\le 5$. Fix a reference
Maslov-zero Lagrangian monomorphism $F_0$ (this exists by \fullref{l:maslov:bdy}); let $M:=M_{df^\C,F_0}:\mathcal{O}\!p\,\Sigma
\ra \gl{n}$ be the map corresponding to $df^\C$ via $F_0$. Using
the retraction $\gl{n}\rightarrow \unitary{n}$, $M$ is homotopic
to a map which takes values in $\unitary{n}$; we will continue to
denote this map $M$. Let $M_{\Sph^1}$ and $M_{\sunitary{n}}$
denote the associated maps into $\Sph^1$ and $\sunitary{n}$
respectively. Then there exists a compatible trivialization if and
only if 
\begin{enumerate}
\item[(i)] $M_{\Sph^1}^*[d\theta] \in \imag(i^*\co H^1(L,\Z)\ra H^1(\Sigma,\Z))$, and
\item[(ii)] $M_{\sunitary{n}}^*\Theta \in \imag(i^*\co H^3(L,\Z)\ra H^3(\Sigma,\Z)),$
\end{enumerate}
where $\Theta$ is the generator of $H^3(\sunitary{n},\Z)$, $i$ is
the natural immersion $\Sigma\rightarrow L$, and the pullback
operation is performed using the restrictions of $M_{\Sph^1}$ and
$M_{\sunitary{n}}$ to $\Sigma$.

Furthermore, there is a natural identification of
$M_{\Sph^1}^*[d\theta]$ with the Maslov map $\mu_f$ of the initial
data.
\end{theorem}
\begin{proof} Notice that by the homotopy extension property (setting $Y:=\gl{n}$ and thinking of $\unitary{n}\subset \gl{n}$)
the map into $\gl{n}$ is extensible if and only if the map into
$\unitary{n}$ is. The idea of the proof is now of course to apply
the previous obstruction theory results to the maps $M_{\Sph^1}$
and $M_{\sunitary{n}}$. Notice however that, as required in \fullref{L:extend:GL}, we want to extend not only the values of these
maps on $\Sigma$ but also on a neighbourhood of $\Sigma$.
Furthermore we want the extension to be smooth. To achieve this,
let us choose a (closed, sufficiently small) tubular neighbourhood
$N\simeq\Sigma\times[0,\epsilon]\subseteq\mathcal{O}\!p\,\Sigma$;
let $\Sigma_\epsilon$ denote the ``inner boundary"
$\Sigma\times\{\epsilon\}\subset N$. We can apply the obstruction
theory results to $M_{\Sph^1}$ and $M_{\sunitary{n}}$ restricted
to $\Sigma_\epsilon$, obtaining obstructions
$M_{\Sph^1}^*[d\theta]$ and $M_{\sunitary{n}}^*\Theta$.

\medskip
\textbf{Note}\qua To be precise, these cocycles live in
$H^1(\Sigma_\epsilon,\Z)$ and $H^3(\Sigma_\epsilon,\Z)$ but we can
identify these with $H^1(\Sigma,\Z)$ and $H^3(\Sigma,\Z)$,
respectively.
\medskip

Conditions (i) and (ii) now follow from
\fullref{c:extend:s1} and \fullref{c:SU:obstruct}. When the conditions are
satisfied we obtain continuous extensions defined on the
complement $L\setminus N$ of the tubular neighbourhood; together
with the given values of $M_{\Sph^1}$ and $M_{\sunitary{n}}$ on
$N$, we now have continuous extensions defined on the whole of
$L$. Standard perturbation results show that we can assume that
these extensions are smooth.

It remains only to show that $M_{\Sph^1}^*[d\theta]$ can be
identified with the Maslov class of the Lagrangian map
$f\co \mathcal{O}\!p\,\Sigma\ra \C^n$. By considering a sufficiently
small neighbourhood of $\Sigma$, we can assume that
$\mathcal{O}\!p\,\Sigma$ is contractible onto $N$. The homotopy
invariance of $\mu_f$ then allows us to restrict $f$ to $N$.
Recall from \fullref{l:maslov:bdy} that, since $F_0$ is
Maslov-zero,  $\mu_f$ is the homotopy class of $\det_\C M:
N\ra\Sph^1$. The claim then follows from the fact that $\det_\C
M=M_{\Sph^1}$ and the identifications $[N,\Sph^1]\simeq
[\Sigma,\Sph^1]\simeq H^1(\Sigma,\Z)$ (where the last isomorphism
is via pullback of $[d\theta]$; see \fullref{l:homotopy:s1}).
\end{proof}
\section{The Prescribed Boundary Problem in low dimensions}\label{s:pbp:lowdim}
\fullref{C:h:princ:rel} showed that solving the Prescribed
Boundary Problem is equivalent to finding a compatible
trivialization, which by \fullref{L:extend:GL} is equivalent to
solving a certain extension problem. \fullref{t:obstruction:summary} shows that in low dimensions the
obstructions to solving this extension problem can be written down
explicitly. We will now examine these cases one by one.

\medskip
\textbf{Note}\qua To simplify the exposition, throughout this section we will assume that $\mathcal{O}\!p\,\Sigma$ is a tubular neighbourhood of $\Sigma\subset L$ so that $\mathcal{O}\!`\,\Sigma\simeq\Sigma\times[0,\epsilon)$. As in the proof of \fullref{t:obstruction:summary}, we will use this to identify $\mu_f\in H^1(\mathcal{O}\!p\,\Sigma,\Z)$ with the corresponding element in $H^1(\Sigma,\Z)$.

\subsection[The Prescribed Boundary Problem in R4]{The Prescribed Boundary Problem in $\R^4$}
\label{ss:pbp:n=2} Suppose we are given an immersion
$i\co \Sigma\rightarrow \R^4$ where $\Sigma$ is a finite union of
oriented circles. Our first task is to check whether $\Sigma$ and
$i$ satisfy \fullref{assumptions}.

It is clear that $\Sigma$ always admits orientable fillings $L$.
If $\Sigma$ is connected then $\Sigma=\partial L$ for any
orientable surface $L$ with connected boundary; if $\Sigma$ is not
connected, the claim now follows via oriented connect sums.
Furthermore, (unlike the case of closed surfaces) any orientable
surface with nonempty boundary is parallelizable. (A proof of
this fact can be seen as follows: any $n$--manifold with nonempty
boundary has the homotopy type of a CW complex of dimension at
most $n-1$. Since any orientable surface with boundary can be
embedded in $\R^3$ it is stably parallelizable. But any stably
trivial bundle of rank $k>d$ is actually trivial over any space
homotopy equivalent to a CW complex of dimension at most $d$
\cite[Corollary IX.1.5]{kosinski}. Hence the tangent bundle of $L$ is
trivial.) Thus assumptions A and B are both satisfied.

For dimensional reasons $i$ is always isotropic. In general
however it will not satisfy the exactness condition C, so we must
make this additional assumption on $i$. \fullref{e:1d:exact}
shows that if $\Sigma$ is connected then exactness of $\Sigma$ is
in fact a necessary condition for the existence of a smooth
Lagrangian extension $f\co L \ra \R^4$.

Before applying the obstruction theory results to this set-up, we
begin by making the following simple observation.
\begin{lemma}
\label{l:disk:maslov:zero} Let $\gamma$ be the boundary of an
immersed oriented Lagrangian surface $f\co L \ra \R^4$ and suppose
that $\gamma$ is connected. Then the Maslov index of $\gamma$,
$\mu_f([\gamma])$, must be zero.
\end{lemma}
\begin{proof}
This is immediate since $\mu_f([\gamma])$ depends only on the
homology class of the loop $\gamma$.
\end{proof}
\begin{remark}
It follows from \fullref{l:disk:maslov:zero} that any connected separating
curve in a closed oriented Lagrangian surface has zero Maslov
index.
\end{remark}
\fullref{l:disk:maslov:zero} shows that when $\Sigma$ is
connected a smooth oriented solution to the Prescribed Boundary
Problem exists only if the initial data $(\gamma, L^2,f)$ is
Maslov-zero. One can give concrete examples of exact initial data
which are not Maslov zero as follows.
\begin{example}
\label{e:non:fillable:data} Given any pair of relatively prime
positive integers $p$, $q$ define a curve $\gamma \subset \Sph^3$
by
$$\gamma_{p,q}(\theta):= \frac{1}{\sqrt{p+q}}\left( \sqrt{q}
e^{ip\theta},i\sqrt{p}e^{-iq\theta}\right), \quad \theta \in
[0,2\pi].$$ Each curve $\gamma_{p,q}$ is a Legendrian curve in
$\Sph^3$, and hence is exact. Since $\gamma_{p,q}$ is Legendrian
the cone over $\gamma_{p,q}$, denoted $C_{p,q}$ is Lagrangian. It
is not difficult to show that the Maslov index of the curve
$\gamma_{p,q}$ in the Lagrangian cone $C_{p,q}$ is $p-q$
\cite[Remark 4.4]{schoen:wolfson:volume}. (The cones $C_{p,q}$ were
studied by Schoen--Wolfson and are important because up to unitary
transformations they account for all the Hamiltonian stationary
cones in $\C^2$.)
Now take $\gamma$ to be $\gamma_{p,q}$ and take the Lagrangian
thickening of $\gamma_{p,q}$ to be a truncation of the cone over
$\gamma$, eg we could take $f\co \gamma \times [0,1) \ra
\R^4$ to be $f(\theta,t)=(1+t)\gamma_{p,q}(\theta)$). Then for any oriented filling $L$
the exact initial data $(\Sigma, L^2, f)$ has nonzero Maslov class
provided $p\neq q$, and hence there is no smooth solution to the
Lagrangian Prescribed Boundary Problem with such initial data.
\end{example}
The main result of this section, \fullref{t:filling:n=2},
provides a converse to \fullref{l:disk:maslov:zero}. The reader
may also like to compare \fullref{l:disk:maslov:zero} with the
following result of Schoen--Wolfson: given any exact smooth Jordan
curve $\gamma$ in $\R^4$ and $m\in \Z$, there exists a
\textit{piecewise\/} smooth Lagrangian disk $D$ bounding $\gamma$ so
that $\gamma$ has Maslov index equal to $m$ \cite[Proposition
4.1]{schoen:wolfson:volume}.

We now state the main result of this section.
\begin{theorem}\label{t:filling:n=2}
Let $L$ be a connected oriented surface with boundary
$\Sigma$ and suppose $(\Sigma,L,f)$ is exact initial data in the
sense of \fullref{d:initial:data}.
\begin{enumerate}
\item If $\Sigma$ is connected then the Prescribed Boundary Problem
is solvable if and only if $f$ has zero Maslov class, ie $\mu_f=0$.
\item In general, the Prescribed Boundary Problem is solvable if
and only if
$$\mu_f\in \imag(i^*\co H^1(L,\Z)\ra H^1(\Sigma,\Z)).$$
Furthermore, if $\mu_f=0$ then one can always find a solution with
zero Maslov class.
\end{enumerate}
\end{theorem}
\begin{proof}
Part (2) is a direct consequence of \fullref{t:obstruction:summary}; in dimension 2 the second condition
there is vacuous, ie the matrix-valued map $M_{SU(2)}$
corresponding to $df^\C$ is always extensible (see also \fullref{c:SU:nice}). The  existence of a Maslov-zero solution follows
from the fact that if $\mu_f$ is zero then the map $M_{\Sph^1}$ is
homotopic to a constant map by \fullref{l:homotopy:s1}.
\fullref{c:extend:s1} then shows that this map admits a
homotopically constant extension.

Part (1) now follows from case (i) of \fullref{c:extend:s1}.
\end{proof}
\begin{example}
Let $L$ be any closed $n$--manifold with $TL^{\C}$ trivial. By \fullref{p:anymaslov:ok}, given
any $\mu \in H^1(L,\Z)$, there
exists an exact Lagrangian immersion $f\co L\ra \C^n$ with
$\mu_f=\mu$. Hence if $H^1(L,\Z)\neq 0$, then we have Maslov
nonzero exact Lagrangian immersions of $L$.

Let $f\co L \ra \C^n$ be any exact Lagrangian immersion with nonzero
Maslov class.
The immersion $C(f): L \times [0,1]\subset \C^n \times \C$ defined
by $C(f) (p,t)=(f(p),t)$ defines an exact Lagrangian ``generalized
cylinder'' in $\C^{n+1}$ and gives a solution to the Prescribed
Boundary Problem determined by the data $C(f): L \times
[0,\epsilon)\cup L\times (1-\epsilon,1] \ra \C^{n+1}$. It is clear
that this initial data has $\mu_{C(f)}\neq 0$, and thus shows that
part (1) of \fullref{t:filling:n=2} is false if we do not
assume that $\Sigma$ is connected.
\end{example}
An immediate corollary of \fullref{t:filling:n=2} is the
following result:
\begin{corollary}
\label{c:SL:n=2} Let $(\Sigma,L^2,f)$ be exact SL initial data.
Then the Prescribed Boundary Problem admits a Maslov-zero
solution.

In particular, let $\Sigma$ be an exact real-analytic curve of
$\R^4$ (eg any real-analytic Legendrian curve in
$\Sph^3$). Choose any oriented surface $L$ bounding $\Sigma$. Then
for any $\theta \in [0,2\pi)$ there exists an exact Lagrangian
immersion $f_{\theta}\co  L\rightarrow \R^4$ with zero Maslov class
which coincides in a neighbourhood of $\Sigma$ with the
$\theta$-SL extension of $\Sigma$.
\end{corollary}
\begin{proof}
If $f$ is SL it is Maslov-zero, so we can apply the last statement
of \fullref{t:filling:n=2}. The statement regarding real
analytic curves is based on \fullref{E:slg:extend}.
\end{proof}
For applications in \fullref{s:desing} it will be useful to
give a second proof of the previous result.

\medskip
\textbf{Second proof of \fullref{c:SL:n=2}}\qua According to \fullref{l:maslov:bdy}, in \fullref{t:obstruction:summary} we
can choose the reference trivialization to be SL. If we start off with SL initial
data, the map $M_{\Sph^1}$ in \fullref{t:obstruction:summary}
will have constant value 1, so we can extend it to a map on $L$
with constant value 1. \fullref{c:SU:nice} shows that there
is no obstruction to extending $M_{\sunitary{2}}$. We thus obtain
a SL monomorphism $\hat{F}$ extending $df$; in particular
$\hat{F}$ is Maslov-zero. We can homotope $\smash{\hat{F}}$ using the
$h$--principle to prove the existence of an exact Maslov-zero
Lagrangian extension.

\begin{remark}\label{r:PBP:MZcase} More generally, the idea in the second proof of \fullref{c:SL:n=2} can also be used to prove part of \fullref{t:filling:n=2} without relying on the full set of results presented in \fullref{S:obstruction}. Assume the initial data $f$ is (exact and) Maslov-zero. Choose a reference SL trivialization $\phi_0$ of $TL^\C$. Again using \fullref{l:maslov:bdy}, it is possible to homotope $df^\C$ to a SL monomorphism of $\mathcal{O}\!p\,\Sigma$ which, via matrices and \fullref{c:SU:nice}, admits an extension to $L$. By the homotopy property of extensions, this means that $df^\C$ is extensible to a Lagrangian monomorphism of $L$. Applying the $h$--principle to this monomorphism allows us to conclude that the Prescribed Boundary Problem admits a Maslov-zero solution for such initial data. Clearly, this works also for $n=3$.
\end{remark}

\subsection[The Prescribed Boundary Problem in R6]{The Prescribed Boundary Problem in $\R^6$}
\label{ss:pbp:n=3}
Let $\Sigma$ be any compact oriented surface not
necessarily connected. If $\Sigma$ is connected with genus $g$
then it admits one obvious ``topological filling'' $L^3$, namely,
the standard genus $g$ handlebody. One can then construct more
complicated fillings of $\Sigma$ using the standard topological
surgeries. Using oriented connect sums one can show that $\Sigma$
admits (connected) fillings even if it is disconnected. Recall
also that any compact orientable $3$--manifold $L$ is
parallelizable. Finally, \fullref{e:exact:surfs} below shows
that any compact oriented surface $\Sigma$ admits exact immersions
$i$ into $\R^6$. Hence \fullref{assumptions} can be
satisfied for any compact orientable surface $\Sigma$.

We now apply our obstruction theory results to this case; once
again, the second obstruction identified in \fullref{t:obstruction:summary} vanishes for dimensional reasons so
the results (and proofs) are largely analogous to those obtained
in dimension 2.
\begin{theorem}\label{t:filling:n=3}
Let $L$ be a connected oriented $3$--manifold with boundary
$\Sigma$ and let $(\Sigma,L,f)$ be exact initial data in the
sense of \fullref{d:initial:data}. Then the Prescribed Boundary Problem is solvable if and only if
$\mu_f\in \imag(i^*\co H^1(L,\Z)\ra H^1(\Sigma,\Z))$.

Furthermore, if $\mu_f=0$ then one can always find a solution with
zero Maslov class.
\end{theorem}
As in the previous subsection, this result has an immediate corollary:
\begin{corollary}\label{c:SL:n=3}
Let $(\Sigma,L^3,f)$ be exact SL initial data. Then the Prescribed
Boundary Problem admits a Maslov-zero solution $\tilde{f}\co L \ra
\R^6$.

In particular, let $\Sigma$ be a real-analytic exact surface in
$\R^6$ (eg any real-analytic Legendrian submanifold of
$\Sph^5$). Choose any oriented $3$--manifold $L$ bounding $\Sigma$.
Then for any $\theta \in [0,2\pi)$ there exists an exact
Lagrangian immersion $f_\theta\co  L\rightarrow \R^6$ with zero
Maslov class which coincides in a neighbourhood of $\Sigma$ with
the $\theta$-SL extension of $\Sigma$.
\end{corollary}
\begin{example}\label{e:exact:surfs}
Since the standard contact structure on $\Sph^5$ with one point
removed is contactomorphic to the standard contact structure on
$\R^5$ \cite{geiges}, a closed surface $\Sigma$ admits Legendrian
immersions into $\Sph^5$ if and only if it admits a Legendrian
immersion into $\R^5$ with its standard contact structure. Hence
by previous remarks $\Sigma$ admits Legendrian immersions into
$\Sph^5$ if and only if it admits exact Lagrangian immersions into
$\R^4$. But we saw in \fullref{s:closed:low:dim} that any
compact orientable surface $\Sigma$ admits such immersions. Thus,
it is possible to find an exact immersion of $\Sigma$ in $\R^6$
for any orientable surface $\Sigma$.
\end{example}
Another special case for \fullref{t:filling:n=3} occurs when
$H^1(L,\Z)=0$: in this case the condition on $\mu_f$ simplifies to
$\mu_f=0$. However, for $n\ge3$ there are constraints on the
topological complexity of any $L$ which bounds $\Sigma$.
\begin{lemma}\label{l:top:filling:low:bd}
Let $L$ be a smooth oriented $(2i+1)$--dimensional manifold with
connected boundary $\Sigma$. Then $\dim H^i(L,\R) \ge
\frac{1}{2}b^i(\Sigma)$.
\end{lemma}
The proof of \fullref{l:top:filling:low:bd} follows from
Poincar\'e--Lefschetz duality for compact orientable manifolds with
boundary; see eg Bredon \cite[page 360]{bredon}. In particular,
when $L$ is $3$--dimensional the special case $H^1(L,\Z)=0$ can
occur only if $\Sigma=S^2$.

\subsection[The Prescribed Boundary Problem in R8 and R10]{The prescribed boundary problem in $\R^8$ and $\R^{10}$}
\label{ss:pbp:n=4or5} We saw in \fullref{s:bdy}, that when $n$
equals $4$ or $5$, not every compact oriented $n$--manifold
$\Sigma$ admits an oriented filling $L$. Furthermore, even if
$\Sigma$ bounds it is not clear that the bounding manifold $L$
will have $TL^{\C}$ trivial or that exact immersions of $\Sigma$
in $\R^{2n}$ exist. Rather than try to analyze further conditions
on $\Sigma$ under which \fullref{assumptions} can be verified we shall
simply assume we are given exact initial data $(\Sigma, L,f)$.
Then we can apply \fullref{t:obstruction:summary} to obtain
the following results.
\begin{theorem}\label{t:filling:n=4,5}
Let $L$ be a connected oriented $n$--manifold with boundary
$\Sigma$ and let $(\Sigma,L,f)$ be exact initial data in the
sense of \fullref{d:initial:data}.
\begin{enumerate}
\item If $n=4$ and $\Sigma$ is connected then the Prescribed Boundary Problem
is solvable if and only if the following two conditions hold:
\begin{align*}\mu_f&\in \imag(i^*\co H^1(L,\Z)\ra H^1(\Sigma,\Z));\\ M_{SU}^*\Theta&=0.\end{align*}
\item In general, if $n=4$ or $n=5$ then the Prescribed Boundary Problem is solvable if and only if the following two conditions hold:
\begin{align*}\mu_f&\in \imag(i^*\co H^1(L,\Z)\ra H^1(\Sigma,\Z)); \\ M_{SU}^*\Theta&\in \imag(i^*\co H^3(L,\Z)\ra H^3(\Sigma,\Z)).\end{align*}
If the Prescribed Boundary Problem is solvable and moreover
$\mu_f=0$, then one can always find a solution with zero Maslov
class.
\end{enumerate}
\end{theorem}
\begin{proof}
The proof is similar to that of \fullref{t:filling:n=2}; when $\Sigma$ is connected and $n=4$ we can prove the stronger condition using \fullref{c:SU:obstruct}.
\end{proof}
Once again, it is clearly useful to know if $\Sigma$ admits simply connected fillings. Using
Surgery Theory one can prove the following simplification result for the fundamental group.
\begin{theorem}{\rm{\cite{milnor:killing}}}\qua
\label{T:kill:pi1}
Let $X$ be any oriented compact manifold (with or without boundary) of dimension
$n$ with $n\ge 4$. Then there exists a simply connected oriented manifold
$X'$ (obtained by a finite sequence of surgeries of type $(1,n-2)$)
with $\partial X = \partial X'$.
In particular, if $L$ is an oriented $n$--manifold which bounds $\Sigma$ and $n\ge 4$ then there is a simply connected
manifold $L'$ which bounds $\Sigma$.
\end{theorem}
The basic point of the proof of \fullref{T:kill:pi1} is that
any nontrivial loop in $\pi_1(X)$ can be represented by an
embedded $S^1$ which has a trivial normal bundle. Then by
performing surgery on these embedded loops one can systematically
kill the fundamental group of $X$. In order to try to simplify the
higher-dimensional homotopy or homology groups one needs to find
higher-dimensional embedded spheres with trivial normal bundles
representing elements of these groups. Hence with further
assumptions about the bounding manifold $L$ one can use Surgery
Theory to simplify its higher homotopy groups; see Kosinski \cite{kosinski} for results in this direction.

\section{Lagrangian desingularizations}\label{s:desing}
We now apply the results of \fullref{s:pbp:lowdim} to
answer the three questions posed in the Introduction to this paper.

\subsection[Lagrangian submanifolds with exact isolated
singularities in Cn]{Lagrangian submanifolds with exact isolated
singularities in $\C^n$} \label{ss:sing:cn}
\begin{definition}
\label{d:isolated:sing}
Let $X$ be a connected topological space with a finite number of
points $\{x_1,\dots,x_m\}$ such that the (not necessarily
connected) space $$X':=X\setminus \{x_1,\dots,x_m\}$$ is an oriented smooth $n$--manifold.
We say that $X$ is an \textit{oriented manifold with isolated
singularities} if, in addition, each $x_i$ has a connected open
neighbourhood $U_i$ such that $U_i\setminus \{x_i\}\simeq
(0,1)\times\Sigma_i$, for some compact (not necessarily connected)
($n-1$)-manifold $\Sigma_i$, and so that $U_i \cap U_j = \emptyset$, if $i \neq
j$.

Let $X$ be an oriented $n$--manifold with isolated singularities.
We say that $f\co X \ra \R^{2n}$ is an
\textit{oriented Lagrangian submanifold with isolated
singularities} if $f$ is continuous and $f|_{X'}$ is a Lagrangian
immersion.

We say that $X$ has \textit{exact singularities\/} if each restriction $f|_{U_i}$ is exact. This ``local exactness'' condition is satisfied automatically near any smooth point of $X$. We say that $X$ is exact or that it has zero Maslov class if this is true for $X'$.
\end{definition}
\begin{example}\label{E:lagr:cones:exact}
Let $X$ be any oriented Lagrangian cone in $\R^{2n}$ with an isolated
singularity at the origin, and let $\Sigma$ denote the
intersection of $X$ with $\Sph^{2n-1}$. $\Sigma$ is an oriented
Legendrian submanifold of $\Sph^{2n-1}$ which is not a totally
geodesic sphere.  The origin is then an isolated singular point of
$X$ which has a neighbourhood in $X$ whose boundary is
diffeomorphic to $\Sigma$; $X'$ has the same number of connected
components as $\Sigma$.

We now give two ways to see that any Lagrangian cone $X$ is exact.
Since $\Sigma$ is Legendrian and the standard contact form
$\alpha$ on $\Sph^{2n-1}$ is just the restriction of $\lambda$ to
$\Sph^{2n-1}$ it follows immediately that $\Sigma$ and hence $X$
is exact. Alternatively, given any closed curve $\gamma \subset X$
and $r\in\R^+$, we can ``slide" $\gamma$ to the curve
$r\cdot\gamma$. As $r\ra 0$ the length of $r\cdot \gamma$ goes to
$0$. Hence, by Stokes' Theorem,
$$ \int_\gamma\lambda=\int_{r\cdot\gamma}\lambda\rightarrow 0\quad \mbox{as }r\rightarrow 0.$$
Since $\int_\gamma\lambda$ is however independent of $r$, it must
be zero.
\end{example}
\begin{example}
Let $X$ be a smooth oriented $n$--manifold and $f\co X\rightarrow
\R^{2n}$ be a continuous map which is a Lagrangian immersion
except at a finite number of points $\{x_1, \ldots, x_m\}$. Then $f$ defines a Lagrangian submanifold with isolated singularities. In this example the singularities arise from $f$ rather than from $X$, so each $\Sigma_i$ is diffeomorphic to $\Sph^{n-1}$.
\end{example}
It is important to find conditions on $f$ ensuring that its
singularities are exact. This can be done fairly easily, as shown
by the following example.

\begin{example}\label{e:metric:exact}
Suppose $f\co X\ra\C^n$ is a Lagrangian submanifold with isolated
singularities $\{x_1,\dots,x_m\}$. Let $g$ denote the pullback
metric on $U_i\setminus \{x_i\}\simeq (0,1)\times\Sigma_i$. Let
$r$ denote the variable on $(0,1)$ and $p$ any point on
$\Sigma_i$.

Suppose that $g$ is a ``$O(r^\beta)$--approximation" of a metric
$g'=dr^2+r^{2\alpha} g_i$, for some $\alpha,\beta>0$ and some
smooth metric $g_i$ on the cross-section $\Sigma_i$. More
precisely, we ask that $g=g'+h$, for some symmetric tensor
$h=h(r,p)$ whose norm (calculated with respect to $g'$) is bounded
by some $r^\beta$. Then the singularities of $f$ are locally
exact. The proof is similar to that given for Lagrangian cones in
\fullref{E:lagr:cones:exact}, but we now use closed curves
$\gamma\subset\Sigma_i$. The condition on the metric ensures that
the length of the curves $r\cdot\gamma$ tends to zero as $r\ra 0$.
\end{example}
In particular, if $X$ is a SL $n$--fold with isolated conical
singularities in the sense of Joyce \cite[Definition
3.6]{joyce:conifolds1} then $X$ is an oriented Lagrangian
submanifold with exact isolated singularities in the sense of
\fullref{d:isolated:sing}.

\begin{definition}
\label{d:desing}
Let $f\co X \ra \R^{2n}$ be a Lagrangian submanifold
with isolated singularities $x_1, \ldots, x_m$, and mutually disjoint connected
neighbourhoods
$U_1, \ldots, U_m$ of these singular points as in
\fullref{d:isolated:sing}. Let $Y$ be a
smooth oriented $n$--manifold (not necessarily connected)
and let $V_1, \ldots, V_m$ be mutually disjoint
(not necessarily connected) open subsets of $Y$ so that
$Y\setminus \smash{\bigcup_{i=1}^{m}}
V_i$ is diffeomorphic  to $X\setminus \smash{\bigcup_{j=1}^m} U_j$. 
We call
 a Lagrangian immersion $\smash{\hat{f}\co Y\rightarrow \R^{2n}}$ an
\textit{oriented Lagrangian desingularization of $X$\/} if
$\hat{f}=f$ on the sets $Y\setminus \bigcup_{i=1}^m V_i=X\setminus
\bigcup_{j=1}^m U_j$.
\end{definition}

\begin{remark}
As a very simple example we can consider the connected $2$--dimensional singular space $X$ with one singular point $x_1$, formed by taking a $2$--sphere and identifying two distinct points.
Then there is a neighbourhood $U_1$ of $x_1$
that has two boundary components each diffeomorphic to $S^{1}$.
If we take $V_1$ to be the disjoint union of two
discs bounded by these two circles, then the resulting manifold $Y$
will be diffeomorphic to $S^2$. If instead we take $V_1$ to be a cylinder
bounding these two circles, the resulting manifold will be a $2$--torus.

More generally, even though we start with a connected space $X$, when we remove the singular points the resulting nonsingular manifold $X'$ will in general be disconnected.
The manifold $Y$ may take different components of $X'$ and
reconnect them into fewer components, but not necessarily into one component only. In particular, different components of the deleted neighbourhood of one singular point $U_i\setminus{\{x_i\}}$ can end up in
different connected components of $Y$.

\end{remark}

Given a Lagrangian submanifold $X$ with isolated singularities, it
is natural to ask if it admits any Lagrangian desingularization $Y$.
Furthermore, if desingularizations do exist then we would like to
know what topology $Y$ can have, and whether the desingularization
process preserves extra properties such as exactness or zero
Maslov class. In the oriented case we can now prove the following
results.
\begin{theorem}
\label{t:lag:desing} For $n=2$ or $n=3$, let $f\co X\ra \C^n$ be an
oriented Maslov-zero Lagrangian submanifold with exact isolated
singularities $\{x_1,\dots,x_m\}$ in the sense of \fullref{d:isolated:sing}. Then $X$ admits oriented Lagrangian
desingularizations in the sense of \fullref{d:desing} with
arbitrarily complicated topology.

Furthermore, let $U_i$ be a connected neighbourhood of $x_i$ with smooth (not necessarily connected) oriented boundary $\Sigma_i$. If either
\begin{itemize}
\item[(i)] each $\Sigma_i$ is connected, or
\item[(ii)] $X$ is SL,
\end{itemize}
then $X$ also admits Maslov-zero oriented desingularizations with
arbitrarily complicated topology. Finally, if each $\Sigma_i$ is
connected and $X$ is also exact, then all these desingularizations
are exact.
\end{theorem}
\begin{proof}
Recall from \fullref{d:isolated:sing} that $U_i\setminus
\{x_i\}\simeq (0,1)\times\Sigma_i$. Choose any oriented filling
$L_i$ of $\Sigma_i$. By assumption $X$ is locally exact, so we can apply the results of \fullref{s:pbp:lowdim} to $f$ restricted to
$\mathcal{O}\!p\,\Sigma_i\simeq \Sigma_i\times (1-\epsilon, 1]$;
this proves that the Prescribed Boundary Problem defined by the
initial data $(\Sigma_i,L_i,f)$ admits an exact Maslov-zero
solution $\hat{f}_i$. Iterating this construction near each
singular point produces an oriented Lagrangian desingularization
$\hat{f}\co Y\ra\C^n$ of $X$ in the sense of \fullref{d:isolated:sing}. $Y$ is smooth because we chose, as initial
data, an open neighbourhood of $\Sigma_i$ rather than $\Sigma_i$
itself. The topology of $Y$ is ``arbitrarily complicated" in the
sense already seen in \fullref{s:pbp:lowdim}: we are allowed
to perform any number of topological surgeries on each $L_i$ so as
to increase its topological complexity.

Suppose now that each $\Sigma_i$ is connected. Let $\gamma$ be a
closed curve in $Y$. Since each $\Sigma_i$ is connected, we can
homologically split $\gamma$ into the sum of $r$ closed curves
$\gamma_k$ such that each $\gamma_k$ is completely contained
either in $X\setminus \smash{\bigcup_{j=1}^m} U_j$ or in some $L_k$. Then,
since on each of these sets $\smash{\hat{f}}$ is Maslov-zero,
$$\mu_{\hat{f}}([\gamma])=\sum_{k=1}^r \,\mu_{\hat{f}}([\gamma_k])=0.$$
Hence $Y$ is Maslov-zero. The same method proves the corresponding
statement for exact desingularizations.

Suppose instead that $X$ is SL. The second proof of \fullref{c:SL:n=2} (which works also for $n=3$) shows that we can
extend the initial data $df$ on each $\mathcal{O}\!p\,\Sigma_i$ to a
SL monomorphism $\hat{F_i}$ defined on $L_i$. The corresponding
Lagrangian immersion of $L_i$, obtained via the $h$--principle, is
clearly Maslov-zero but has the additional virtue of being SL on
$\mathcal{O}\!p\,\Sigma_i$. Let us now choose any closed curve
$\gamma$ in $Y$. Let $G(\gamma)$ denote the corresponding curve of
tangent planes, defined as $G(\gamma)(t):=T_{\gamma(t)}Y$. This is
a closed curve in the Lagrangian Grassmannian bundle
$\textrm{Gr}_{\textrm{lag}}(\R^{2n})$. According to \fullref{l:homotopy:s1}, to show that $Y$ has zero Maslov class it is
equivalent to show that for any such $\gamma$ the corresponding
curve $\det_\C\circ G(\gamma)$ is homotopically trivial in
$\Sph^1$.

If $\gamma$ is completely contained in $X\setminus \bigcup_{j=1}^m
U_j$, then $G(\gamma)$ is by hypothesis a curve of SL planes so
$\det_{\C}\circ G(\gamma)\equiv 1$; in particular, it is
homotopically trivial.  The same is true if $\gamma$ is completely
contained in some $L_i$.
Finally, assume $\gamma$ enters some $L_i$ at a point
$p_i\in\Sigma_i$; it is then forced to exit at a point
$q_i\in\Sigma_i$. Consider the portion of $\gamma$ contained in
$L_i$. The $h$--principle gives us a homotopy between the
distribution of SL planes determined by the SL monomorphism
$\hat{F}_i$ and the tangent planes of $L_i$.  Restricting this
homotopy to the corresponding portion of $G(\gamma)$, shows that
our curve can be homotoped rel $\{p_i,q_i\}$ to a curve of SL
planes. We can repeat this procedure for each $L_i$; the remaining
portion of $G(\gamma)$ is already special Lagrangian. This proves
that also in this case $\det_{\C}\circ G(\gamma)$ is homotopically
trivial.
\end{proof}
\begin{remark}
Our desingularization method starts with a Maslov-zero $X$ and
replaces a neighbourhood of each singular point with a Maslov-zero
$L_i$. We should not expect \textit{any\/} such desingularization to
be Maslov-zero. After all, this procedure is purely local and any
Lagrangian submanifold is locally Maslov-zero, though not
necessarily globally so. The same is true for the exactness
condition.
\end{remark}
\begin{corollary} \label{c:SL:cones}
Let $X$ be a SL cone in $\C^3$. Then $X$ admits connected oriented exact Maslov-zero desingularizations of arbitrarily complicated topology.
\end{corollary}
\begin{proof} Recall from \fullref{E:lagr:cones:exact} or \fullref{e:metric:exact} that $X$ is exact. Let $\Sigma$ denote the (not necessarily connected)
link of $X$, ie the intersection of $X$ with the sphere
$\Sph^{5}\subset\C^3$. Choose any connected filling $L$ of
$\Sigma$; then the methods of \fullref{t:lag:desing} prove the
existence of a Maslov-zero Lagrangian desingularization $Y$. Since
any curve in $Y$ is homotopic to a curve in $L$ and $L$ is exact by construction, then $Y$ is exact
even if $\Sigma$ is not connected.
\end{proof}
\begin{example}
The result of \fullref{c:SL:cones} also holds for $n=2$, but
in this case it is well-known that SL desingularizations can be found explicitly as follows. Any $2$--dimensional minimal cone (and in particular any
SL cone) is the union of $2$--planes. The ``hyperk\"ahler rotation"
trick shows that a submanifold of $\R^4$ is SL if and only if it
is a complex submanifold with respect to one particular
(nonstandard) complex structure on $\R^4$. Accordingly, the
$i$--th plane $\pi_i$ can be described by some complex linear
equation $a_{i1}z_1+a_{i2}z_2=0$, so our cone will be given by the
equation $\prod_{i=1}^n (a_{i1}z_1+a_{i2}z_2)=0$. For any nonzero
$\epsilon\in\C$ the set $\prod_{i=1}^n
(a_{i1}z_1+a_{i2}z_2)=\epsilon$ describes a smooth complex surface
in $\R^4$. Notice that these surfaces do not coincide with the
original cone outside a compact set, but do become asymptotic to
the cone at infinity. In this sense, these
surfaces are SL desingularizations of the original SL cone.
\end{example}
\subsection{Desingularization of Lagrangian singularities in almost Calabi--Yau manifolds}
\label{s:general:ambient} The goal of this section is to extend
the results of \fullref{ss:sing:cn} to more general ambient
manifolds. Let $(M^{2n},\omega)$ be a smooth symplectic manifold.
By Darboux's Theorem $M$ is locally symplectomorphic to $\R^{2n}$;
notice that \fullref{d:isolated:sing} relies on purely
local properties of $f$ and $X$, so it is simple to extend this
definition so that it includes the notion of Lagrangian
submanifolds $f\co X\ra M$ with isolated singularities. Unless $M$ is
\textit{exact\/}, ie $\omega=d\lambda$ for some $1$--form
$\lambda$, it is not possible to define a global exactness
condition for $X$ (unless $X$ happens to be contained inside one
Darboux chart). However, it still makes sense to discuss whether
$X$ is locally exact near each of its singular points.

The generalization of the Maslov class to Lagrangian submanifolds
of general symplectic manifolds is more complicated but is
well-known. For completeness we describe this
generalization briefly and indicate the simplifications which
occur when the symplectic manifold is an \textit{almost Calabi--Yau
manifold}.

Let $L^n$ be an oriented Lagrangian submanifold of $(M,\omega)$. As in
\fullref{S:h:princ:closed} denote by
$\smash{\textrm{Gr}^+_{\textrm{lag}}}(M)$  the oriented Lagrangian
Grassmannian bundle of $M$. The main issue is that in general
$\textrm{Gr}^+_{\textrm{lag}}(M)$ is now a nontrivial
$\unitary{n}/\sorth{n}$ fibre bundle over $M$.

Let us fix a compatible almost complex structure $J$ on $M$ and let
$g$ be the metric defined by $\omega$ and $J$. Let
$\hat{\gamma}\co (D^2,\partial D)\ra (M,L)$ be a disk in $M$ with
boundary $\gamma:=\hat{\gamma}(\partial D) \subset L$. Since $D$
is contractible there exists a complex volume form $\Omega$ defined on a neighbourhood of $D$ in $M$.
Any choice of $\Omega$ defines an evaluation map
$$\Omega\co {\textrm{Gr}^+_{\textrm{lag}}(M)|_D}\rightarrow \C^*,$$
defined as follows.
Let $\{e_i\}_{i=1}^n$ be a $g$--orthonormal basis for $T_pL$.
Since $T_pL$ is Lagrangian $\{e_i\}_{i=1}^n$ is a basis for the
complex vector space $(T_p M,J)$. In particular, $\{e_i\}_{i=1}^n$
is linearly independent over $\C$ and hence defines an element
$e_1 \wedge \ldots \wedge e_n \neq 0$ in $\Lambda^{n,0}(T_p M,J)$.
Hence we obtain a map from $T_p L$ to $\C^*$ by evaluating
$\Omega$ on this nonzero $n$--vector. It is easy to check that
this map does not depend on the choice of the orthonormal
basis of $T_p L$.

The Gauss map associates to $\gamma$ a loop $G(\gamma)$ of
Lagrangian tangent spaces defined by $G(\gamma)(t):=T|_{\gamma(t)}L$; we can thus define the
\textit{Maslov index\/} $\mu_L(\hat{\gamma})$ of $\hat{\gamma}$ to
be the degree of the map $$\frac{\Omega}{|\Omega|}
\co G(\gamma)\subset\textrm{Gr}^+_{\textrm{lag}}(M)\rightarrow \Sph^1.$$
This number is independent of the choice of $J$ and
$\Omega$: this follows from the contractibility of $D$ and of the
space of compatible almost complex structures. For the same reason
it depends only on the relative homotopy class of $\hat{\gamma}$;
we have thus defined a map on the relative homotopy group
$$\mu_L\co \pi_2(M,L)\rightarrow \Z.$$
Let us now assume that $(M,\omega,g,J)$ is a K\"ahler manifold
with trivial canonical bundle $K_M$ and fix a global never-vanishing
holomorphic section $\Omega$ of $K_M$. In the special Lagrangian
literature the data $(M,\omega,g,J,\Omega)$ often goes under the
name \textit{almost Calabi--Yau manifold\/}; the manifold is
\textit{Calabi--Yau\/} if $\Omega$ is covariantly constant and
satisfies a normalization condition. Almost Calabi--Yau manifolds
provide a particularly convenient framework for discussing Maslov
indices, because the disk $D$ is no longer needed to guarantee the
existence of the complex volume form in a neighbourhood of a loop
$\gamma\subset L$.  Hence the Maslov index can be thought of as a
map   $$\mu_L\co \pi_1(L)\rightarrow \Z.$$ Alternatively, we could
proceed as in \fullref{s:maslov:class} using $\Omega/|\Omega|$
instead of $\det_\C$ to define a \textit{Maslov class\/} $\mu_L\in
H^1(L;\Z)$. Once again, homotopy considerations show that the
Maslov class is locally independent of the particular $\Omega,J$
used in this construction. In particular, suppose we are given a
Lagrangian submanifold of $M$ with isolated singularities. Then,
in a Darboux neighbourhood of a singularity, it is Maslov-zero
according to the definition just given if and only if it is
Maslov-zero with respect to $\det_\C$, ie according to
the definition of \fullref{s:maslov:class}.

We say that a Lagrangian submanifold $L$ of an almost
Calabi--Yau manifold is \textit{special Lagrangian\/} if
$$ \textrm{Im}\,\Omega|_{L}=0 \quad \textrm{and} \quad \textrm{Re}\,\Omega|_{L} >0.$$
Equivalently, notice that the data $(\omega,g,J,\Omega)$  defines
a $\textrm{SU}(n)$ principal fibre bundle over $M$, given by the
set of all $(g,J)$--unitary frames $\{e_1,\dots,e_n\}$ such that
$\Omega(e_1,\dots,e_n)\in \R^+$; $L$ is special Lagrangian if, for
each $p\in L$, we can find an orthonormal basis $\{e_1,\dots,
e_n\}$ of $T_pL$ which, as a unitary frame of $T_pM$, belongs to
this fibre bundle. Any SL submanifold is automatically oriented by
the volume form $\textrm{Re}(\Omega)|_{L}$, and clearly has zero
Maslov class.

The above considerations show that in Darboux coordinates any SL
submanifold $X$ with isolated singularities is locally Maslov-zero
in the sense of \fullref{s:maslov:class}. As seen in \fullref{e:metric:exact}, it is easy to use either the ambient metric
$g$ or the local metric induced by Darboux coordinates to find
conditions ensuring that $X$ is locally exact. Finally, the
\textit{$C^0$--dense\/} version of the $h$--principle \cite[\S
6.2.D]{eliashberg} proves that the methods used in the previous
Sections are completely local; in other words the
desingularizations that we build there, near any singular point
$x_i\in X$, can be made to live in any small neighbourhood of
$x_i$. In particular, the construction can take place completely
inside any Darboux coordinate chart. Our methods thus apply
verbatim, proving the following result.
\begin{corollary}
\label{c:SL:desing:CY} Let $M$ be an almost Calabi--Yau manifold of
dimension 2 or 3. Let $X\subset M$ be a SL submanifold with
isolated exact singularities. Then $X$ admits oriented Lagrangian
desingularizations with zero Maslov class and arbitrarily
complicated topology.
\end{corollary}
\subsection{Soft asymptotically conical smoothings of SL cones}\label{s:other:ends}
Finally, in this section we address the third question posed in
the Introduction to this paper. We begin with the following
definition, adapted from Joyce \cite[Definition 7.1]{joyce:conifolds1}.

\begin{definition}
Let $C^n$ be a regular oriented cone in $\C^n$ with oriented (not
necessarily connected) link $\Sigma^{n-1}$. Let $i:\Sigma\times
(0,\infty)\ra\C^n$ denote the corresponding immersion. We say that
an immersion $f\co \Sigma\times (1-\epsilon,\infty)\ra \C^n$ is an
\textit{asymptotically conical (AC) end with decay $\lambda$ and
cone $C$} if
$$|\nabla^k(f-i)|=O(r^{\lambda-1-k}) \ \ \mbox{as $r\rightarrow\infty$ for $k=0,1$}.$$
Here, $\nabla$ and $|\cdot|$ are defined using the natural metric
$g'$ on $C$ and we assume $\lambda<2$.
\end{definition}
Now suppose both $i$ and $f$ are Lagrangian. Even though $i$ is
automatically exact, in general there is no reason for $f$ to be
exact. However, it is simple to find additional assumptions
ensuring that $f$ is also exact. One such assumption was described
by Joyce \cite[Proposition 7.3]{joyce:conifolds1}: if $\lambda<0$ then
$f$ is exact. Notice also that, under this assumption, if $C$ is
Maslov-zero then this will be true also for $f$ because the two
manifolds are $C^1$--asymptotically close at infinity. Another
possible assumption is as follows. Suppose $\psi$ is a Hamiltonian
diffeomorphism of $\C^n$ and that $f:=\psi\circ i$ (where we
restrict $i$ to $\Sigma\times [1-\epsilon,\infty)$ and assume that
$\psi$ decays appropriately at infinity (so as to ensure the AC
condition on $f$). Then $f$ is exact and it is Maslov-zero if and
only if $C$ is.

Suppose now that $C$ is SL and that $\Sigma$ satisfies 
\fullref{assumptions}; as we have already
remarked, these assumptions can be satisfied in low dimensions via any
compact oriented choice of filling $L$. We can then try to solve
the Prescribed Boundary Problem for the exact initial data
$(\Sigma,L,f)$ where now $f$ is restricted to $\Sigma\times
(1-\epsilon,1]$. Since $C$ is Maslov-zero, under either of the
above assumptions $f$ will be Maslov-zero also. We can thus
conclude that in dimensions 2 and 3 the Prescribed Boundary
Problem is solvable. The resulting submanifold $Y$ still has the
AC ends given by $f$, and $Y$ is thus an exact Maslov-zero
Lagrangian submanifold asymptotic to the SL cone $C$.

\appendix
\section{Comparisons with other results in the literature}
\setobjecttype{App}
\label{appendix:compare} The goal of this Appendix is to fit our
results into a broader context by describing some similarities,
differences and relationships with other constructions in the
Lagrangian and SL literature.

\subsection{Lagrangian cobordism groups}
Arnold \cite{arnold} defined various types of Lagrangian and
Legendrian cobordism groups. In \cite{audin:cobord,audin:book},
Audin provides a detailed study of the oriented and nonoriented
\textit{exact Lagrangian cobordism groups\/}.

The abstract set-up is by now standard for cobordism-type
theories: given a specific category of geometric objects, one
defines an equivalence relationship such that, under disjoint
union, the equivalence classes inherit a group structure.
Following Thom, one then looks for an algebraic formulation of
these groups in terms of homotopy groups of certain universal
classifying spaces. The calculation of these homotopy groups is
often possible but relies on rather sophisticated methods in
homotopy theory.

In Audin's case, the geometric objects in question are pairs
$(L,f)$ where $L$ is a closed (not necessarily connected) smooth
manifold with $TL^\C$ trivial and $f\co L \ra \C^n$ is an exact
Lagrangian immersion. Two such pairs $(L_0,f_0)$ and $(L_1,f_1)$
are considered equivalent if
\begin{enumerate}
\item $L_0,L_1$ are cobordant in the usual sense; ie
there exists a manifold $Y$ such that $\partial Y=L_0\cup L_1$ (in
the oriented category, we also require that the orientations be
compatible in the sense that $\partial Y=-L_0\cup L_1$;
\item $Y$ has $TY^\C$ trivial and admits an exact Lagrangian immersion $F\co Y\ra\C^{n+1}$ which,
on the boundary, restricts to the isotropic immersions $(f_0,0,0)$
and $(f_1,1,0)$;
\item let $n$ denote the inward unit conormal of $L_i \subset Y$.
Then, using the standard notation for complex variables
$z_j=x_j+iy_j$, we require that $dF(n)=\partial x_{n+1}$ along
$L_0$ and $dF(n)=-\partial x_{n+1}$ along $L_1$.
\end{enumerate}
The third condition should be thought of as a transversality
condition which implies a canonical form for $F$ near the boundary
similar to our ``cylindrical thickenings" in \fullref{e:maslov:thickening}. A more thorough discussion of these
groups can be found in \cite[page 21]{audin:book}.

\cite{audin:book} provides an algebraic reformulation of these
Lagrangian cobordism groups (both oriented and nonoriented) in
terms of a Thom-type set-up. Using this reformulation Audin
computes the oriented cobordism groups in dimensions up to $10$,
often identifying explicit generators. Along the way, Audin
defines and studies the Maslov-zero analogues of these groups
(``SL cobordism groups") as a tool to compute the full oriented
cobordism groups \cite[page XII]{audin:book}. In the unoriented case
matters simplify considerably: if $L_0$ and $L_1$ are cobordant as
smooth unoriented manifolds and both admit Lagrangian immersions
$f_0, f_1$ in $\C^n$, then $(L_0,f_0)$ and $(L_1,f_1)$ are exact
Lagrangian cobordant \cite[Corollary 2.2]{audin:cobord}

Audin's results apply to the Lagrangian Cobordism Problem stated
in \fullref{ss:lagr:cobordism} as follows. Using the notation
of that section, suppose that the exact isotropic immersion $i:
\Sigma \ra \R^{2n}$ is in fact an exact Lagrangian immersion into
$\R^{2n-2}$. If the pair $(\Sigma,i)$ is zero in the exact
Lagrangian cobordism group then our Lagrangian Cobordism Problem
with boundary data $(\Sigma,i)$ is solvable. A priori, this does
not allow us to prescribe the topology of the filling $L$ nor the
data of the Lagrangian thickening. A deeper understanding of the
relationship between the Lagrangian cobordism groups and our
Prescribed Boundary Problem requires many of the tools described
in \fullref{S:h:princ:closed} and \fullref{s:bdy}.

We expect that the more direct methods used in this paper and the
explicit results of \fullref{s:pbp:lowdim} and \fullref{s:desing}
are closer to the needs of differential geometers interested in SL
singularities.

\subsection{Lagrangian surgeries} \label{ss:lagn:surgery} The
fillings and desingularizations obtained in \fullref{s:pbp:lowdim} and \fullref{s:desing} are in general only
immersed, not embedded. Moreover, as already mentioned in \fullref{r:embeddings}, Lagrangian embeddings of a fixed manifold $L$
do not satisfy any type of $h$--principle. However, a generic
Lagrangian immersion of $L$ will contain only isolated transverse
double points and there exists a standard Lagrangian surgery
procedure which replaces a neighborhood of any such
self-intersection with an embedded Lagrangian $1$--handle
\cite{polterovich}. If the number of double points is finite,
an iteration of this procedure will thus generate an embedded
Lagrangian submanifold $\tilde{L}$. Clearly, this procedure
changes the topology of $L$: $L$ and $\tilde{L}$ have different
cohomology groups in dimensions $1$ and $n-1$. It will thus affect
the Maslov class of $L$ and usually also its orientability.
However, when $n$ is odd-dimensional one can always choose this
surgery so that $\tilde{L}$ is orientable provided $L$ is. Thus,
for example, in the $n=3$ case the only problem introduced by
these surgeries concerns the Maslov class. However it is not clear
how to deal with this issue. In particular, it is not clear if it
is possible to apply the Lagrangian surgery procedure to the
immersed Maslov-zero desingularizations obtained in \fullref{s:desing} and get embedded Maslov-zero desingularizations. It
seems very likely that there would be further topological
obstructions to being able to do this.

Polterovich's surgery procedure relies on the existence of an
explicit Lagrangian model for desingularizing the union of any two
Lagrangian planes. It is precisely the lack of analogous models
for other types of singularities, such as most Lagrangian cones,
that motivates our use of nonexplicit $h$--principle techniques to
prove the existence of (immersed) desingularizations.

\subsection{Fu's SL moment conditions} It is well-known that in a Calabi--Yau manifold
$(M,\Omega,g,J,\omega)$, an
orientable submanifold $L$
is SL (with the correct choice of
orientation)  if and only if $\omega$ and $\beta:= \imag{\Omega}$
restrict to zero on $L$. Thus if a compact Lagrangian submanifold
$L^n$ (with or without boundary) is SL, or is even only homologous
to a SL, it must satisfy the condition $\int_L\beta=0$.

One can immediately derive other constraints that any SL
submanifold of a Calabi--Yau manifold must satisfy by considering
the \textit{differential ideal\/} $\mathcal{I}$ generated by
$\omega$ and $\beta$. Let $\Omega(M)$ denote the algebra of real differential forms on $M$. Recall that a differential ideal on a smooth
manifold $M$ is an ideal $\mathcal{I}$ of $\Omega(M)$ which is
also $d$--closed; in other words, $\mathcal{I}$ is a differential
ideal if $\alpha \in \mathcal{I}$ and $\beta \in \Omega(M)$
implies that $\alpha \wedge \beta \in \mathcal{I}$ and $d\alpha
\in \mathcal{I}$. Using the standard properties of forms under
pullback and exterior differentiation we see that if $f$ is a
SL-immersion, then $f^*\sigma=0$ for any $\sigma \in
\mathcal{I}(\omega,\beta)$, the differential ideal generated by
$\omega$ and $\beta$. In particular, if $L$ is SL then
$\int_L\alpha=0$ for \textit{any\/} $\alpha\in\mathcal{I}^n$, where
$\mathcal{I}^n$ is the degree $n$ part of $\mathcal{I}$.

Fu \cite{fu} used these constraints to determine necessary
conditions for a compact orientable $(n-1)$--dimensional isotropic
submanifold $\Sigma \subset \C^n$ to bound a SL $n$--fold. In many
ways this is the SL analogue of the Lagrangian Cobordism Problem,
so we now briefly summarize Fu's construction. For notational
simplicity we will write $\Omega^k$ for $\Omega^k(M)$.

Consider the following complex:
$$(\Omega^0,\mathcal{I}^0) \stackrel{d}{\ra} \dots\ra(\Omega^{k-1},\mathcal{I}^{k-1})\stackrel{d}{\ra}(\Omega^k,\mathcal{I}^k)\ra\dots \stackrel{d}{\ra}(\Omega^m,\mathcal{I}^m).$$
Let $\mathcal{H}^k$ denote the corresponding (relative) cohomology
groups. That is, we define
$$\mathcal{Z}^k:=\{\alpha\in\Omega^k: d\alpha\in\mathcal{I}^{k+1}\} \quad\text{and}\quad \mathcal{B}^k:=d\Omega^{k-1}+\mathcal{I}^{k},$$
and then $\mathcal{H}^k$ is defined as the quotient space
$\mathcal{H}^k:=\mathcal{Z}^k/\mathcal{B}^k$. Using the facts that
$\Sigma$ is isotropic and that elements $\gamma$ of
$\mathcal{I}^{n-1}$ are of the form $\gamma=\omega\wedge\gamma'$,
it follows that $\int_\Sigma[\alpha]$ is well-defined for all
$[\alpha]\in \mathcal{H}^{n-1}$. If $\Sigma$ admits a smooth SL
filling, then $\Sigma$ satisfies the \textit{SL moment conditions\/}
$$\int_\Sigma[\alpha]=0, \quad \textrm{for all} \quad [\alpha]\in \mathcal{H}^{n-1}.$$
In fact, even if $\Sigma$ is only the boundary of a (possibly
singular) $n$--dimensional SL-rectifiable current then $\Sigma$
still satisfies these moments conditions. Furthermore,  if
$\Sigma$ admits a (possibly singular) SL filling, then
\textit{any\/} Lagrangian filling $L$ of $\Sigma$ must satisfy the
condition
$$\int_Ld[\alpha]=0, \quad \textrm{for any} \quad  [\alpha]\in
\mathcal{H}^{n-1}.$$ Hence the SL moment conditions are not
sensitive enough to detect the difference between the existence of
a smooth SL filling and a singular one. For example, the SL moment
conditions will be satisfied even if $\Sigma$ can only be filled
by a singular SL cone. In this sense Fu's moment conditions are
strictly weaker than what one would need to solve the ``smooth SL
Cobordism Problem". Furthermore, using an explicit description of the SL moment
conditions in $\C^n$ (see below), Fu showed that there are isotropic manifolds
$\Sigma$ which satisfy all the SL moment conditions and yet do not
bound any SL rectifiable current. This shows that the SL moment conditions are not sufficient for the ``singular SL Cobordism Problem" either.

When $M=\C^n$, Fu was able to express the SL moment conditions
very explicitly. In particular he showed that the space
$\mathcal{H}^{n-1}$ is isomorphic to the space of functions $f$
satisfying the condition $\mathcal{L}_{X_f}\beta=0$, where
$\mathcal{L}$ denotes Lie differentiation and $X_f$ is the
Hamiltonian vector field generated by $f$. He then determined all
such functions explicitly;
for $n\ge 3$ they are exactly the hermitian harmonic quadratics
defined following \fullref{harmonic:momentmap}
and arising from the Hamiltonian
action of $\sunitary{n} \ltimes \C^n$ on $\C^n$. Hence the SL
moment conditions amount to exactly $n(n+2)$ independent
conditions on any compact orientable isotropic $(n-1)$--manifold
$\Sigma \subset \C^n$.

\subsection{Joyce's SL gluing results} In the papers
\cite{joyce:conifolds3,joyce:conifolds4} Joyce presents a series
of results on the desingularization of certain kinds of singular
SL submanifolds, which we have called SL conifolds.
More specifically, these papers deal with a
compact SL $n$--fold  $X$ in $M$ with isolated singularities
$x_1,\dots,x_m$ (as in our \fullref{d:isolated:sing}),
satisfying two additional constraints. First, he assumes that each
singularity $x_i$ is modelled metrically (in a precise sense) on a
regular SL cone $C_i\subset \C^n$. Second, he assumes that there
exists a smooth SL submanifold $L_i\subset\C^n$ which is
asymptotic (again, in a precise sense) to the cone $C_i$. He then
proves that, under certain conditions, it is possible to remove a
neighbourhood of each singular point $x_i$ in $X$ and glue in a
copy of $L_i$ so that, after a global perturbation, the resulting
submanifold $Y$ is SL in $M$.

The simplest case in which Joyce obtains desingularization results
is presented in \cite[Theorem 6.13]{joyce:conifolds3}. The first step
in his process is to produce a Lagrangian gluing of each $L_i$ to
$X$ as follows. Let $\Sigma_i$ denote the link of the SL cone
$C_i$. The use of Lagrangian neighbourhoods allows him to work
inside the cotangent bundle $T^*(\Sigma_i\times(0,\infty))$,
endowed with its standard symplectic structure. More precisely,
the gluing happens inside a small region of this bundle,
$T^*(\Sigma_i\times[t^\tau,2t^\tau])$ (depending on parameters
$t,\tau$). There, the relevant regions of $L_i$ and $X$ appear as
the graphs of closed $1$--forms $\chi_i$ and $\eta_i$ respectively.
Producing the required Lagrangian gluing reduces to proving that
there exists a closed $1$--form which interpolates between $\chi_i\in
T^*(\Sigma_i\times \{t^\tau\})$ and $\eta_i\in T^*(\Sigma_i\times
\{2t^\tau\})$. There is an obvious obstruction to doing this: the
cone-like form of the metric on $X$ close to each singular point
$x_i$ implies (as in our \fullref{e:metric:exact}) that
$\eta_i$ is exact, so such an interpolation exists if and only if
$\chi_i$ is also exact. The exactness of $\chi_i$ is an additional
condition on $L_i$ which is encoded in Joyce's cohomological
invariant $Y(L_i)$ (see \cite[Definition 7.2]{joyce:conifolds1}). Thus
the Lagrangian gluing procedure just described is possible if and
only if $Y(L_i)=0$.

There are some strong similarities between the Lagrangian gluing
procedure outlined above and the methods we use, eg in
\fullref{c:SL:desing:CY}. First of all, both are completely
local in the sense that they perturb the original $X$ only in a
neighborhood of the singularity, leaving the rest untouched.
Second, in both cases (and for the same reasons) the singularities
are locally exact. The main difference is of course that Joyce
relies on an \textit{a priori\/} smoothing $L_i$ of $C_i$ while we
create our own smoothings using $h$--principle techniques. Thus in
our case the role played by exactness is different from in Joyce's
work; in our setting it is a necessary condition for applying the
$h$--principle.

There are SL smoothings $L$ of SL cones $C$ which do not satisfy
$Y(L)=0$ \cite[Example 6.8]{joyce:conifolds:summary}. For
desingularizations based on these models, Joyce provides a second
type of Lagrangian gluing process which perturbs $X$ globally.
Again, the gluing is possible only if a certain cohomological
compatibility holds between the appropriate regions of $L_i$ and
$X$.

It may also be useful to point out that the reason why Joyce is
not concerned with Maslov-class type problems is that, by being
careful in the Lagrangian gluing process, the Lagrangian
submanifolds he produces this way are almost SL (in a precise
sense). Thus they are automatically Maslov-zero. In conclusion, it
is the whole set-up Joyce begins with that allows him not to worry
about the kind of topological questions we address here, and
instead concentrate on the analytic aspects of the SL condition.

\begin{remark}
Joyce also has results for a SL analogue of the Lagrangian surgery
procedure explained in this Appendix, \fullref{ss:lagn:surgery}. Once again, the starting point for this is
provided by explicit local models due to Lawlor. These models
furnish local SL desingularizations of certain pairs of SL planes,
satisfying an angle condition. Again one has no Maslov class
problems because one can create an almost SL smoothing which
therefore has zero Maslov class. Joyce shows that any connected SL
submanifold with transverse self-intersection points modelled by
such pairs of planes can be desingularized, leading to an embedded
SL submanifold. The angle condition is always verified in the
$n=3$ case. Joyce also provides results for a ``SL connect sum" of
smooth SL submanifolds which intersect one another transversely.
Related results were obtained by A\,Butscher, Y\,I\,Lee and D\,Lee
\cite{butscher,y:lee,d:lee}.
\end{remark}

\bibliographystyle{gtart}
\bibliography{link}

\begin{thebibliography}{}
\providecommand\bibmarginpar{\leavevmode\marginpar}
\def\urlstyle#1{{\tt #1}}

\bibitem{arnold}
\textbf{V\,I Arnold}, \emph{Lagrange and {L}egendre cobordisms I}, Funktsional.
  Anal. i Prilozhen. 14 (1980) 1--13, 96 \xox{MR}{583797}

\bibitem{auckly}
\textbf{D Auckly}, \textbf{L Kapitanski},
  \href{http://dx.doi.org/10.1007/s00220-003-0901--x} {\emph{Holonomy and
  {S}kyrme's model}}, Comm. Math. Phys. 240 (2003) 97--122 \xox{MR}{2004981}

\bibitem{audin:cobord}
\textbf{M Audin}, \href{http://www.numdam.org/item?id=AIF_1985__35_3_159_0}
  {\emph{Quelques calculs en cobordisme lagrangien}}, Ann. Inst. Fourier
  $($Grenoble$)$ 35 (1985) 159--194 \xox{MR}{810672}

\bibitem{audin:book}
\textbf{M Audin}, \emph{Cobordismes d'immersions lagrangiennes et
  legendriennes}, Travaux en Cours 20, Hermann, Paris (1987) \xox{MR}{903652}

\bibitem{audin}
\textbf{M Audin}, \textbf{J Lafontaine} (editors), \emph{Holomorphic curves in
  symplectic geometry}, Progress in Mathematics 117, Birkh\"auser Verlag, Basel
  (1994) \xox{MR}{1274923}

\bibitem{bott:tu}
\textbf{R Bott}, \textbf{L\,W Tu}, \emph{Differential forms in algebraic
  topology}, Graduate Texts in Mathematics 82, Springer, New York (1982)
  \xox{MR}{658304}

\bibitem{bredon}
\textbf{G\,E Bredon}, \emph{Topology and geometry}, Graduate Texts in
  Mathematics 139, Springer, New York (1993) \xox{MR}{1224675}

\bibitem{butscher}
\textbf{A Butscher}, \emph{Regularizing a singular special {L}agrangian
  variety}, Comm. Anal. Geom. 12 (2004) 733--791 \xox{MR}{2104075}

\bibitem{mcintosh:carberry}
\textbf{E Carberry}, \textbf{I McIntosh}, \emph{Minimal {L}agrangian 2--tori in
  $\Bbb C\Bbb P^2$ come in real families of every dimension}, J. London Math.
  Soc. $(2)$ 69 (2004) 531--544 \xox{MR}{2040620}

\bibitem{castro:urbano:sym}
\textbf{I Castro}, \textbf{F Urbano}, \emph{On a new construction of special
  {L}agrangian immersions in complex {E}uclidean space}, Q. J. Math. 55 (2004)
  253--265 \xox{MR}{2082092}

\bibitem{cheng}
\textbf{B\,N Cheng}, \emph{Area-minimizing cone-type surfaces and coflat
  calibrations}, Indiana Univ. Math. J. 37 (1988) 505--535 \xox{MR}{962922}

\bibitem{don:kron}
\textbf{S\,K Donaldson}, \textbf{P\,B Kronheimer}, \emph{The geometry of
  four-manifolds}, Oxford Mathematical Monographs, The Clarendon Press Oxford
  University Press, New York (1990) \xox{MR}{1079726}, Oxford Science
  Publications

\bibitem{eliashberg}
\textbf{Y Eliashberg}, \textbf{N Mishachev}, \emph{Introduction to the
  $h$--principle}, Graduate Studies in Mathematics 48, American Mathematical
  Society, Providence, RI (2002) \xox{MR}{1909245}

\bibitem{fu}
\textbf{L Fu},
  \href{http://projecteuclid.org/getRecord?id=euclid.dmj/1077285157} {\emph{On
  the boundaries of special {L}agrangian submanifolds}}, Duke Math. J. 79
  (1995) 405--422 \xox{MR}{1344766}

\bibitem{geiges}
\textbf{H Geiges}, \emph{Contact geometry}, from: ``Handbook of differential
  geometry. Vol. II'', Elsevier/North-Holland, Amsterdam (2006)  315--382
  \xox{MR}{2194671}

\bibitem{gromov:jholo}
\textbf{M Gromov}, \href{http://dx.doi.org/10.1007/BF01388806}
  {\emph{Pseudoholomorphic curves in symplectic manifolds}}, Invent. Math. 82
  (1985) 307--347 \xox{MR}{809718}

\bibitem{gromov:pdr}
\textbf{M Gromov}, \emph{Partial differential relations}, volume~9 of
  \emph{Ergebnisse der Mathematik und ihrer Grenzgebiete (3) [Results in
  Mathematics and Related Areas (3)]}, Springer, Berlin (1986) \xox{MR}{864505}

\bibitem{gross:fibrations1}
\textbf{M Gross}, \emph{Special {L}agrangian fibrations I: {T}opology}, from:
  ``Integrable systems and algebraic geometry (Kobe/Kyoto, 1997)'', World Sci.
  Publ., River Edge, NJ (1998)  156--193 \xox{MR}{1672120}

\bibitem{gross:fibrations2}
\textbf{M Gross}, \emph{Special {L}agrangian fibrations II: {G}eometry. {A}
  survey of techniques in the study of special {L}agrangian fibrations}, from:
  ``Surveys in differential geometry: differential geometry inspired by string
  theory'', Surv. Differ. Geom. 5, Int. Press, Boston (1999)  341--403
  \xox{MR}{1772274}

\bibitem{gross:examples}
\textbf{M Gross}, \emph{Examples of special {L}agrangian fibrations}, from:
  ``Symplectic geometry and mirror symmetry (Seoul, 2000)'', World Sci. Publ.,
  River Edge, NJ (2001)  81--109 \xox{MR}{1882328}

\bibitem{hardt:simon}
\textbf{R Hardt}, \textbf{L Simon}, \emph{Area minimizing hypersurfaces with
  isolated singularities}, J. Reine Angew. Math. 362 (1985) 102--129
  \xox{MR}{809969}

\bibitem{harvey:lawson}
\textbf{R Harvey}, \textbf{H\,B Lawson, Jr}, \emph{Calibrated geometries}, Acta
  Math. 148 (1982) 47--157 \xox{MR}{666108}

\bibitem{haskins:invent}
\textbf{M Haskins}, \emph{The geometric complexity of special {L}agrangian {$T^
  2$}-cones}, Invent. Math. 157 (2004) 11--70 \xox{MR}{2135184}

\bibitem{haskins:ajm}
\textbf{M Haskins}, \emph{Special {L}agrangian cones}, Amer. J. Math. 126
  (2004) 845--871 \xox{MR}{2075484}

\bibitem{haskins:kapouleas}
\textbf{M Haskins}, \textbf{N Kapouleas}, \emph{{S}pecial {L}agrangian cones
  with higher genus links}, to appear in Inventiones Mathematicae
  \xox{arXiv}{math.DG/0512178}

\bibitem{hatcher}
\textbf{A Hatcher}, \emph{Algebraic topology}, Cambridge University Press,
  Cambridge (2002) \xox{MR}{1867354}

\bibitem{hirsch}
\textbf{M\,W Hirsch},
  \href{http://links.jstor.org/sici?sici=0002-9947(195911)93:2%3C242:IOM%3E2.0%
.CO%3B2--Y} {\emph{Immersions of manifolds}}, Trans. Amer. Math. Soc. 93 (1959)
  242--276 \xox{MR}{0119214}

\bibitem{joyce:syz}
\textbf{D Joyce}, \emph{Singularities of special {L}agrangian fibrations and
  the {SYZ} conjecture}, Comm. Anal. Geom. 11 (2003) 859--907 \xox{MR}{2032503}

\bibitem{joyce:conifolds5}
\textbf{D Joyce}, \emph{Special {L}agrangian submanifolds with isolated conical
  singularities V: {S}urvey and applications}, J. Differential Geom. 63 (2003)
  279--347 \xox{MR}{2015549}

\bibitem{joyce:conifolds:summary}
\textbf{D Joyce}, \emph{Singularities of special {L}agrangian submanifolds},
  from: ``Different faces of geometry'', Int. Math. Ser. (N. Y.) 3,
  Kluwer/Plenum, New York (2004)  163--198 \xox{MR}{2102996}

\bibitem{joyce:conifolds1}
\textbf{D Joyce}, \href{http://dx.doi.org/10.1023/B:AGAG.0000023229.72953.57}
  {\emph{Special {L}agrangian submanifolds with isolated conical singularities
  I: {R}egularity}}, Ann. Global Anal. Geom. 25 (2004) 201--251
  \xox{MR}{2053761}

\bibitem{joyce:conifolds2}
\textbf{D Joyce}, \href{http://dx.doi.org/10.1023/B:AGAG.0000023230.21785.8d}
  {\emph{Special {L}agrangian submanifolds with isolated conical singularities
  II: {M}oduli spaces}}, Ann. Global Anal. Geom. 25 (2004) 301--352
  \xox{MR}{2054572}

\bibitem{joyce:conifolds3}
\textbf{D Joyce}, \href{http://dx.doi.org/10.1023/B:AGAG.0000023231.31950.cc}
  {\emph{Special {L}agrangian submanifolds with isolated conical singularities
  III: {D}esingularization, the unobstructed case}}, Ann. Global Anal. Geom. 26
  (2004) 1--58 \xox{MR}{2054578}

\bibitem{joyce:conifolds4}
\textbf{D Joyce}, \emph{Special {L}agrangian submanifolds with isolated conical
  singularities IV: {D}esingularization, obstructions and families}, Ann.
  Global Anal. Geom. 26 (2004) 117--174 \xox{MR}{2070685}

\bibitem{kosinski}
\textbf{A\,A Kosinski}, \emph{Differential manifolds}, Pure and Applied
  Mathematics 138, Academic Press, Boston (1993) \xox{MR}{1190010}

\bibitem{d:lee}
\textbf{D\,A Lee},
  \href{http://projecteuclid.org/getRecord?id=euclid.cag/1090526521}
  {\emph{Connected sums of special {L}agrangian submanifolds}}, Comm. Anal.
  Geom. 12 (2004) 553--579 \xox{MR}{2128603}

\bibitem{y:lee}
\textbf{Y-I Lee}, \emph{Embedded special {L}agrangian submanifolds in
  {C}alabi--{Y}au manifolds}, Comm. Anal. Geom. 11 (2003) 391--423
  \xox{MR}{2015752}

\bibitem{lees}
\textbf{J\,A Lees}, \emph{On the classification of {L}agrange immersions}, Duke
  Math. J. 43 (1976) 217--224 \xox{MR}{0410764}

\bibitem{mcintosh}
\textbf{I McIntosh}, \emph{Special {L}agrangian cones in $\Bbb C^3$ and
  primitive harmonic maps}, J. London Math. Soc. $(2)$ 67 (2003) 769--789
  \xox{MR}{1967705}

\bibitem{milnor:killing}
\textbf{J\,W Milnor}, \emph{A procedure for killing homotopy groups of
  differentiable manifolds.}, from: ``Proc. Sympos. Pure Math., Vol. III'',
  American Mathematical Society, Providence, R.I (1961)  39--55
  \xox{MR}{0130696}

\bibitem{milnor:stasheff}
\textbf{J\,W Milnor}, \textbf{J\,D Stasheff}, \emph{Characteristic classes},
  Annals of Mathematics Studies 76, Princeton University Press, Princeton, NJ
  (1974) \xox{MR}{0440554}

\bibitem{naitoh}
\textbf{H Naitoh}, \emph{Totally real parallel submanifolds in $P^{n}(c)$},
  Tokyo J. Math. 4 (1981) 279--306 \xox{MR}{646040}

\bibitem{ohnita:slg}
\textbf{Y Ohnita}, \emph{Stability and rigidity of special {L}agrangian cones
  over certain minimal {L}egendrian orbits}, to appear in Osaka J Math
  (2006)~Osaka City Preprint Series 2006, 06-2

\bibitem{polterovich}
\textbf{L Polterovich}, \href{http://dx.doi.org/10.1007/BF01896378} {\emph{The
  surgery of {L}agrange submanifolds}}, Geom. Funct. Anal. 1 (1991) 198--210
  \xox{MR}{1097259}

\bibitem{schoen:wolfson:volume}
\textbf{R Schoen}, \textbf{J Wolfson}, \emph{The volume functional for
  {L}agrangian submanifolds}, from: ``Lectures on partial differential
  equations'', New Stud. Adv. Math. 2, Int. Press, Somerville, MA (2003)
  181--191 \xox{MR}{2055848}

\bibitem{smale}
\textbf{S Smale}, \emph{The classification of immersions of spheres in
  {E}uclidean spaces}, Ann. of Math. $(2)$ 69 (1959) 327--344 \xox{MR}{0105117}

\bibitem{steenrod}
\textbf{N Steenrod}, \emph{The topology of fibre bundles}, Princeton Landmarks
  in Mathematics, Princeton University Press, Princeton, NJ (1999)
  \xox{MR}{1688579}Reprint of the 1957 edition, Princeton Paperbacks

\bibitem{syz}
\textbf{A Strominger}, \textbf{S-T Yau}, \textbf{E Zaslow}, \emph{Mirror
  symmetry is $T$--duality}, Nuclear Phys. B 479 (1996) 243--259
  \xox{MR}{1429831}

\bibitem{thom}
\textbf{R Thom}, \emph{Quelques propri\'et\'es globales des vari\'et\'es
  diff\'erentiables}, Comment. Math. Helv. 28 (1954) 17--86 \xox{MR}{0061823}

\bibitem{wall}
\textbf{C\,T\,C Wall}, \emph{Determination of the cobordism ring}, Ann. of
  Math. $(2)$ 72 (1960) 292--311 \xox{MR}{0120654}

\bibitem{whitehead:1940}
\textbf{J\,H\,C Whitehead},
  \href{http://links.jstor.org/sici?sici=0003-486X(194010)2:41:4%3C809:OC%3E2.%
0.CO%3B2-2} {\emph{On $C^1$--complexes}}, Ann. of Math. $(2)$ 41 (1940)
  809--824 \xox{MR}{0002545}

\bibitem{whitney:hard}
\textbf{H Whitney},
  \href{http://links.jstor.org/sici?sici=0003-486X(194404)2:45:2%3C220:TSOASN%%
3E2.0.CO%3B2-4} {\emph{The self-intersections of a smooth $n$--manifold in
  $2n$--space}}, Ann. of Math. $(2)$ 45 (1944) 220--246 \xox{MR}{0010274}

\end{thebibliography}

\end{document}